\documentclass[american,onefignum,onetabnum]{article}
\usepackage[T1]{fontenc}
\usepackage[latin9]{inputenc}
\usepackage[english]{babel}

\usepackage{amsfonts}
\usepackage{graphicx}
\usepackage{epstopdf}
\usepackage{algorithm}
\usepackage{algorithmic}
\ifpdf
  \DeclareGraphicsExtensions{.eps,.pdf,.png,.jpg}
\else
  \DeclareGraphicsExtensions{.eps}
\fi

\usepackage{array}
\usepackage{float}
\usepackage{units}
\usepackage{mathtools}
\usepackage{bm}
\usepackage{mathtools}
\usepackage{bm}
\usepackage{amsmath}
\usepackage{amsthm}
\usepackage{stmaryrd}
\usepackage{wasysym}
\usepackage{bbm}
\usepackage{url}

\newcommand{\noun}[1]{\textsc{#1}}

\providecommand{\tabularnewline}{\\}
\floatstyle{ruled}
\newfloat{algorithm}{tbp}{loa}
\providecommand{\algorithmname}{Algorithm}
\floatname{algorithm}{\protect\algorithmname}

\numberwithin{equation}{section}
\numberwithin{figure}{section}
\theoremstyle{plain}
\newtheorem{theorem}{\protect\theoremname}
\theoremstyle{definition}
\newtheorem{definition}[theorem]{\protect\definitionname}
\theoremstyle{plain}

\theoremstyle{remark}
\newtheorem{remark}[theorem]{\protect\remarkname}
\theoremstyle{plain}
\newtheorem{lemma}[theorem]{\protect\lemmaname}
\newtheorem{corollary}[theorem]{\protect\corollaryname}

\usepackage{amsopn}
\DeclareMathOperator{\diag}{diag}

\makeatother

\usepackage{babel}
\addto\captionsamerican{\renewcommand{\lemmaname}{lemma}}
\addto\captionsamerican{\renewcommand{\assumptionname}{Assumption}}
\addto\captionsamerican{\renewcommand{\definitionname}{definition}}
\addto\captionsamerican{\renewcommand{\remarkname}{remark}}
\addto\captionsamerican{\renewcommand{\theoremname}{Theorem}}
\addto\captionsenglish{\renewcommand{\lemmaname}{lemma}}
\addto\captionsenglish{\renewcommand{\algorithmname}{Algorithm}}
\addto\captionsenglish{\renewcommand{\assumptionname}{Assumption}}
\addto\captionsenglish{\renewcommand{\definitionname}{definition}}
\addto\captionsenglish{\renewcommand{\remarkname}{remark}}
\addto\captionsenglish{\renewcommand{\theoremname}{Theorem}}
\providecommand{\lemmaname}{lemma}
\providecommand{\assumptionname}{Assumption}
\providecommand{\definitionname}{definition}
\providecommand{\remarkname}{remark}
\providecommand{\theoremname}{Theorem}
\providecommand{\corollaryname}{Corollary}

\begin{document}

\selectlanguage{english}%
\global\long\def\R{\mathbb{R}}%

\global\long\def\e{{\mathbf{e}}}%

\global\long\def\et#1{{\e(#1)}}%

\global\long\def\ef{{\mathbf{\et{\cdot}}}}%

\global\long\def\x{{\mathbf{x}}}%

\global\long\def\w{{\mathbf{w}}}%

\global\long\def\m{{\mathbf{m}}}%

\global\long\def\xt#1{{\x(#1)}}%

\global\long\def\xf{{\mathbf{\xt{\cdot}}}}%

\global\long\def\d{{\mathbf{d}}}%

\global\long\def\b{{\mathbf{b}}}%

\global\long\def\k{{\mathbf{k}}}%

\global\long\def\a{{\mathbf{a}}}%

\global\long\def\u{{\mathbf{u}}}%

\global\long\def\y{{\mathbf{y}}}%

\global\long\def\yt#1{{\y(#1)}}%

\global\long\def\yf{{\mathbf{\yt{\cdot}}}}%

\global\long\def\z{{\mathbf{z}}}%

\global\long\def\w{{\mathbf{w}}}%

\global\long\def\v{{\mathbf{v}}}%

\global\long\def\h{{\mathbf{h}}}%

\global\long\def\s{{\mathbf{s}}}%

\global\long\def\c{{\mathbf{c}}}%

\global\long\def\p{{\mathbf{p}}}%

\global\long\def\f{{\mathbf{f}}}%

\global\long\def\rb{{\mathbf{r}}}%

\global\long\def\rt#1{{\rb(#1)}}%

\global\long\def\rf{{\mathbf{\rt{\cdot}}}}%

\global\long\def\mat#1{{\ensuremath{\bm{\mathrm{#1}}}}}%

\global\long\def\matN{\ensuremath{{\bm{\mathrm{N}}}}}%

\global\long\def\matA{\ensuremath{{\bm{\mathrm{A}}}}}%

\global\long\def\matB{\ensuremath{{\bm{\mathrm{B}}}}}%

\global\long\def\matC{\ensuremath{{\bm{\mathrm{C}}}}}%

\global\long\def\matD{\ensuremath{{\bm{\mathrm{D}}}}}%

\global\long\def\matP{\ensuremath{{\bm{\mathrm{P}}}}}%

\global\long\def\matU{\ensuremath{{\bm{\mathrm{U}}}}}%

\global\long\def\matV{\ensuremath{{\bm{\mathrm{V}}}}}%

\global\long\def\matW{\ensuremath{{\bm{\mathrm{W}}}}}%

\global\long\def\matM{\ensuremath{{\bm{\mathrm{M}}}}}%

\global\long\def\matZ{\ensuremath{{\bm{\mathrm{Z}}}}}%

\global\long\def\matR{\mat R}%

\global\long\def\matQ{\mat Q}%

\global\long\def\matS{\mat S}%

\global\long\def\matSb{\mat{S_{B}}}%

\global\long\def\matSw{\mat{S_{w}}}%

\global\long\def\matY{\mat Y}%

\global\long\def\matX{\mat X}%

\global\long\def\matYhat{\hat{\mat Y}}%

\global\long\def\matXhat{\hat{\mat X}}%

\global\long\def\matI{\mat I}%

\global\long\def\matzero{\mat 0}%

\global\long\def\matJ{\mat J}%

\global\long\def\matZ{\mat Z}%

\global\long\def\matL{\mat L}%

\global\long\def\S#1{{\mathbb{S}_{N}[#1]}}%

\global\long\def\IS#1{{\mathbb{S}_{N}^{-1}[#1]}}%

\global\long\def\PN{\mathbb{P}_{N}}%

\global\long\def\TNormS#1{\|#1\|_{2}^{2}}%

\global\long\def\TNorm#1{\|#1\|_{2}}%

\global\long\def\InfNorm#1{\|#1\|_{\infty}}%

\global\long\def\FNorm#1{\|#1\|_{F}}%

\global\long\def\FNormS#1{\|#1\|_{F}^{2}}%

\global\long\def\UNorm#1{\|#1\|_{\matU}}%

\global\long\def\UNormS#1{\|#1\|_{\matU}^{2}}%

\global\long\def\UINormS#1{\|#1\|_{\matU^{-1}}^{2}}%

\global\long\def\ANorm#1{\|#1\|_{\matA}}%

\global\long\def\BNorm#1{\|#1\|_{\mat B}}%

\global\long\def\ANormS#1{\|#1\|_{\matA}^{2}}%

\global\long\def\AINormS#1{\|#1\|_{\matA^{-1}}^{2}}%

\global\long\def\BINormS#1{\|#1\|_{\matB^{-1}}^{2}}%

\global\long\def\BINorm#1{\|#1\|_{\matB^{-1}}}%

\global\long\def\ONorm#1#2{\|#1\|_{#2}}%

\global\long\def\T{\textsc{T}}%

\global\long\def\pinv{\textsc{+}}%

\global\long\def\Expect#1{{\mathbb{E}}\left[#1\right]}%

\global\long\def\ExpectC#1#2{{\mathbb{E}}_{#1}\left[#2\right]}%

\global\long\def\dotprod#1#2{(#1,#2)}%

\global\long\def\dotprodX#1#2#3{(#1,#2)_{#3}}%

\global\long\def\dotprodM#1#2{(#1,#2)_{\matM}}%

\global\long\def\dotprodsqr#1#2{(#1,#2)^{2}}%

\global\long\def\Trace#1{{\bf Tr}\left(#1\right)}%

\global\long\def\dist#1{{\bf dist}\left(#1\right)}%

\global\long\def\vectorization#1{{\bf vec}\left(#1\right)}%

\global\long\def\vecskew#1{{\bf vec_{skew}}\left(#1\right)}%

\global\long\def\nnz#1{{\bf nnz}\left(#1\right)}%

\global\long\def\blockdiag#1{{\bf blkdiag}\left(#1\right)}%

\global\long\def\vol#1{{\bf vol}\left(#1\right)}%

\global\long\def\rank#1{{\bf rank}\left(#1\right)}%

\global\long\def\diag#1{{\bf diag}\left(#1\right)}%

\global\long\def\grad#1{{\bf grad}#1}%

\global\long\def\hess#1{{\bf Hess}#1}%

\global\long\def\sym#1{{\bf sym}#1}%

\global\long\def\skew#1{{\bf skew}#1}%

\global\long\def\st{\,\,\,\text{s.t.}\,\,\,}%

\global\long\def\elp{\mathbb{S}^{\matB}}%

\global\long\def\elpCCA{\mathbb{S}_{\x\y}}%

\global\long\def\elpa{\mathbb{S}^{\mat A}}%

\global\long\def\elplad{\mathbb{S}^{\mat{\mat{S_{w}}+\lambda\matI_{d}}}}%

\global\long\def\Stiefellda{{\bf St}_{(\mat{S_{w}}+\lambda\matI_{d})}(p,d)}%

\global\long\def\ldaB{\mat{S_{w}}+\lambda\matI_{d}}%

\global\long\def\elpsigx{\mathbb{S}^{\mat{\Sigma_{\x\x}}}}%

\global\long\def\elpsigy{\mathbb{S}^{\mat{\Sigma_{\y\y}}}}%

\global\long\def\elpparam#1{\mathbb{S}^{#1}}%

\global\long\def\poly#1{{\bf poly}\left(#1\right)}%

\global\long\def\id{{\bf id}}%

\global\long\def\stiefel{{\bf St}}%

\global\long\def\stiefelB{{\bf St}_{\matB}}%

\global\long\def\qf#1{{\bf qf}\left(#1\right)}%

\global\long\def\qfm#1#2{{\bf qf}_{#2}\left(#1\right)}%

\global\long\def\qfmsmall#1#2{{\bf qf}_{#2}(#1)}%

\global\long\def\fcca{f_{{\bf CCA}}(\matZ)}%

\global\long\def\justfcca{f_{{\bf CCA}}}%

\global\long\def\sigmacca{\Sigma_{\nabla^{2}f_{{\bf CCA}}}}%

\global\long\def\flda{f_{{\bf FDA}}(\matW)}%

\global\long\def\justflda{f_{{\bf FDA}}}%

\global\long\def\nicehalf{\nicefrac{1}{2}}%

\newcommand{\BS}[1]{{\color{blue}BS: #1}}
\newcommand{\UM}[1]{{\color{olive}UM: #1}}
\newcommand{\HA}[1]{{\color{green}HA: #1}}

\newcommand{\reply}[1]{{\color{blue}#1}}
\newcommand{\DONE}{\reply{Done.}}

\title{Faster Randomized Methods for Orthogonality Constrained Problems}

\author{Boris Shustin and Haim Avron\\
Tel Aviv University}

\maketitle
\begin{abstract}
Recent literature has advocated the use of randomized methods for
accelerating the solution of various matrix problems arising in
machine learning and data science. One popular strategy for leveraging randomization in numerical linear algebra is to use it as a way to reduce problem size. However, methods based on this strategy lack sufficient 
accuracy for some applications. Randomized preconditioning is another approach for leveraging randomization in numerical linear algebra, which provides higher accuracy. The main challenge in using randomized preconditioning is the need for an underlying iterative method, thus 
randomized preconditioning so far has been applied almost exclusively to solving regression problems and linear systems. In this article, we show how to expand the application of randomized
preconditioning to another important set of problems prevalent in machine learning: optimization
problems with (generalized) orthogonality constraints. We demonstrate our approach, which is based on the framework of Riemannian optimization and Riemannian preconditioning, on the problem
of computing the dominant canonical correlations and on the Fisher
linear discriminant analysis problem. More broadly, our method is designed for problems with input matrices featuring one dimension much larger than the other (e.g., the number of samples much larger than the number of features). For both problems, we evaluate
the effect of preconditioning on the computational costs and asymptotic convergence
and demonstrate empirically the utility of our approach.
\end{abstract}

\section{Introduction}

Matrix sketching has recently emerged as a powerful tool for accelerating
the solution of many important matrix computations, with widespread
use throughout machine learning and data science. Important examples
of matrix computations that have been accelerated using matrix sketching
include linear regression, low rank approximation, and principal component
analysis (see recent surveys \cite{woodruff2014sketching,yang2016implementing}).
Matrix sketching is one of the main techniques used in so-called \emph{Randomized
Numerical Linear Algebra}.

Roughly speaking, matrix sketching provides a transformation that embeds
a high dimensional space in a lower dimensional space, while preserving
some desired properties of the high dimensional space \cite{woodruff2014sketching}.
There are several ways in which such an embedding can be used. The
most popular approach is \emph{sketch-and-solve}, in which matrix sketching
is used to form a smaller problem. That is, sketch-and-solve based
algorithms attempt to find a ``good'' approximate
solution by sketching the input problem so that with high probability
the exact solution of the sketched problem is a good approximate solution
to the original problem. For example, the sketch-and-solve approach
for solving unconstrained overdetermined linear regression problem,
i.e., $\min_{\w}\TNorm{\matX\w-\y}$, where $\matX\in\R^{n\times d}$
($n\geq d$) is assumed to be full rank matrix, is to randomly generate
a matrix $\matS\in\R^{s\times n}$ and solve the reduced size problem $\min_{\w}\TNorm{\matS\matX\w-\matS\y}$~\cite{drineas2011faster}. 
If $\matS$ is chosen appropriately, then with high
probability the solution of the sketched problem, $\hat{\w}$, is
close to the exact solution, $\w^{\star}$, in the sense that the
inequality $\TNorm{\matX\hat{\w}-\y}\leq(1+\varepsilon)\TNorm{\matX\w^{\star}-\y}$
holds, and that $\matS\matX$ and $\matS\y$ can be computed quickly.

Sketch-and-solve algorithms have been proposed for a wide spectrum
of linear algebra problems relevant for machine learning applications: linear regression \cite{drineas2011faster},
principal component analysis \cite{kannan2014principal}, canonical
correlation analysis (CCA) \cite{ABTZ14}, kernelized methods \cite{ANW14},
low-rank approximations \cite{CW17}, structured decompositions \cite{boutsidis2017optimal},
etc. However, there are two main drawbacks to the sketch-and-solve
approach. First, it is unable to deliver highly accurate results
(typically, for sketch-and-solve algorithms, the running time dependence
on the accuracy parameter $\epsilon$ is $\Theta(\epsilon^{-2})$).
The second drawback is that sketch-and-solve algorithms typically
have only Monte-Carlo type guarantees, i.e., they return a solution
within the prescribed accuracy threshold only with high probability
(on the positive side, the running time is deterministic).

These drawbacks have prompted researchers to develop a second approach,
typically termed \emph{sketch preconditioning} or \emph{randomized
preconditioning. }The main idea in randomized preconditioning is to
use an iterative method which, in turn, uses a preconditioner that
is formed using a sketched matrix. For example, consider again the
unconstrained overdetermined linear regression problem problem. It
is possible to accelerate the solution of $\min_{\w}\TNorm{\matX\w-\y}$
by first sketching the matrix $\matX$ to form $\matS\matX$, and
then using a factorization of $\matS\matX=\matQ\matR$ to form a preconditioner,
$\matR$, for an iterative Krylov method (e.g., LSQR). By choosing
$\matS$ properly, with high probability the preconditioner $\matR$
is such that the condition number governing the convergence of the
Krylov method, $\kappa(\matX\matR^{-1})$, is bounded by a small constant~\cite{RT08,avron2010blendenpik,meng2014lsrn,CW17,GOS16}.
Thus, when using a Krylov method to solve $\min_{\w}\TNorm{\matX\matR^{-1}\w-\y}$
only a small number of iterations are necessary for convergence.

More generally, by using an iterative method, it is typically possible
to reduce the running time dependence on the accuracy parameter to
be logarithmic instead of polynomial. Furthermore, since we can control
the stopping criteria of the iterative methods, sketch preconditioning
algorithms typically entertain Las-Vegas type guarantees, i.e., they
return a solution within the accuracy threshold, albeit at the cost of 
probabilistic running time.

The sketch-and-solve approach is more prevalent in the literature
than sketch preconditioning. Indeed, in one way or the other, almost
all sketch preconditioning methods have essentially been designed for
linear regression or solving linear systems. The main reason is that
sketch preconditioning requires an iterative method that
can be preconditioned, and such a method is not always known for the
various problems addressed by sketching. Indeed, linear regression
and solving linear systems are cases where the use of preconditioning
is straightforward.  

The goal of this paper is to go beyond linear regression, and design
sketch preconditioning algorithms for another important class of problems:
optimization problems under generalized orthogonality constraints. We would like to emphasize that in the context of optimization, preconditioning can be performed on some algorithms, e.g., algorithms based on conjugate gradients ~\cite{vandereycken2010riemannian}. However, our approach aims to be more general; rather than preconditioning each optimization method on its own, we aim to precondition the problem itself through randomized preconditioning.
In general, we are interested in solving problems of the form 
\begin{equation}
\min f(\matX_{1},\dots,\matX_{k})\st\matX_{i}^{\T}\matB_{i}\matX_{i}=\matI_{p}\quad(i=1,\dots,k)\label{eq:general_problem}
\end{equation}
where $f(\matX_{1},\dots,\matX_{k}):\R^{d_{1}\times p}\times\cdots\times\R^{d_{k}\times p}\to\R$
is a smooth function, and $\matB_{i}\in\R^{d_{i}\times d_{i}}$ are
fixed symmetric positive definite (SPD) matrices. An important example
of such problems is the problem of finding the dominant subspace:
given a symmetric matrix $\matA\in\R^{n\times n}$, if $f(\matX_{1})=-\Trace{\matX_{1}^{\T}\matA\matX_{1}}$
and $\matB_{1}=\matI_{n},$ then Problem~(\ref{eq:general_problem})
corresponds to finding a basis for the dominant eigenspace.

Problems of the form of Eq.~(\ref{eq:general_problem}) frequently
appear in machine learning. In this paper, we focus on a specific important subset of these problems; the constraint matrices $\matB_{i}$ are not given explicitly, however, as is typically happening in various machine learning problems, we assume that they can be written as a Gram matrix
of a tall-and-skinny data matrix $\matZ_{i}$, i.e., $\matB_{i}=\matZ_{i}^{\T}\matZ_{i}$
for $\matZ_{i}\in\R^{n_{i}\times d_{i}}$. More generally, we allow the following form $\matB_{i}=\matZ_{i}^{\T}\matZ_{i}+\lambda_{i}\matI_{d}$
where $\lambda_{i}\geq0$ is some regularization parameter. Moreover, we assume that $n_{i} \gg d_{i}$. Indeed, that is often the case that $\matZ_{i}$ represent data matrices of samples stacked in rows, such that the number of samples, $n_{i}$, is much larger than the dimension of each data point, $d_{i}$. We aim to address the class of problems where $\matZ_{1},\dots,\matZ_{k}$ are given as inputs and not $\matB_{1}, \dots, \matB_{k}$, and $n_{i} \gg d_{i}$. Both requirements are fundamental to our paper. Thus, our goal is to avoid explicitly forming $\matB_{i}$ or factorizing it (e.g., Cholesky decomposition) computations that require $O(nd^2)$ operations. We aim to design
algorithms that use $o(nd^{2})$ operations to find the preconditioner,
and use $O(nd)$ operations per iteration, where  $n=\max_{i}n_{i}$ and $d=\max_{i}d_{i}$. 

Two important unsupervised machine learning methods that reduce to
Problem~(\ref{eq:general_problem}) are canonical correlation analysis
(CCA) and Fisher linear discriminant analysis (FDA). 
We illustrate our approach, demonstrating its effectiveness both theoretically and empirically,  
on both of these problems.
In particular, we improve on the $\Theta(nd^{2})$ running time possible for both problems
using direct methods. Moreover, $\Theta(nd^{2})$ is also the cost of forming $\matB_{i}$, which is one of the reasons for prohibiting the explicit formation of $\matB_{i}$ in our setting.  For CCA and FDA, the proposed approach can be viewed as a novel way of solving an old problem.

The underlying iterative methods we precondition are based on the
framework of \emph{Riemannian optimization} \cite{AMS09,boumal2022intromanifolds}.
Riemannian optimization is well suited for problems with manifold
constraints, e.g., under generalized orthogonality constraints \cite{EAS98}.
It makes use of the Riemannian geometry components associated with
the constraining manifolds, which in the case of Eq.~(\ref{eq:general_problem}) are
products of generalized Stiefel manifolds. \emph{Riemannian Preconditioning}~\cite{MS16},
which is a technique for preconditioning Riemannian optimization algorithms based on carefully choosing the Riemannian metric is another component of our proposed method. 
By combining randomized preconditioning with Riemannian preconditioning we obtain randomized preconditioners and
faster methods for solving Eq.~(\ref{eq:general_problem}). Specifically, we propose a constant randomized preconditioning scheme $\matX\mapsto\matM_{\matX}\coloneqq\matM$ that is an SPD matrix $\matM\in \R^{d\times d}$ which defines the Riemannian metric on $\stiefelB(p,d)$ (see Section \ref{sec:rand-precond}). 

\begin{remark}\label{rem:non-constant-metric}
    We note that our proposed preconditioner is based on the search space, as we assume that the main computational burden lies in forming $\matB_{i}$. Nevertheless, our method can be combined with the approach presented in~\cite{MS16}, i.e., also include the a cost function dependent component, which would make the preconditioner non-constant. We leave this for future work, and in particular one example of this kind of preconditioner can be found in \cite{mor2020solving}.
\end{remark}

\subsection{\label{subsec:Contributions}Contributions}
In this subsection, we emphasize our main contributions. In this paper, we propose to expand the application of randomized preconditioning to optimization problems with (generalized) orthogonality constraints. In order to achieve that goal we use the framework of Riemannian optimization \cite{AMS09,boumal2022intromanifolds}, and perform Riemannian preconditioning via Riemannian metric selection \cite{MS16}. Thus, we utilize the Riemannian components of the generalized Stiefel manifold developed in \cite{shustin2021randomized}. 
\begin{itemize}
\item Our main contribution, presented in Section \ref{sec:rand-precond}, is designing a constant randomized preconditioning scheme which is incorporated via the Riemannian metric for the generalized Stiefel manifold (and products of it) with $\matB_{i}=\matZ_{i}^{\T}\matZ_{i}+\lambda\matI_{d}$ where $\matZ_{i}$ are given. Moreover, we exemplify a specific sketching transformation, \noun{CountSketch} \cite{CCF04}, list the corresponding computational costs of the geometric components required for sketching based Riemannian preconditioning in Table \ref{tab:costs} (we also provide costs for the \emph{Subsampled Randomized Hadamard Transform} (SRHT) \cite{tropp2011improved,boutsidis2013improved, CW17}), and present Lemma \ref{lem:sketching} (and its analogue for SRHT, Lemma \ref{lem:sketching_SRHT}), which would later aid in assessing the effectiveness of our proposed preconditioning.
\item In sections \ref{sec:Randomized-preconditioning-for-CCA} and \ref{sec:Randomized-preconditioning-for-LDA}, we exemplify our approach by developing preconditioned iterative methods
for CCA and FDA correspondingly. We identify the corresponding critical points, and analyze their stability. We theoretically analyze the effect of a generic preconditioner on the asymptotic
convergence rate, and identify the optimal preconditioner. Then, we present end-to-end sketching based algorithms
for CCA and FDA, analyze their asymptotic convergence and evaluate the computational costs. Finally, we demonstrate numerically our randomized preconditioning approach in Section \ref{sec:Numerical-Experiments}.
\end{itemize}
\subsection{\label{subsec:Related-Work}Related Work}

We begin with the related work on \textbf{matrix sketching}. The literature on sketching has so far mostly focused on the sketch-and-solve approach. In particular, in the context of solving problems with generalized orthogonality
constraints, a sketch-and-solve based approach for CCA was developed
in \cite{ABTZ14}. Sketch preconditioning was predominantly applied
to linear least-squares regression problems (e.g., \cite{avron2010blendenpik,meng2014lsrn}),
and also to non least-squares variants (e.g., $\ell_{1}$-regression \cite{MM13}). For a regularized version of FDA, a randomized iterative method
was developed in \cite{chowdhury2018randomized}. Even though the method in \cite{chowdhury2018randomized} does not use sketch preconditioning per se, it can be viewed as a preconditioned Richardson iteration. 
A randomized iterative method for solving LP problems was proposed in~\cite{CLAD20},
where sketch preconditioning is used to find the Newton search direction as part 
of an interior point method. To the
best of our knowledge, sketch preconditioning has neither been applied
to CCA before nor has it been used in the context of Riemannian optimization. We remark, that one approach to solving CCA is via alternating least squares (ALS) presented in \cite{golub1995canonical} (see Algorithm 5.2 in their paper). Thus, one natural idea is to use randomized preconditioning for least squares for each of the sub-problems in ALS approach. While a series of least squares solved exactly is equivalent to a specific choice of Riemannian algorithm, and a prescribed step size, our preconditioning approach allows for a wide variety of algorithms (e.g., first, and second-order methods). In addition, solving a series of least squares requires to solve each of the sub-problems, whereas in our approach one does not need to solve sub-problems. Finally, on a conceptual level our approach aims to precondition the problem directly, rather than performing a two-level solution of replacing the problem with a series of sub-problems, and then preconditioning these.


Next, we recall the related work on \textbf{Riemannian optimization}. Recent works introducing Riemannian optimization are \cite{AMS09,boumal2022intromanifolds}.
Earlier works are \cite{luenberger1972gradient,gabay1982minimizing,smith1994optimization},
and specifically, using Riemannian optimization to solve problems
under orthogonality constraints is presented in the seminal work \cite{EAS98}.

In particular, we mention related work on \textbf{Riemannian preconditioning}. Preconditioning of Riemannian optimization methods based on the cost
function alone is presented in several works (see e.g., \cite{ngo2012scaled,mishra2014r3mc,shi2016low,zhou2016riemannian}).
Most of the aforementioned works attempt to perform preconditioning
by approximating the Euclidean Hessian of the cost function. However,
the Riemannian Hessian and the Euclidean Hessian are quite different
even for simple examples (see \cite[Section 3.5]{shustin2021randomized}).
The Riemannian Hessian is related to the Hessian of the Lagrangian,
and indeed in \cite{MS16} it is shown that selecting
the Riemannian metric inspired by the Hessian of the Lagrangian affects
convergence of Riemannian steepest-descent, coining the term \emph{Riemannian
preonditioning }for judiciously choosing the Riemannian metric in
order to accelerate convergence. Unlike in \cite{MS16}, we present a randomized
preconditioning strategy, and analyze the condition number of the
Riemannian Hessian at the optimum for specific examples (CCA and FDA),
which allows us to quantify the quality of the proposed preconditioner. Recently, a new paper, \cite{gao2023optimization}, that appeared after our paper was sent for initial review, presented a preconditioning technique for product manifolds by defining the Riemannian metric via the diagonal blocks of the Riemannian Hessian at the optimum. In addition, Gao et. al analyzed the effect of preconditioning via the condition number of the Riemannian Hessian at the optimum, in a similar manner to the analysis in our paper. They demonstrated their technique on CCA, truncated singular value decomposition, and matrix and tensor completion problems. Unlike in \cite{gao2023optimization}, we present a randomized preconditioning strategy via matrix sketching.

Additional works regarding preconditioning of Riemannian methods are a Riemannian analogue
to the preconditioned Richardson method~\cite{kressner2016preconditioned},
and the works ~\cite[Algorithm 5]{baker2008riemannian} and ~\cite{vandereycken2010riemannian} where a Riemannian
version of preconditioning the trust-region sub-problem was proposed. To some extent that work can be also viewed as Riemannian metric selection. The aforementioned preconditioning approaches are essentially two-stage approaches, in which a specific choice of a Riemannian algorithm is made, and then a preconditioning is applied specific to this method. In contrast, our approach aims to allow a wider choice of algorithms (zero, first, and second-order methods), without the need to solve sub-problems. Moreover, from a conceptual approach, we aim to precondition the problem itself, regardless of the algorithm at hand, and avoid two-stage procedures.

Another work on Riemannian preconditioning for the trust-region sub-problem is \cite{mor2020solving}.
We remark that apart from the last work mentioned, no previous work proposed to use randomized preconditioners in the context of Riemannian preconditioning.


Finally, we also recall related work on \textbf{iterative methods for CCA and FDA}. Several iterative methods for solving CCA have been proposed in the
literature. In \cite{golub1995canonical}, an iterative method for CCA was presented based
on ALS. Each iteration
requires the solution of two least squares problems. The authors suggest
using LSQR for that task. In \cite{WangEtAl16}, it was proposed to replace LSQR with
either accelerated gradient descent, stochastic variance reduce gradient
(SVRG) or accelerated SVRG. They also proposed
a different approach based on shift-and-invert preconditioning. In \cite{MLF15}, an algorithm for CCA based on augmented approximate
gradients was developed. In \cite{ge2016efficient}, an iterative algorithm
for the generalized eigenvalue problem was provided, and a standard reduction
of CCA to generalized eigenvalue problems was used to derive an algorithm for
CCA. They assume a fast black box access
to an approximate linear system solver.

Convergence bounds of all the aforementioned algorithms depend on
the condition number of the input matrices, which might be large.
In contrast, the condition number bounds for our proposed sketching
based algorithms are independent of the conditioning of the input
matrices. As aside, in \cite{yger2012adaptive}, a Riemannian method for
adaptive CCA was proposed.

Recently, in \cite{chapmanunconstrained}, a stochastic method to solve large-scale CCA problems was proposed. The method is based on unifying multiview CCA and self-supervised learning. To make the method feasible for large-scale problems a family of fast stochastic algorithms was proposed for a novel unconstrained objective function for generalized eigenvalue problems (which CCA generalizes).

In the context of FDA, in a recent work, \cite{chowdhury2018randomized},
an iterative, sketching-based algorithm for regularized FDA was proposed.

Note that both CCA and FDA are closely related to the singular value decomposition (SVD) and generalized eigenvalue problems correspondingly. In the context of optimization on manifolds, an analysis of SVD and generalized eigenvalue problems correspondingly, as dynamical systems, i.e., gradient flows on manifolds (of unitary matrices) was done in \cite{helmke1992singular} and \cite{brockett1991dynamical}. In particular, the stability of the equilibrium point was established, which we use in Subsections \ref{sec:CCAstability} and \ref{sec:LDAstability}.

\section{\label{sec:Preliminaries}Preliminaries}
In this section, we present notation and recall preliminaries on Riemannian optimization and preconditioning. 

\subsection{Notation and Basic Definitions}

Scalars are denoted by lower case letters $\alpha,\beta,...,x,y,\dots$,
while vectors are denoted by bold letters $\boldsymbol{\alpha,\beta},...,\boldsymbol{x,y},\dots$,
and matrices are denoted by bold uppercase English letters $\matA,\matB\dots$
or upper case Greek letters. We use the convention that vectors are
column-vectors. We describe
a diagonal matrix using $\diag{\cdot}$,
and a block diagonal matrix using $\blockdiag{\cdot}$.

We denote by $\dotprod{\cdot}{\cdot}_{\matM}$ the inner-product with
respect to a matrix $\matM$: for $\matU$ and $\matV$, $\dotprodM{\matU}{\matV}\coloneqq\Trace{\matU^{\T}\matM\matV}$
where $\Trace{\cdot}$ denotes the trace operator. The $s\times s$
identity matrix is denoted $\matI_{s}$. The $s\times s$ zero matrix
is denoted $\matzero_{s}$. We denote by ${\cal {\cal S}_{\text{sym}}}(p)$
and ${\cal {\cal S}_{\text{skew}}}(p)$ the set of all symmetric and
skew-symmetric matrices (respectively) in $\R^{p\times p}$. The symmetric
and skew-symmetric components of $\matA$ are denoted by $\sym{(\matA)}\coloneqq\left(\matA+\matA^{\T}\right)/2$
and $\skew{(\matA)}\coloneqq\left(\matA-\matA^{\T}\right)/2$ respectively.

For a SPD matrix $\matB\in\R^{d\times d}$, we denote by $\matB^{\nicehalf}$
the unique SPD matrix such that $\matB=\matB^{\nicehalf}\matB^{\nicehalf}$,
obtained by keeping the same eigenvectors and taking the square root
of the eigenvalues. We denote the inverse of $\matB^{\nicehalf}$
by $\matB^{-\nicehalf}$.

The eigenvalues of a symmetric $d\times d$ matrix $\matA$ are denoted
by $\lambda_{1}(\matA)\geq\lambda_{2}(\matA)\geq\dots\geq\lambda_{d}(\matA)$.
The condition number of $\matA$, i.e. the ratio between the largest and smallest
eigenvalues in absolute value, is denoted by $\kappa({\matA})$.
The generalized condition number of a matrix pencil $(\matA,\matB)$, which is the ratio between the largest and the smallest generalized eigenvalues, is denoted by $\kappa(\matA,\matB)$, and it holds that if $\matB$ is a SPD matrix, then  $\kappa(\matA,\matB)=\kappa(\matB^{-\nicehalf}\matA\matB^{-\nicehalf})$.

For differential geometry objects we use the following notations:
tangent vectors (of a manifold) are denoted using lower case Greek
letters with a subscript indicating the point of the manifold for
which they correspond (e.g., $\eta_{\x}$), while normal vectors (of
a manifold) are denoted using bold lower and upper case letters with
a subscript (e.g., $\u_{\x}$). We denote by $\stiefelB(p,d)$ the
\emph{generalized Stiefel manifold},  which is a submanifold of $\R^{d\times p}$ defined by
\[
\stiefelB(p,d)\coloneqq\left\{ \matX\in\R^{d\times p}:\ \matX^{\T}\matB\matX=\matI_{p}\right\} \,.
\]
Given a function defined on $\stiefelB(p,d)$, a smooth extension
of that function to the entire $\R^{d\times p}$ is denoted by a bar
decorator. For example, given a smooth function $f:\stiefelB(p,d)\to\R$,
we use $\bar{f}:\R^{d\times p}\to\R$ to denote a smooth function
on $\R^{d\times p}$ such that on $\stiefelB(p,d)$ that function
agrees with $f$.

\subsection{\label{subsec:Preconditioned-Geometry-for}Riemannian Optimization
and Preconditioning for the Generalized Stiefel Manifold}

In this subsection we detail several necessary ingredients required for preconditioning strategy: Riemannian optimization, Riemannian preconditioning, and
the geometric optimization components of the generalized Stiefel manifold with a
non-standard metric. In following sections
we make use of these components to present our randomized Riemannian
preconditioning strategy.

Riemannian optimization \cite{AMS09,boumal2022intromanifolds} is
a framework of designing algorithms for solving constrained optimization problems, 
where the constraints form a smooth manifold ${\cal M}$. The general idea of
these algorithms is to make use of the differential geometry components
of ${\cal M}$ in order to generalize iterative methods for unconstrained
optimization. Iterative algorithms for smooth problems such as gradient
methods, trust-region, and conjugate gradient (CG), are adapted to
the Riemannian setting using the following components: a \emph{Riemannian
metric }which is a smoothly varying inner product $\x\mapsto g_{\x}$
on the tangent bundle $T{\cal M}$ such that $({\cal M},g)$ becomes
a Riemannian manifold, a \emph{retraction mapping} $R_{\x}:T_{\x}{\cal M}\to{\cal M}$
which allows to take a step at point $\x\in{\cal M}$ in a direction
$\xi_{\x}\in T_{\x}{\cal M}$, a \emph{vector transport }${\cal T}_{\eta_{\x}}:T_{\x}{\cal M}\to T_{R_{\x}(\eta_{\x})}{\cal M}$
which allows operations between tangent vectors from two different
tangent spaces, a \emph{Riemannian gradient} $\grad{f(\x)\in T_{\x}{\cal M}}$,
and a \emph{Riemannian Hessian $\hess{f(\x)}:T_{\x}{\cal M}\to T_{\x}{\cal M}$.
}Another important notion in the context of this paper is the notion
of Riemannian submanifold, which allows to easily compute the aforementioned
geometric components if they are known for the ambient manifold.\emph{
}Usually, a Riemannian optimization algorithm is built from iterations
on the tangent space $T_{\x}{\cal M}$ which are then retracted to
the manifold. For example, Riemannian gradient methods are given by
the formula $\x_{k+1}=R_{\x_{k}}(\tau_{k}\grad{f(\x_{k})})$ where $\tau_{k}$ is the step size. Many of the Riemannian algorithms and common manifolds are implemented
in \noun{manopt,} which is a \noun{matlab} library \cite{manopt}.
There is also a \noun{python} parallel for \noun{manopt} called \noun{pymanopt}
\cite{townsend2016pymanopt}. The experiments reported in Section
\ref{sec:Numerical-Experiments} use the \noun{manopt} library.

Preconditioning of iterative methods is often challenging since it
is not initially clear how to actually precondition the problem for an iterative method. Riemannian preconditioning \cite{MS16} generalizes
preconditioning for Riemannian optimization methods. In Riemannian
preconditioning, preconditioning is performed via a Riemannian metric
selection, i.e., the preconditioner is incorporated via the Riemannian
metric.

A motivation for metric selection is the observation that the condition number of the Riemannian Hessian at the
optimum affects the asymptotic convergence rate of Riemannian optimization, e.g.,
\cite[Theorem 4.5.6, Theorem 7.4.11 and Equation 7.50]{AMS09}, and
\cite[Section 4.6]{boumal2022intromanifolds}.
In the case of convex objective function in the Riemannian sense
\cite[Chapter 3.2]{udriste2013convex} there are also global convergence
results (e.g., \cite[Chapter 7, Theorem 4.2]{udriste2013convex}).
Thus, selecting a metric such that this condition number is lowered
should improve convergence.

In this paper, we propose to utilize randomized preconditioning
to accelerate the solution of orthogonality constrained problems by
using Riemannian preconditioning. Optimization problems under generalized orthogonality constraints can be written as Riemannian optimization problems on the generalized
Stiefel manifold. The standard Riemannian metric for the generalized
Stiefel manifold (e.g., \cite{EAS98,manopt}) is 
\[
g_{\matX}(\xi_{\matX},\eta_{\matX})\coloneqq\dotprod{\xi_{\matX}}{\eta_{\matX}}_{\matB}=\Trace{\xi_{\matX}^{\T}\matB\eta_{\matX}}.
\]
The use of this metric has two possible shortcomings. First, it might not be the optimal metric with respect to convergence of Riemannian
optimization methods. Second, and more relevant for this paper, the
computation of most Riemannian components requires taking products
of $\matB^{-1}$ with a vector. Computing $\matB$ and/or factorizing
is in many cases as costly as solving the problem directly or even solving it via standard iterative methods, e.g., for
CCA (see Section~\ref{sec:Randomized-preconditioning-for-CCA}).

In order to overcome these shortcomings, we leverage our recently developed
geometric components for the generalized Stiefel manifold with a Riemannian
metric which is based on a preconditioning scheme $\matX\mapsto\matM_{\matX}$ that maps a $\matX\in\stiefelB(p,d)$ to an SPD matrix $\matM_{\matX}\in \R^{d\times d}$, instead using the standard metric defined by $\matB$ \cite{shustin2021randomized}. In this paper,
we apply a preconditioning scheme independent of $\matX$, i.e., $\matM_{\matX}\coloneqq\matM$
for all $\matX\in\stiefel_{\matB}(p,d)$ (see Section \ref{sec:rand-precond}).
As a reference we summarize (without proofs) in Table \ref{tab:gen_Stiefel_comps} the various geometric components 
developed in~\cite{shustin2021randomized}.

In addition, we also address optimization problems that are constrained on a product
of generalized Stiefel manifolds such as CCA. In this case, the Riemannian
components of the generalized Stiefel manifold can be generalized
to be the Riemannian components of a product of generalized Stiefel
manifolds in a straightforward way~\cite{shustin2021randomized}.

\begin{table}[t]
\caption{\label{tab:gen_Stiefel_comps}Summary of the Preconditioned various Riemannian components
of the generalized Stiefel manifold \cite{shustin2021randomized}, i.e., a generalized Stiefel manifold with a Riemannian metric defined via an SPD matrix $\matM\in \R^{d\times d}$. Note that, the matrix $\matS_{\xi_{\matX}}$ denotes a solution for a Sylvester equation, arising from the orthogonal projection of $\xi_{\matX}$ on the tangent space.}

\centering{}{\scriptsize{}}%
\begin{tabular}{|>{\centering}p{5cm}|>{\centering}p{9.5cm}|}
\hline 
\textbf{\footnotesize{}Component} & \textbf{\footnotesize{}Formula}\tabularnewline
\hline 
\hline 
{\scriptsize{}Tangent space $T_{\matX}\stiefel_{\matB}(p,d)$} & {\scriptsize{}$\left\{ \matZ\in\R^{d\times p}\,:\,\matZ^{\T}\matB\matX+\matX^{\T}\matB\matZ=\mat 0_{p}\right\} $}\tabularnewline
\hline 
{\scriptsize{}Tangent space $T_{\matX}\stiefel_{\matB}(p,d)$ alternative
form} & {\scriptsize{}$\left\{ \matZ=\matX\mat{\Omega}+\matX_{\matB\perp}\mat K\in\R^{d\times p}\,:\,\mat{\Omega}=-\mat{\Omega}^{\T},\ \mat K\in\R^{(d-p)\times p}\right\} $
such that $\matX_{\matB\perp}^{\T}\matB\matX_{\matB\perp}=\matI_{d-p},\ \matX_{\matB\perp}^{\T}\matB\matX=\mat 0_{(d-p)\times p}$}\tabularnewline
\hline 
{\scriptsize{}The Riemannian metric} & {\scriptsize{}$g_{\matX}(\xi_{\matX},\eta_{\matX})=\dotprodM{\xi_{\matX}}{\eta_{\matX}}=\Trace{\xi_{\matX}^{\T}\matM\eta_{\matX}}$}\tabularnewline
\hline 
{\scriptsize{}Normal space with respect to the Riemannian metric} & {\scriptsize{}$\left(T_{\matX}\text{St}_{\matB}(p,d)\right)^{\perp}=\left\{ \matM^{-1}\matB\matX\matS\ :\ \matS=\matS^{\T}\right\} $}\tabularnewline
\hline 
{\scriptsize{}Polar-based retraction } & {\scriptsize{}$R_{\matX}^{\text{polar}}(\xi_{\matX})=(\matX+\xi_{\matX})(\mat I_{p}+\xi_{\matX}^{\T}\mat B\xi_{\matX})^{-\nicehalf}$}\tabularnewline
\hline 
{\scriptsize{}QR-based retraction} & {\scriptsize{}$R_{\matX}^{\text{QR}}(\xi_{\matX})=\matB^{-\nicehalf}\qf{\matB^{\nicehalf}\left(\matX+\xi_{\matX}\right)}$}\tabularnewline
\hline 
{\scriptsize{}Orthogonal projection on the normal space } & {\scriptsize{}$\mat{\Pi}_{\matX}^{\perp}\left(\xi_{\matX}\right)=\matM^{-1}\matB\matX\matS_{\xi_{\matX}}$ }\tabularnewline
\hline 
{\scriptsize{}Orthogonal projection on the tangent space} & {\scriptsize{}$\mat{\Pi}_{\matX}\left(\xi_{\matX}\right)=\left(\id_{T_{\matX}\text{St}_{\matB}(p,d)}-\mat{\Pi}_{\matX}^{\perp}\right)\left(\xi_{\matX}\right)=\xi_{\matX}-\matM^{-1}\matB\matX\matS_{\xi_{\matX}}$ }\tabularnewline
\hline 
{\scriptsize{}Sylvester equation for $\matS_{\xi_{\matX}}$ (orthogonal
projection)} & {\scriptsize{}$\left(\matX^{\T}\matB\matM^{-1}\matB\matX\right)\matS_{\xi_{\matX}}+\matS_{\xi_{\matX}}\left(\matX^{\T}\matB\matM^{-1}\matB\matX\right)=\matX^{\T}\matB\xi_{\matX}+\left(\matX^{\T}\matB\xi_{\matX}\right)^{\T}$}\tabularnewline
\hline 
{\scriptsize{}Vector transport} & {\scriptsize{}$\tau_{\eta_{\matX}\xi_{\matX}}=\mat{\Pi}_{R_{\matX}(\eta_{\matX})}\left(\xi_{\matX}\right)$}\tabularnewline
\hline 
{\scriptsize{}Riemannian gradient} & {\scriptsize{}$\grad{f(\matX)}=\mat{\Pi}_{\matX}\left(\matM^{-1}\nabla\bar{f}(\matX)\right)$}\tabularnewline
\hline 
{\scriptsize{}Riemannian Hessian applied on a tangent vector} & {\scriptsize{}$\hess{f(\matX)[\eta_{\matX}]}=\mat{\Pi}_{\matX}(\matM^{-1}(\nabla^{2}\bar{f}(\matX)\eta_{\matX}-\matB\eta_{\matX}(\matX^{\T}\nabla\bar{f}(\matX)-\matX^{\T}\matM\grad f(\matX))))$}\tabularnewline
\hline 
\end{tabular}{\scriptsize\par}
\end{table}

\section{\label{sec:rand-precond}{Randomized Riemannian
Preconditioning}}

In this section we present our main contribution: randomized Riemannian
preconditioners for optimization problems which feature generalized
orthogonality constraints defined by Gram matrices.
That is, our goal is to address constraints of the form $\matX_{i}^{\T}\matB_{i}\matX_{i}=\matI_{p}\quad(i=1,\dots,k)$,
where $\matB_{i}$ are not given explicitly but can be written as a regularized Gram matrix of
a given tall-and-skinny matrix $\matZ_{i}\in\R^{n_{i}\times d_{i}}$, $n_{i} \gg d_{i}$, with
a regularization parameter $\lambda\geq0$, i.e., $\matB_{i}=\matZ_{i}^{\T}\matZ_{i}+\lambda\matI_{d}$.
As alluded to earlier, our proposed solution incorporates a preconditioner
by selecting a non-standard metric for ${\bf St}_{\matB_{1}}\times\dots\times{\bf St}_{\matB_{k}}$.

For simplicity, let us focus our exposition on the case of $k=1$
and drop the subscript from $\matB_{1},\matZ_{1}$ etc. Generalization
to an arbitrary $k$ is straightforward. We define the metric on $\stiefelB(p,d)$
using a constant preconditioning scheme $\matM_{\matX}\coloneqq\matM$ for all $\matX\in\stiefelB(p,d)$ formed from sketching $\matZ$ prior to computing
the Gram matrix. Our proposed construction of $\matM$ requires $o(nd^{2}$)
operations, which is cheaper than computing $\matB$ as required when the standard metric is used.

Our randomized construction of $\matM$ is based on the following
observation: in many cases, if $\matM$ approximates $\matB$ in the
sense that the condition number $\kappa(\matB,\matM)$ is small, then
convergence will be fast. Conceptually, the last observation is synonymous
with the observation that usually the standard metric is a good choice
iteration complexity-wise, albeit a computationally expensive choice,
and thus we should aim at cheaply approximating it. Mathematically,
the underlying reason is that in many cases the condition number of
the Riemannian Hessian at the optimum is bounded by the product of
$\kappa(\matB,\matM)$ and a problem/input dependent quantity which
is typically small. It is known in the Riemannian optimization literature
that the local convergence rate can be analyzed by inspecting the
condition number of the Riemnnian Hessian at the optimum, see \cite[Section 4.1 and Theorem 4.5.6]{AMS09} and
\cite[Section 4.6]{boumal2022intromanifolds}.

We demonstrate the observation that when $\matM$ approximates $\matB$
well, then the Riemannian Hessian at the optimum is well conditioned
on two important use-cases: CCA and FDA (see Theorems~\ref{thm:CCAHessiantheorem}
and \ref{thm:LDAHessiantheorem}). Indeed, we show that under certain
assumptions, both for CCA and FDA we can bound the condition number
at the optimum by the product of some baseline condition number that
depends on eigengaps, and the condition number of $\kappa(\matB,\matM$).

We propose to construct $\matM$ using the technique of sketching \cite{woodruff2014sketching}.
Let $\matS\in\R^{s\times n}$ be some \emph{sketching matrix} (a certain
distribution on matrices; we discuss a concrete choice in the next
paragraph), where $s< n$. We then use $\matM=\matZ^{\T}\matS^{\T}\matS\matZ+\lambda\matI_{d}$.
However, we do not propose to actually compute $\matM$. Using the
metric defined by $\matM$ requires only taking products with $\matM$
and $\matM^{-1}$ (see Table~\ref{tab:gen_Stiefel_comps}). Suppose
we have already computed $\matS\matZ$. The product of $\matM$ with
a vector can be computed in $O(\nnz{\matS\matZ})$ operations. As
for $\matM^{-1}$, by computing a QR factorization of $\bigl[\begin{array}{c}
\matS\matZ\\
\sqrt{\lambda}\matI_{d}
\end{array}\bigr]$ ($O(sd^{2})$ operations), we can obtain a Cholesky factorization
of $\matM$, and then taking the product of $\matM^{-1}$ with a vector
requires $O(d^{2})$. So our goal is to design sketching matrices
$\matS$ such that $\matS\matZ$ is cheap to compute, and for which
$\kappa(\matB,\matM)$ is bounded by a constant.

There are quite a few sketching distributions proposed in the literature,
and most of them are good choices for $\matS$. For concreteness,
we describe a specific choice: the \noun{CountSketch}\footnote{Note that SRHT
is also a good choice, in particular for dense data sets. In Appendix \ref{sec:analogue_of_lem_SRHT}, we provide a similar analysis for SRHT.} transformation
\cite{CCF04}. \noun{CountSketch}
is specified by a random hash function $h:\{1,\dots,d\}\to\{1,\dots,s\}$
and random sign function $g:\{1,\dots,d\}\to\{-1,+1\}$. Applying
$\matS$ to a vector $\x\in \R^n$ is given by the formula 
\[
(\matS\x)_{i}=\sum_{j|h(j)=i}g(j)x_{j}\,.
\]
It is easy to see that $\matS$ is a random matrix where the $j$th
column contains a single nonzero entry $g(j)$ in the $h(j)$th row.
Clearly, $\matS\matZ$ can be computed using $O(\nnz{\matZ})=O(nd)$
arithmetic operations. Thus, it only remains to bound the condition
number. The following lemma shows that if the sketch size is large
enough, then with high probability the condition number is bounded
by a constant.
\begin{lemma}
\label{lem:sketching} Assume that $\lambda>0$ or that $\matZ\in\R^{n\times d}$
has full column rank. Let $s_{\lambda}(\matZ)\coloneqq\Trace{(\matZ^{\T}\matZ+\lambda\matI)^{-1}\matZ^{\T}\matZ}$.
Suppose that $\matS\in \R^{s \times n}$ is a \noun{CountSketch} matrix with $s\geq20s_{\lambda}(\matZ)^{2}/\delta$
for some $\delta\in(0,1)$. Then with probability of at least $1-\delta$
we have that all the generalized eigenvalues of the pencil $(\matZ^{\T}\matZ+\lambda\matI,\matZ^{\T}\matS^{\T}\matS\matZ+\lambda\matI)$
are contained in the interval $[1/2,3/2]$ and $\kappa(\matZ^{\T}\matZ+\lambda\matI,\matZ^{\T}\matS^{\T}\matS\matZ+\lambda\matI)\leq3$.
\end{lemma}

\begin{proof}
The argument is rather standard and appeared in similar forms in the
literature. To show that the generalized eigenvalues of the pencil
$(\matZ^{\T}\matZ+\lambda\matI,\matZ^{\T}\matS^{\T}\matS\matZ+\lambda\matI)$
are contained in the interval $[1/2,3/2]$ and that $\kappa(\matZ^{\T}\matZ+\lambda\matI,\matZ^{\T}\matS^{\T}\matS\matZ+\lambda\matI)\leq3$
to hold, it is enough to show that 
\[
\frac{1}{2}(\matZ^{\T}\matZ+\lambda\matI)\preceq\matZ^{\T}\matS^{\T}\matS\matZ+\lambda\matI\preceq\frac{3}{2}(\matZ^{\T}\matZ+\lambda\matI)\,.
\]
Let $\matZ=\matQ\matR$ be a $\lambda$-QR factorization of $\matZ$, i.e., $\matQ$ is a full-rank matrix and $\matR$ is upper triangular such that $\matR^{\T}\matR = \matZ^{\T}\matZ + \lambda \matI$ \cite[Definition 28]{ACW16b}. Note that such a factorization always exists \cite[remark following Definition 28]{ACW16b}. Left-multiplying by $\matR^{-\T}$ and right-multiplying by $\matR^{-1}$
on both sizes, we find it suffices to show that with probability of
at least $1-\delta$ we have
\[
\frac{1}{2}\matI_{d}\preceq\matQ^{\T}\matS^{\T}\matS\matQ+\lambda\matR^{-\T}\matR^{-1}\preceq\frac{3}{2}\matI_{d}
\]
 or, equivalently, 
\[
\TNorm{\matQ^{\T}\matS^{\T}\matS\matQ+\lambda\matR^{-\T}\matR^{-1}-\matI_{d}}\leq\frac{1}{2}\,.
\]
Since $\matQ^{\T}\matS^{\T}\matS\matQ+\lambda\matR^{-\T}\matR^{-1}-\matI_{d}=\matQ^{\T}\matS^{\T}\matS\matQ-\matQ^{\T}\matQ$
\cite[Fact 29]{ACW16b} and the spectral norm is dominated by the Frobenius
norm, it is enough to show that 
\[
\FNorm{\matQ^{\T}\matS^{\T}\matS\matQ-\matQ^{\T}\matQ}\leq\frac{1}{2}\,.
\]
It is known \cite[Lemma 2]{ANW14} that for any two fixed matrices $\matA$
and $\matB$, and a \noun{CountSketch} matrix $\matS_{0}$ with $m\geq5/(\epsilon^{2}\delta)$
rows, we have that with probability of at least $1-\delta$, 
\[
\FNorm{\matA^{\T}\matS_{0}^{\T}\matS_{0}\matB-\matA^{\T}\matB}\leq\epsilon\cdot\FNorm{\matA}\cdot\FNorm{\matB}\,.
\]
Since $\FNormS{\matQ}=s_{\lambda}(\matZ)$ \cite[Fact 30]{ACW16b}, then with
$s\geq20s_{\lambda}(\matZ)^{2}/\delta$ we have 
\[
\FNorm{\matQ^{\T}\matS^{\T}\matS\matQ-\matQ^{\T}\matQ}\leq\frac{1}{2}
\]
with probability of at least $1-\delta$.
\end{proof}
The last lemma justifies the use of \noun{CountSketch} to form the
preconditioner $\matM$. In practice, additional heuristics that improve
running time and robustness can be inserted into the construction
of randomized preconditioners, and these can improve running time considerably;
see~\cite{avron2010blendenpik}.

Furthermore, our sketching-based preconditioner construction naturally
allows for a warm-start. While this is not captured by our theory,
heuristically (and empirically) the Riemannian optimization part of
our proposed algorithm converges faster if the starting vectors are
close to the optimum. Our sketching approach lets us quickly compute
good starting vectors (these are the sketch-and-solve approximations).
We demonstrate warm-start numerically in Section \ref{sec:Numerical-Experiments}
for CCA and FDA.

In Table \ref{tab:costs} we detail the computational cost, measured
in terms of arithmetic operations, of computing the Riemannian components
of Table \ref{tab:gen_Stiefel_comps}, for our construction of $\matM$
as a preconditioner. Table \ref{tab:costs} is based on \cite[Table 1]{shustin2021randomized}.
Note that all the costs are for operations in ambient coordinates.
In the table, we denote by $T_{\nabla\bar{f}}$ and by $T_{\nabla^{2}\bar{f}}$
the cost of computing the Euclidean gradient and the cost of applying
the Euclidean Hessian to a tangent vector. Instead of committing to
a specific sketch size, we use $s$ for sketch size (number of rows),
and consider two possible sketching distributions: \noun{CountSketch}
and SRHT. 

To put Table \ref{tab:costs} in the contexts of the total computation complexities, we can use the results presented in \cite{boumal2019global}. For example, Riemannian gradient descent with a fixed step-size (which requires a Riemannian gradient and a retraction computation per iteration) takes at most $O(\varepsilon^{-2})$ iterations to reach to an $\varepsilon$-approximate first-order KKT point \cite[Theorem 2.5]{boumal2019global}. Using our proposed preconditioned components, preprocessing takes $O(\nnz{\matZ} + sd^2)$ operations for \noun{CountSketch}, and $O(ndlog(s) + sd^2)$ operations for SRHT. On the other hand, using the standard metric would require $O(nd^2)$ operations in preprocessing. In addition, computing a Riemannian gradient and a retraction is unaffected of the choice of the Riemannian metric after preprocessing is done. Thus, in the worst-case scenario, using our preconditioned components for a fixed-step Riemannian gradient descent would result in a total of $O(\nnz{\matZ}+sd^{2}+\varepsilon^{-2}(ndp+d^{2}p+dp^{2}+T_{\nabla\bar{f}}))$ operations for \noun{CountSketch}, and $O(nd\log(s)+sd^{2}+\varepsilon^{-2}(ndp+d^{2}p+dp^{2}+T_{\nabla\bar{f}}))$ operations for SRHT, while using the standard metric would result in a total of $O(nd^{2}+\varepsilon^{-2}(ndp+d^{2}p+dp^{2}+T_{\nabla\bar{f}}))$ operations. As we demonstrate in Section \ref{sec:Numerical-Experiments}, we have $s \ll n$, thus our preconditioned algorithm is expected to be faster, at least in the worst-case scenario. As for the asymptotic convergence, we demonstrate for CCA and FDA how our preconditioner affects the asymptotic convergence via analysis of the condition number of the Riemannian Hessian at the optimum in Subsection \ref{subsec:Randomized-CCA} and Subsection \ref{subsec:Randomized-LDA}.

\begin{table}[t]
\caption{\label{tab:costs}Summary of the cost of various components for using
sketching based Riemannian preconditioning for generalized orthogonality
constraints.}

\centering{}{\scriptsize{}}%
\begin{tabular}{|>{\centering}p{4cm}|>{\centering}p{5cm}|>{\centering}p{5cm}|}
\hline 
\textbf{\scriptsize{}Operation} & \textbf{\scriptsize{}Cost using}{\scriptsize{} }\noun{\scriptsize{}CountSketch} & \textbf{\scriptsize{}Cost using }{\scriptsize{}SRHT}\tabularnewline
\hline 
\hline 
{\scriptsize{}Preprocessing: computing $\mat S\matZ$} & {\scriptsize{}$O(\nnz{\matZ})=O(nd)$} & {\scriptsize{}$O(nd\log(s))$}\tabularnewline
\hline 
{\scriptsize{}Preprocessing: given $\mat S\matZ$ forming the inverse
of $\matM=\matZ^{\T}\matS^{\T}\matS\matZ+\lambda\matI_{d}$} & {\scriptsize{}$O(sd^{2})$} & {\scriptsize{}$O(sd^{2})$}\tabularnewline
\hline 
{\scriptsize{}Applying $\matM^{-1}$ on a vector} & {\scriptsize{}$O(d^{2})$} & {\scriptsize{}$O(d^{2})$}\tabularnewline
\hline 
{\scriptsize{}Retraction} & {\scriptsize{}$O(ndp+dp^{2})$} & {\scriptsize{}$O(ndp+dp^{2})$}\tabularnewline
\hline 
{\scriptsize{}Inner product on the tangent space (Riemannian metric)} & {\scriptsize{}$O(\nnz{\matS\matZ}p+dp)$} & {\scriptsize{}$O(\nnz{\matS\matZ}p+dp)$}\tabularnewline
\hline 
{\scriptsize{}Orthogonal projections on the tangent/normal space} & {\scriptsize{}$O(ndp+d^{2}p+dp^{2})$} & {\scriptsize{}$O(ndp+d^{2}p+dp^{2})$}\tabularnewline
\hline 
{\scriptsize{}Vector Transport} & {\scriptsize{}$O(ndp+d^{2}p+dp^{2})$} & {\scriptsize{}$O(ndp+d^{2}p+dp^{2})$}\tabularnewline
\hline 
{\scriptsize{}Riemannian gradient} & {\scriptsize{}$O(ndp+d^{2}p+dp^{2}+T_{\nabla\bar{f}})$} & {\scriptsize{}$O(ndp+d^{2}p+dp^{2}+T_{\nabla\bar{f}})$}\tabularnewline
\hline 
{\scriptsize{}Applying the Riemannian Hessian to a tangent vector} & {\scriptsize{}$O(ndp+d^{2}p+\nnz{\matS\matZ}p+dp^{2}+T_{\nabla\bar{f}}+T_{\nabla^{2}\bar{f}})$} & {\scriptsize{}$O(ndp+d^{2}p+\nnz{\matS\matZ}p+dp^{2}+T_{\nabla\bar{f}}+T_{\nabla^{2}\bar{f}})$}\tabularnewline
\hline 
\end{tabular}{\scriptsize\par}
\end{table}

\begin{remark}
In Table \ref{tab:costs} we assumed, for simplicity, that $s\geq d$
. However, if $\lambda$ is sufficiently large, it is possible for
the prescribed values of $s$ ($s\geq20s_{\lambda}(\matZ)^{2}/\delta$)
to be smaller than $d$. In such cases, we can reduce the $O(sd^{2})$
term in the complexity to $O(sd\min(s,d))$ by employing the Woodbury
formula. We omit the details.
\end{remark}

\section{\label{sec:Randomized-preconditioning-for-CCA}Sketched Iterative
CCA}

CCA, originally introduced by Hotelling in 1936~\cite{hotelling1936relations},
is a well-established method in statistical learning with numerous
applications (e.g., \cite{sun2010scalable,kim2007discriminative,su2012discriminant,dhillon2012two,dhillon2011multi,chaudhuri2009multi}). In CCA, the relation between a pair of data sets in matrix form is
analyzed, where the goal is to find the directions of maximal correlation
between a pair of observed variables. In the language of linear algebra,
CCA measures the similarities between two subspaces spanned by the
columns of the two matrices, whereas in the geometric point of
view, CCA computes the cosine of the principal angles between the
two subspaces. We consider a regularized version of CCA defined below\footnote{The definition is formulated as a linear algebra problem. While the
problem can be motivated, and described, in the language of statistics,
the linear algebraic formulation is more convenient for our purposes.}:
\begin{definition}
\label{def:CCA for Stiefel}Let $\matX\in\R^{n\times d_{\x}}$ and
$\matY\in\R^{n\times d_{\y}}$ be two data matrices, and $\lambda_{\x},\lambda_{\y}\geq0$
be two regularization parameter. Let $q=\max\left(\rank{\matX^{\T}\matX+\lambda_{\x}\matI_{d_{\x}}},\rank{\matY^{\T}\matY+\lambda_{\y}\matI_{d_{\y}}}\right)$.
The $(\lambda_{\x},\lambda_{\y})$\emph{-canonical correlations} $\sigma_{1}\geq\dots\geq\sigma_{q}$
and the $(\lambda_{\x},\lambda_{\y})$\emph{-canonical weights} $\u_{1},\dots,\u_{q}$, $\v_{1},\dots,\v_{q}$
are the arguments that maximize $\Trace{\matU^{\T}\matX^{\T}\matY\matV}$
subject to 
\begin{equation*}
\matU^{\T}(\matX^{\T}\matX+\lambda_{\x}\matI_{d_{\x}})\matU  =  \matI_{q},\ \matV^{\T}(\matY^{\T}\matY+\lambda_{\y}\matI_{d_{\y}})\matV  =  \matI_{q},\ \matU^{\T}\matX^{\T}\matY\matV  =  \diag{\sigma_{1},\dots,\sigma_{q}},
\end{equation*}
where $\matU=\left[\begin{array}{ccc}
\u_{1} & \dots & \u_{q}\end{array}\right]\in\R^{d_{\x}\times q}$ and $\matV=\left[\begin{array}{ccc}
\v_{1} & \dots & \v_{q}\end{array}\right]\in\R^{d_{\y}\times q}$.
\end{definition}

\subsection{\label{subsec:Canonical-Correlation-Analysis}CCA as an Optimization Problem on a Product of Generalized Stiefel Manifolds}

We focus on finding $\sigma_{1},...,\sigma_{p},\u_{1},...,\u_{p}$
and $\v_{1},...,\v_{p}$, where $p\leq q$ is a parameter, i.e., on
finding the top $p$-canonical correlations and the corresponding left and right vectors. For convenience, we use
the following notations: 
\[
\Sigma_{\x\x}\coloneqq\matX^{\T}\matX+\lambda_{\x}\matI_{d_{\x}},\quad\Sigma_{\y\y}\coloneqq\matY^{\T}\matY+\lambda_{\y}\matI_{d_{\y}},\quad\Sigma_{\x\y}\coloneqq\matX^{\T}\matY\,.
\]
With these notations, we can reformulate the problem of finding the
top $p$-canonical correlations succinctly in the following way: 
\begin{eqnarray}
\max\Trace{\matU^{\T}\Sigma_{\x\y}\matV},&\st& \matU^{\T}\Sigma_{\x\x}\matU=\matI_{p},\ \matV^{\T}\Sigma_{\y\y}\matV=\matI_{p},\nonumber\\
& & \matU^{\T}\Sigma_{\x\y}\matV\,\text{is diagonal with non-increasing diagonal}
\label{eq:CCAprob_pre}
\end{eqnarray}
Notice that \emph{without} the last constraint, Problem~(\ref{eq:CCAprob_pre})
is a maximization over the product of two generalized Stiefel manifolds.
However, without this constraint the solution is not unique. The reason
for that is that the trace operator and the constraint set are invariant
to multiplication by orthonormal matrices, and so there are optimal
values that are non-diagonal (and so additional steps are needed to
extract the canonical correlations from such values). In order to
circumvent this issue, we use a well known method to modify such problems
so to make the solution unique. The modification is based on the \emph{von
Neumann cost function} \cite{von1962some}. That is, we replace the
objective function $\Trace{\matU^{\T}\Sigma_{\x\y}\matV}$ with $\Trace{\matU^{\T}\Sigma_{\x\y}\matV\matN}$
where $\matN=\diag{\mu_{1},...,\mu_{p}}$ and we take arbitrary $\mu_{1},...,\mu_{p}$
such that $\mu_{1}>...>\mu_{p}>0$. In other words, the problem we
wish solve is 
$$
\max\Trace{\matU^{\T}\Sigma_{\x\y}\matV\matN},\st \matU^{\T}\Sigma_{\x\x}\matU=\matI_{p},\ \matV^{\T}\Sigma_{\y\y}\matV=\matI_{p}.
$$

In the next subsection, we detail the Riemannian components which allow to solve the CCA problem. We show that critical points of the corresponding objective function consist of coordinated left and right canonical correlation vectors not necessarily
on the same phase. In particular, the optimal solutions are critical points consisting of coordinated left and right top $p$-canonical correlation vectors on the same phase. 

\subsection{\label{subsec:Preconditioned-CCA-Algorithm}Preconditioned Riemannian
Components for CCA}

In this subsection we derive the Riemannian components associated
with the CCA problem. The CCA problem is a constraint maximization
on the product of two generalized Stiefel manifolds: $\stiefel_{\Sigma_{\x\x}}(p,d_{\x})$
and $\stiefel_{\Sigma_{\y\y}}(p,d_{\y})$. We denote the search space
by $\elpCCA\coloneqq\stiefel_{\Sigma_{\x\x}}(p,d_{\x})\times\stiefel_{\Sigma_{\y\y}}(p,d_{\y})$.
We consider the use of Riemannian optimization for solving the CCA
problem, while exploiting the geometry of the preconditioned generalized
Stiefel manifold and use the notion of product manifold. To make the
calculations easier, we denote $d=d_{\x}+d_{\y}$ and $\matZ\coloneqq\left[\begin{array}{cc}
\matU^{\T} & \matV^{\T}\end{array}\right]^{\T}\in\R^{d\times p}$ where $\matU\in\stiefel_{\Sigma_{\x\x}}(p,d_{\x})$ and $\matV\in\stiefel_{\Sigma_{\y\y}}(p,d_{\y})$.
Henceforth, we abuse notation and view $\elpCCA$ as a subset of $\R^{d\times p}$
given by this coordinate split, and also write $\matZ=(\matU,\matV)$
as a shorthand for $\matZ=\left[\begin{array}{cc}
\matU^{\T} & \matV^{\T}\end{array}\right]^{\T}$. With these conventions, the optimization problem can be rewritten
in the following way:

\begin{equation}
\min_{\matZ\in\elpCCA}\fcca\ ,\ \fcca\coloneqq-\frac{1}{2}\Trace{\matZ^{\T}\left[\begin{array}{cc}
 & \Sigma_{\x\y}\\
\Sigma_{\x\y}^{\T}
\end{array}\right]\matZ\matN}\label{eq:costriemannianccaBrockett}
\end{equation}

For a product of disjoint Generalized Stiefel manifolds, if the number of columns of the matrices that belong to each of the manifolds is the same,
the various Riemannian components can be computed separately on each
of the manifolds and then stacked on top of each other. We use two preconditioning schemes $\matU\mapsto\matM^{(\x\x)}_{\matU}$ and $\matV\mapsto\matM^{(\y\y)}_{\matV}$ to define Riemannian
metrics on $\stiefel_{\Sigma_{\x\x}}(p,d_{\x})$ and $\stiefel_{\Sigma_{\y\y}}(p,d_{\y})$
respectively, and these, in turn, define a Riemannian metric on
$\elpCCA$ via $\matZ\mapsto\matM_{\matZ}\coloneqq\blockdiag{\matM^{(\x\x)}_{\matU},\matM^{(\y\y)}_{\matV}}$.
For $\matU\in\stiefel_{\Sigma_{\x\x}}(p,d_{\x})$, let $\mat{\Pi}_{\matU}\left(\cdot\right)$
denote the projection on $T_{\matU}\stiefel_{\Sigma_{\x\x}}(p,d_{\x})$,
and similarly for $\mat{\Pi}_{\matV}\left(\cdot\right)$ where $\matV\in\stiefel_{\Sigma_{\y\y}}(p,d_{\y})$.
Similarly, let $\Pi_{\matU}^{\perp}(\cdot)$ and $\Pi_{\matV}^{\perp}(\cdot)$
be the projections on the corresponding normal spaces. Given $\matZ\in\elpCCA$
the orthogonal projection on the tangent space $T_{\matZ}\elpCCA$
is $\mat{\Pi}_{\matZ}\left(\xi_{\matZ}\right)=\left(\mat{\Pi}_{\matU}\left(\xi_{\matU}\right),\mat{\Pi}_{\matV}\left(\xi_{\matV}\right)\right)$
where $\xi_{\matZ}\coloneqq\left(\xi_{\matU},\xi_{\matV}\right)\in\R^{d\times p}$.

Let $\bar{f}_{{\bf CCA}}$ be $\justfcca$ extended smoothly to be
defined on $\R^{d\times p}$ where $\bar{f}_{{\bf CCA}}$ is defined
by Eq.~(\ref{eq:costriemannianccaBrockett}) as well. The following
are analytical expressions for the Riemannian gradient and the Riemannian
Hessian in ambient coordinates:
\begin{equation}
\grad{\fcca}=\mat{\Pi}_{\matZ}\left(\matM^{-1}_{\matZ}\nabla\bar{f}_{{\bf CCA}}(\matZ)\right)=-\left[\begin{array}{c}
\mat{\Pi}_{\matU}\left(\left(\matM^{(\x\x)}_{\matU}\right)^{-1}\Sigma_{\x\y}\matV\matN\right)\\
\mat{\Pi}_{\matV}\left(\left(\matM^{(\y\y)}_{\matV}\right)^{-1}\Sigma_{\x\y}^{\T}\matU\matN\right)
\end{array}\right],\label{eq:gradCCA}
\end{equation} 
\begin{multline}
\hess{\fcca}  =  \mat{\Pi}_{\matZ}\left(\matM^{-1}_{\matZ}\left[-\sigmacca\xi_{\matZ}\matN+\Sigma\left[\begin{array}{c}
\xi_{\matU}\matU^{\T}\Sigma_{\x\y}\matV\matN\\
\xi_{\matV}\matV^{\T}\Sigma_{\x\y}^{\T}\matU\matN
\end{array}\right]+\right.\right. \\
    \quad\quad\quad\quad\quad\quad\left.\left.+\Sigma\left[\begin{array}{cc}
\xi_{\matU}\matU^{\T}\\
  \xi_{\matV}\matV^{\T}
\end{array}\right]\matM_{\matZ}\grad{\fcca}\right]\right),\label{eq:HessCCA}
\end{multline}
$$
\sigmacca\coloneqq\left[\begin{array}{cc}
 & \Sigma_{\x\y}\\
\Sigma_{\x\y}^{\T}
\end{array}\right],\ \Sigma\coloneqq\blockdiag{\Sigma_{\x\x},\Sigma_{\y\y}},
$$
where Eq. \eqref{eq:HessCCA} is valid for critical
points or if $\matM_{\matZ}\coloneqq\matM=\blockdiag{\matM^{(\x\x)},\matM^{(\y\y)}}$.

Along with formulas for the retraction and vector transport, various Riemannian optimization algorithms can be applied to solve Problem~(\ref{eq:costriemannianccaBrockett}). In the next theorem, we summarize the critical points. Note that an important outcome of this theorem is that we can obtain $\sigma_i$ from $\matU^{\T}\Sigma_{\x\y}\matV$ for critical $\matU$ and $\matV$. In the theorem statement, a pair of left and right canonical correlation vectors $\u$ and $\v$ are {\em on the same phase} if $\u^\T \matX^\T\matY\v \geq 0$. The proof is delegated to Appendix \ref{subsec:proofcriticalpointofCCA}.
\begin{theorem}
\label{thm:criticalpointofCCA} A point $\matZ=(\matU,\matV)\in\elpCCA$
is a critical point of $\fcca$ on $\elpCCA$
if and only if the columns of $\matU$ and $\matV$ are left and right coordinated canonical correlation vectors not necessarily on the same phase. 

The optimal solutions of minimizing $\fcca$ on $\elpCCA$ are critical points $\matZ$ such that the columns of $\matU$ and $\matV$ are coordinated left and right top $p$-canonical correlation vectors on the same phase. Moreover, the optimal solution is unique up to sign of the columns of $\matU$ and $\matV$ if $\sigma_{1}>\sigma_{2}>...>\sigma_{p+1}\geq0$. 
\end{theorem}

\subsection{\label{sec:CCAstability}Stability of the Critical Points of $f_{{\bf CCA}}$}

Naturally, we want our proposed optimization algorithm to converge
to an optimal point. In Theorem \ref{thm:criticalpointofCCA} we characterized
all the critical points of $\fcca$ on $\elpCCA$.
In general, a typical guarantee for Riemannian optimization algorithms is
that for a sequence of iterates, all the accumulation points of the
sequence are critical points (e.g., \cite[Theorem 4.3.1]{AMS09}).
 Unfortunately, this guarantee does not specify to which of the critical
points the convergence is to. However, we can utilize the fact that in practice, for a sufficiently close initial guess, Riemannian optimization methods converge to the stable critical points and do not converge to unstable critical points~\cite[Section 4.4]{AMS09}. We analyze the stability of the various critical points of $\fcca$ on $\elpCCA$.  Note that solving CCA on the generalized Stiefel manifold (and its products) can be reduced to an SVD problem and that the discussion of the stability of critical points is independent of the choice of the Riemannian metric. In \cite[Theorem 3.4]{helmke1992singular}, the stability of the global minimum of SVD problem is established assuming distinct singular values. In the context of Riemannian optimization, the SVD problem was discussed by \cite{sato2013riemannian}. However, apart from \cite[Theorem 4.6.3]{AMS09} where the stability of the critical points of finding the extreme eigenvalue, i.e., $p=1$, none of the aforementioned references provides an analysis of the stability of the critical points in the context of Riemannian optimization. Thus, for completeness of the presentation, we present our analysis in this paper and delegate the proofs to the appendix. In Theorem \ref{thm:stabilityCCA}, we show that under reasonable assumptions the critical
points which solve Problem~(\ref{eq:costriemannianccaBrockett})
are asymptotically stable. Moreover, we show that critical points
which are saddle points or local maxima of Problem~(\ref{eq:costriemannianccaBrockett})
are unstable. The proof of Theorem \ref{thm:stabilityCCA} is in Appendix \ref{subsec:proofstabilityCCA}.

To do so, we use \cite[Proposition 4.4.1]{AMS09} and \cite[Proposition 4.4.2]{AMS09}.
First, recall the definition of a descent mapping from~\cite[Chapter 4.4]{AMS09}:
we say that $F$ is a \emph{descent mapping }for a cost function $f$
on ${\cal M}$ if $f(F(\x))\leq f(\x)$ for all $\x\in{\cal M}$.
Now, \cite[Proposition 4.4.1]{AMS09} shows that if we use any Riemannian
algorithm that induces a descent mapping, and for which for every
starting point the series of points generated by the algorithm has
only accumulation points that are critical points, then any critical
point which is not a local minimum with a compact neighborhood where
the cost function achieves the same value for all other critical points
is \emph{unstable}. Additionally, \cite[Proposition 4.4.2]{AMS09}
shows that if the same conditions hold, and the distance on the manifold
between iterations goes to zero as the algorithm approaches a local
minimum, then if this minimum is an isolated critical point, it
is an \emph{asymptotically stable} critical point. 

\begin{theorem}
\label{thm:stabilityCCA}Consider using Riemannian optimization to
minimize $\fcca$ subject to $\matZ\in\elpCCA$, and assume that the
mapping defined by the algorithm is a descent mapping. Assume that 
$\sigma_{1}>\sigma_{2}>...>\sigma_{p+1}\geq0$, then $\matZ$
that minimize $\fcca$ on $\elpCCA$ are asymptotically stable. Furthermore,
critical points which are not a local minimum of Problem~(\ref{eq:costriemannianccaBrockett})
are unstable.
\end{theorem}

\begin{remark}
Consider the case that one or more of the top $p$-canonical correlations is not simple, i.e., some $\sigma_{i}=\sigma_{j}$
for $1\leq i,j\leq p$. Then there is no longer an isolated minimum
of the form $\matZ=(\matU,\matV)\in\elpCCA$ such that the columns
of $\matU$ and $\matV$ are the left and right top $p$-canonical vectors, because permutations and linear combinations of the columns of $\matU$ and $\matV$ associated with $\sigma_{i}$ which maintain the phase
do not change the optimal value. In such case, we can still guarantee
that there exists a neighborhood of the space of the left and right top $p$-canonical correlation spaces for which all the starting points converge. However, note that linear convergence of Riemannian gradient descent is no longer guaranteed, as the Riemannian Hessian is only positive semi-definite, but theorems such as \cite[e.g., Theorem 4.5.6]{AMS09}, require a positive definite Riemannian Hessian. \\
Indeed, any neighborhood of these spaces contains a sublevel
set of $\justfcca$ where the only critical points of $\justfcca$
in this sublevel set belongs to the space of the left
and right top $p$-canonical correlation spaces. Thus, if a Riemannian optimization algorithm
which induces a descent mapping is started with an initial point within
such a sublevel set, and assuming all accumulation points are critical
points of $\justfcca$, then it converges to the space of the left
and right top $p$-canonical correlation spaces.
\end{remark}

In the next theorem, we show that under certain assumptions,
critical points which do not solve Problem~(\ref{eq:costriemannianccaBrockett})
are saddle points or local maximizers. Furthermore, the algorithm is likely to converge to the desired global minimizer since under some assumptions it is
the only local minimizer (up to the signs of the columns) among the
critical points, thus making it the only asymptotically stable critical
point. The proof relies on the proof of Theorem \ref{thm:CCAHessiantheorem}, Thus, its proof appears in Appendix \ref{subsec:proofmorestabilityCCA}, after the proof of Theorem \ref{thm:CCAHessiantheorem} (Appendix \ref{subsec:proofCCAHessiantheorem}).

\begin{theorem}\label{thm:morestabilityCCA}
Consider using Riemannian optimization
to minimize $\fcca$ subject to $\matZ\in\elpCCA$, where we use the
Riemannian metric defined by $\matM_{\matZ}=\blockdiag{\matM^{(\x\x)}_{\matU},\matM^{(\y\y)}_{\matV}}$
where $\matM^{(\x\x)}_{\matU}\in\R^{d_{\x}\times d_{\x}}$ and $\matM^{(\y\y)}_{\matV}\in\R^{d_{\y}\times d_{\y}}$
are given preconditioning schemes. Assume that for all $i=1,...,q$ the values $\sigma_{i}$
are distinct,
and that $\Sigma$ is a SPD matrix. Then the global minimizer of
$\fcca$ subject to $\matZ\in\elpCCA$, denoted by $\matZ^{\star}=(\matU^{\star},\matV^{\star})$, is the only local minimizer of $\fcca$ on $\elpCCA$, up to the signs
of the columns, and it is also strict, and all
other critical points are either saddle points or global maximizers.
Thus, $\matZ^{\star}$ is the only asymptotically stable
critical point, and all other critical points are unstable.
\end{theorem}

\subsection{The Effect of Preconditioning on the Convergence}

The following theorem allows us to reason about the quality of a preconditioner for the CCA problem. In this theorem, we
bound the condition number of the Riemannian Hessian at the optimum
based on how well the preconditioner approximates a specific matrix
($\Sigma$). The proof appears in Appendix \ref{subsec:proofCCAHessiantheorem}.
\begin{theorem}
\label{thm:CCAHessiantheorem} Consider using Riemannian optimization
to minimize $\fcca$ subject to $\matZ\in\elpCCA$, where we use the
Riemannian metric defined by $\matM_{\matZ}=\blockdiag{\matM^{(\x\x)}_{\matU},\matM^{(\y\y)}_{\matV}}$
where $\matM^{(\x\x)}_{\matU}\in\R^{d_{\x}\times d_{\x}}$ and $\matM^{(\y\y)}_{\matV}\in\R^{d_{\y}\times d_{\y}}$
are given preconditioning schemes. Also assume that $\sigma_{1}>\sigma_{2}>...>\sigma_{p+1}\geq0$
and that $\Sigma$ is a SPD matrix. Then at the global minimizer of
$\fcca$ subject to $\matZ\in\elpCCA$, denoted by $\matZ^{\star}=(\matU^{\star},\matV^{\star})$,
the following bound on the condition number of the Riemannian Hessian
at $\matZ^{\star}$ holds
\[
\kappa(\hess{f_{{\bf CCA}}(\matZ^{\star}))\leq}\kappa_{{\bf CCA}}^{\star}\cdot\kappa\left(\Sigma,\matM_{\matZ^{\star}}\right)
\]
where 
\[
\kappa_{{\bf CCA}}^{\star}\coloneqq\frac{\max\left\{ \mu_{1}(\sigma_{1}+\sigma_{p+1}),\frac{1}{2}(\mu_{1}+\mu_{2})(\sigma_{1}+\sigma_{2})\right\} }{\min\left\{ \mu_{p}(\sigma_{p}-\sigma_{p+1}),\min_{1\leq j<p}\frac{1}{2}\left(\mu_{j}-\mu_{j+1}\right)\left(\sigma_{j}-\sigma_{j+1}\right)\right\} }
\]
and $\mu_{1}>\dots>\mu_{p}>0$.
If $\matM_{\matZ^{\star}}=\Sigma$ then $\kappa(\hess{f_{{\bf CCA}}(\matZ^{\star}))=}\kappa_{{\bf CCA}}^{\star}$.
\end{theorem}

Note that from Theorem \ref{thm:CCAHessiantheorem}, we have a cumbersome connection between $\matN$ and the condition number. Nevertheless, from the expression for $\kappa_{{\bf CCA}}^{\star}$, one can conclude that to balance the numerator and denominator the differences between the values on the diagonal of $\matN$ should not be too small nor too large. 

The condition number for the case $p=1$ does not involve the values
in $\matN$, and thus is simple and illuminating:
\begin{corollary}\label{cor:CCA_p=1}
For $p=1$, the condition number of the Riemannian
Hessian at the optimum is at most $\frac{\sigma_{1}+\sigma_{2}}{\sigma_{1}-\sigma_{2}}\cdot\kappa(\Sigma,\matM_{\z^{\star}})$. If $\matM_{\z^{\star}}=\Sigma$ then the condition number of the Riemannian
Hessian at the optimum is $\frac{\sigma_{1}+\sigma_{2}}{\sigma_{1}-\sigma_{2}}$.
\end{corollary}

The condition number bound from Corollary \ref{cor:CCA_p=1} decomposes into two components: the first
is the relative eigengap ($(\sigma_{1}+\sigma_{2})/(\sigma_{1}-\sigma_{2})$),
which forms a natural condition number for the problem (if the first
and second correlations are very close, it is very hard to distinguish
between them) that almost always appear in problems of this form,
and a second component which measures how close the preconditioner-defined
metric approximates the natural metric for the constraints. The optimal
preconditioner, according to the bound, is $\matM=\Sigma$. However,
using this preconditioner requires explicitly computing it in $O(nd^{2})$
time. This is too expensive since the exact correlations can be computed
analytically in $O(nd^{2})$ time as well~\cite{bjorck1973numerical}.

\subsection{\label{subsec:Randomized-CCA}Randomized Preconditioning for CCA}

The condition number bound in Theorem \ref{thm:CCAHessiantheorem}
separates two factors: $\kappa_{{\bf CCA}}^{\star}$ and $\kappa\left(\Sigma,\matM\right)$.
The first, $\kappa_{{\bf CCA}}^{\star}$, depends on the gap between
the $p+1$ largest canonical correlations and on the differences between
the values in $\matN$, which are parameters for the CCA problem as
a Riemannian optimization problem. The dependence on the gap almost
always appears in problems of this form, since the more the singular
values are distinct it is easier to distinguish between them. The
second component, $\kappa\left(\Sigma,\matM_{\matZ^{\star}}\right)$, measures how
close the preconditioner, which defines the Riemannian metric, approximates
$\Sigma$. A preconditioner that minimizes the bound in Theorem~\ref{thm:CCAHessiantheorem} is such that $\matM_{\matZ^{\star}}=\Sigma$.
However, using that preconditioner, requires
explicitly computing a factorization of $\Sigma$ which classically requires $\Omega(nd^{2})$
arithmetic operations which is non-beneficial in light of direct solution
methods.

Thus, Theorem \ref{thm:CCAHessiantheorem} provides an argument in
favor of our proposed randomized preconditioner, i.e., easy to factorize
$\matM_{\matZ}$ such that $\kappa\left(\Sigma,\matM_{\matZ^{\star}}\right)$ is bounded.
In order to achieve this goal, the preconditioning schemes $\matU\mapsto\matM^{(\x\x)}_{\matU}$
and $\matV\mapsto\matM^{(\y\y)}_{\matV}$ should approximate $\Sigma_{\x\x}=\matX^{\T}\matX+\lambda_{\x}\matI_{d_{\x}}$
and $\Sigma_{\y\y}=\matY^{\T}\matY+\lambda_{\y}\matI_{d_{\y}}$ respectively (at least at the optimum). Thus, 
as described in Section \ref{sec:rand-precond} we propose constant randomized preconditioning scheme $\matM\coloneqq\blockdiag{\matM^{(\x\x)}, \matM^{(\y\y)}}$, i.e., using sketching
for $\matX$ and $\matY$ correspondingly to define $\matM^{(\x\x)}_{\matU}\coloneqq\matM^{(\x\x)}$
and $\matM^{(\y\y)}_{\matV}\coloneqq\matM^{(\y\y)}$. A pseudocode description
of an end-to-end randomized preconditioned CCA algorithm with warm-start
appears in Algorithm \ref{alg:CCAopt}. The following corollary summarizes
our theoretical results regarding the proposed algorithm. We remark
that \noun{CountSketch} can possibly be replaced with other sketching
transforms (such as SRHT), and Riemannian CG can be replaced with
other Riemannian optimization methods, although the bounds in the
corollary might change. We prove Corollary \ref{cor:CCA} in Appendix \ref{subsec:proofofcorCCA}.

\begin{algorithm}[t]
\caption{\label{alg:CCAopt}Sketched Riemannian Iterative CCA with warm-start.}

\begin{algorithmic}[1]

\STATE\textbf{ Input: }$\matX\in\R^{n\times d_{\x}}$, $\matY\in\R^{n\times d_{\y}}$,
$s\geq\max(d_{\x},d_{\y})$, $\lambda_{\x},\lambda_{\y}\geq0$.

\STATE \textbf{Generate random }$h:\{1,\dots,d\}\to\{1,\dots,s\}$
and $g:\{1,\dots,d\}\to\{-1,+1\}$. Let $\matS$ denote the corresponding
\noun{CountSketch} matrix.

\STATE $\matX_{\matS}\gets\matS\matX$, $\matY_{\matS}\gets\matS\matY$.

\STATE $\tilde{\matU},\tilde{\matV}\gets\mathrm{exact-cca}(\matX_{\matS},\matY_{\matS})$.

\STATE $\matM^{(\x\x)}\gets\matX_{\matS}^{\T}\matX_{\matS}+\lambda_{\x}\matI_{d_{\x}}$,
$\matM^{(\y\y)}\gets\matY_{\matS}^{\T}\matY_{\matS}+\lambda_{\y}\matI_{d_{\y}}$.

\STATE \textbf{Notation: }$\Sigma_{\x\x}=\matX^{\T}\matX+\lambda_{\x}\matI$,
$\Sigma_{\y\y}=\matY^{\T}\matY+\lambda_{\y}\matI$. Do not compute
these matrices (algorithms only require taking products with them).

\STATE \textbf{Choose: }any $\matN=\diag{\mu_{1},...,\mu_{p}}$ s.t.
$\mu_{1}>...>\mu_{p}>0$.

\STATE Using Riemannian CG, solve $\max\Trace{\matU^{\T}\Sigma_{\x\y}\matV\matN}$
s.t. $\matU\in\stiefel_{\Sigma_{\x\x}}(p,d_{\x})$, $\matV\in\stiefel_{\Sigma_{\y\y}}(p,d_{\y})$.
Use $\matM^{(\x\x)}$ and $\matM^{(\y\y)}$ for the metric. Start
the iteration from $\qfm{\tilde{\matU}}{\Sigma_{\x\x}}$ and $\qfm{\tilde{\matV}}{\Sigma_{\y\y}}$.

\end{algorithmic}
\end{algorithm}

\begin{corollary}
\label{cor:CCA}Consider Algorithm \ref{alg:CCAopt}. Let $\delta\in(0,1)$
and denote $s_{\lambda}=\max(s_{\lambda_{\x}}(\matX),s_{\lambda_{\y}}(\matY))$.
If $s=\max(\left\lceil 40s_{\lambda}^{2}/\delta\right\rceil ,d)$,
then with probability of at least $1-\delta$, the
condition number of the Riemannian Hessian at the optimum is bounded
by $3\kappa_{{\bf CCA}}^{\star}$, regardless of the condition number
of $\Sigma_{\x\x}$ and $\Sigma_{\y\y}$. Furthermore, assuming we use Riemannian CG, $n\geq d\geq p$,
and all computations are done in ambient $\R^{d\times p}$ coordinates,
then the preprocessing steps take $O(\nnz{\matX}+\nnz{\matY})=O(nd)$
and $O(sd^{2})$. Assuming a bounded number of line-search steps
in each iteration then each iteration takes $O\left(p\left(\nnz{\matX}+\nnz{\matY}\right)+dp^{2}+d^{2}p\right)$ operations.
\end{corollary}

\section{\label{sec:Randomized-preconditioning-for-LDA}Sketched Iterative
FDA}

Fisher's linear discriminant analysis (FDA), introduced in \cite{fisher1936use}, is a well-known method for classification
\cite{fisher1936use,mika1999fisher}, and more commonly for dimensionality
reduction before classification \cite{chen2012communications}. The
latter is achieved by finding an embedding such that simultaneously
the between-class scatter is maximized and the within-class scatter
is minimized. In this paper, we consider a regularized version of FDA as defined
below:
\begin{definition}
\label{def:LDA multiclass-1}\cite[Section 4.11]{duda2001pattern}
Let $\x_{1}^{(1)},...,\x_{n_{1}}^{(1)},\x_{1}^{(2)},...,\x_{n_{2}}^{(2)}...,\x_{1}^{(2)},...,\x_{n_{l}}^{(l)}\in\R^{d}$
be samples from $l\leq d$ different classes, and denote by $\x_{1},...,\x_{n}$
the union of the different classes (the entire data set in a sequential
index). For $i=1,\dots,n$, let $y_{i}$ denote the label corresponding
to $\x_{i}$, i.e., $y_{i}=k$ if $\x_{i}=\x_{j}^{(k)}$ for some
$j$. Let $\mathbf{m}_{k}$, for $k=1,\dots,l$, denote the sample
mean of class $k$ (i.e., $\mathbf{m}{}_{k}\coloneqq n_{k}^{-1}\sum_{i=1}^{n_{k}}\x_{i}^{(k)}$),
and $\mathbf{m}\coloneqq n^{-1}\sum_{i=1}^{n}\x_{i}=n^{-1}\sum_{k=1}^{l}n_{k}\m_{k}$
denote the data set sample mean of the entire data set. Let $\matS_{\matB}$
and $\matS_{\w}$ be the \emph{between-class}\footnote{For $l=2$ the matrix $\matS_{\matB}$ is defined by $\mat{S_{B}}\coloneqq(\mathbf{m}_{1}-\mathbf{m}_{2})(\mathbf{m}_{1}-\mathbf{m}_{2})^{\T}$.
The definitions coincide after multiplying $\mat{S_{B}}=(\mathbf{m}_{1}-\mathbf{m}_{2})(\mathbf{m}_{1}-\mathbf{m}_{2})^{\T}$
by $2(n_{1}n_{2})/n$. } and \emph{within-class} scatter matrices (respectively): 
\[
\mat{S_{B}}\coloneqq\sum_{k=1}^{l}n_{k}(\mathbf{m}_{k}-\mathbf{m})(\mathbf{m}_{k}-\mathbf{m})^{\T}\quad\mat{S_{w}}\coloneqq\sum_{i=1}^{n}(\x_{i}-\mathbf{m}_{y_{i}})(\x_{i}-\mathbf{m}_{y_{i}})^{\T}
\]
Let $\lambda\geq0$ be a regularization parameter. The $l-1$ FDA weight vectors $\w_{1},...,\mathbf{w}_{l-1}$ are the columns of $\mat W\in\R^{d\times(l-1)}$
such that $\matW$ is the maximizer of the following cost function:
\begin{equation}
J(\mat W)=\frac{\det(\mat W^{\T}\mat{S_{B}}\mat W)}{\det(\mat W^{\T}(\ldaB)\mat W)}.\label{eq:LDAdefcost}
\end{equation}
\end{definition}

\subsection{FDA as an Optimization Problem on a Generalized Stiefel Manifold}

It is well known that the solution of maximizing Eq. (\ref{eq:LDAdefcost}) (i.e., finding the FDA weight vectors) is equivalent to finding a matrix $\matW$ such that its columns are the leading $l-1$
generalized eigenvectors of the matrix pencil $(\mat{S_{B}},\ldaB)$
\cite[Section 4.11]{duda2001pattern}. Note that this generalized
eigenproblem has at most $l-1$ nonzero generalized eigenvalues since
the matrix $\ldaB$ is a SPD matrix, and $\mat{S_{B}}$ is the sum
of $l$ matrices of rank one or less, where only $l-1$ of these are
independent, thus, $\mat{S_{B}}$ is of rank $l-1$ or less. We denote
the eigenvalues of the matrix pencil $\left(\mat{S_{B}},\mat{S_{w}}+\lambda\matI_{d}\right)$ with correspondence to the FDA weight vectors
by $\rho_{1}\geq\rho_{2}\geq...\geq\rho_{d}\geq0$.

We focus on finding the $p$ leading FDA weight vectors, i.e., $\w_{1},...,\mathbf{w}_{p}$ corresponding to $\rho_{1}\geq\rho_{2}\geq...\geq\rho_{p}$ where $p\leq l-1$.
For the purpose of describing and analyzing our algorithm, it is useful
to write,
\[
\mat{S_{w}}=\hat{\matX}^{\T}\hat{\matX},\quad\hat{\matX}\coloneqq\matX-\matY,\quad\mat{S_{B}}=\hat{\matY}^{\T}\hat{\matY}\ ,
\]
where $\hat{\matX}\in\R^{n\times d}$ is a matrix such that each $i$-th
row of $\hat{\matX}$ is $(\x_{i}-\mathbf{m}_{y_{i}})^{\T}$, $\mat X\in\R^{n\times d}$
is a matrix such that each $i$-th row of $\matX$ is $\x_{i}^{\T}$,
$\matY\in\R^{n\times d}$ is a matrix such that each $i$-th row is
of the form $\mathbf{m}_{y_{i}}$ (thus, there are at most $l$ different
rows in $\matY$), and $\hat{\matY}\in\R^{l\times d}$ is a matrix
such that each $k$th row of $\hat{\matY}$ is $\sqrt{n_{k}}(\mathbf{m}_{k}-\mathbf{m})^{\T}$.
With these notations, we can reformulate the problem of finding the
$p$ leading FDA weight vectors as a Riemannian optimization problem
on the generalized Stiefel manifold (see, \cite[Section 10.2, Eq. (10.5)]{fukunaga2013introduction}),
i.e., finding the generalized eigenvalues of the matrix pencil $(\mat{S_{B}},\ldaB)$.
We use the \emph{Brockett cost function} \cite{brockett1991dynamical}
and obtain the following optimization problem 
\begin{equation}
\max_{\matW\in\R^{d\times p}}\Trace{\mat W^{\T}\hat{\matY}^{\T}\hat{\matY}\mat W\matN},\st
\mat W^{\T}(\hat{\matX}^{\T}\hat{\matX}+\lambda\matI_{d})\mat W=\matI_{p},
\label{eq:LDAprob}
\end{equation}
where $\matN=\diag{\mu_{1},...,\mu_{p}}$ where we take arbitrary $\mu_{1},...,\mu_{p}$
such that $\mu_{1}>...>\mu_{p}>0$.

\begin{remark}
Problem (\ref{eq:LDAprob}) is actually a relaxation of the \emph{trace-ratio
} problem \cite{wang2007trace,ngo2012trace}, which provides superior
results in various tasks (e.g., classification and clustering \cite{wang2007trace}).
We focus on Problem (\ref{eq:LDAprob}) since it is amenable to our
preconditioning strategy.
\end{remark}

In the next subsection, we detail the Riemannian components which allow to solve the FDA problem. We show that critical points of the corresponding objective function are matrices $\matW$ such that the columns are some $p$ FDA weight vectors. In particular, the optimal solutions are critical points consisting of the $p$ leading FDA weight vectors. We further show that if $\rho_{1},\dots,\rho_{p+1}$ are
distinct, the optimal solution is unique up to sign of the columns
of $\matW$.

\subsection{\label{subsec:Preconditioned-LDA-Algorithm}Preconditioned Riemannian
Components for FDA}

We first transform Problem (\ref{eq:LDAprob}) into a minimization
problem: 
\begin{equation}
\min_{\Stiefellda}\flda\ ,\ \flda\coloneqq-\frac{1}{2}\Trace{\mat W^{\T}\mat{S_{B}}\mat W\matN}.\label{eq:lda_riemannian}
\end{equation}
We use the components of the generalized Stiefel manifold with a Riemannian metric defined by a preconditioning scheme $\matM_{\matW}\in\R^{d\times d}$ to apply Riemannian optimization to solve Problem~(\ref{eq:lda_riemannian}).

Let $\bar{f}_{{\bf FDA}}$ be $\justflda$ extended smoothly to be
defined on $\R^{d\times p}$ where $\bar{f}_{{\bf FDA}}$ is defined
by Eq.~(\ref{eq:lda_riemannian}) as well. 

For $\matW\in\Stiefellda$, let $\mat{\Pi}_{\matW}\left(\cdot\right)$ denote the projection on
$T_{\matW}\Stiefellda$. Similarly, let $\mat{\Pi}_{\matW}^{\perp}\left(\cdot\right)$
be the projection on the corresponding normal space. The following
are analytical expressions for the Riemannian gradient and 
Hessian in ambient coordinates: 
\begin{equation}
\grad{\flda}=\mat{\Pi}_{\matW}\left(\matM^{-1}_{\matW}\nabla\bar{f}_{{\bf FDA}}(\mat W)\right)-\mat{\Pi}_{\matW}\left(\matM^{-1}_{\matW}\mat{S_{B}}\mat W\matN\right),\label{eq:LDAgrad}
\end{equation}
\begin{multline}
\hess{\flda}[\xi_{\mat W}] = \mat{\Pi}_{\matW}(\matM^{-1}_{\matW}[-\mat{S_{B}}\xi_{\mat W}\matN+\\(\ldaB)\xi_{\mat W}(\mat W^{\T}\mat{S_{B}}\mat W\matN+\grad{f(\mat W)})]),\label{eq:LDAhess}
\end{multline}
where Eq. \eqref{eq:LDAhess} is valid for critical
points or if $\matM_{\matW}\coloneqq\matM$.

Along with formulas for the retraction and vector transport, various Riemannian optimization algorithms can be applied to solve Problem~(\ref{eq:lda_riemannian}). In the next theorem, we summarize the critical points. Note that an important outcome of this theorem is that we can obtain $\rho_i$ from $\matW^{\T}\matSb \matW$ with a critical $\matW$. The proof is almost the same as the proof of Theorem \ref{thm:criticalpointofCCA}, and presented in Appendix \ref{subsec:proofcriticalpointsofLDA}.

\begin{theorem}\label{thm:criticalpointsofLDA} A point $\matW\in\Stiefellda$ is
a critical point of $\flda$ on the manifold $\Stiefellda$ if and only
if the columns of $\matW$ are some $p$ FDA weight vectors.

The optimal solutions of minimizing $\flda$ on $\Stiefellda$ are critical points $\matW$ such that the columns are the $p$ leading FDA weight vectors. Moreover, the optimal solution is unique up to sign of the columns if $\rho_{1}>\rho_{2}>...>\rho_{p+1}\geq0$.
\end{theorem}

\subsection{\label{sec:LDAstability}Stability of the Critical Points of $f_{{\bf FDA}}$}

In the FDA problem, we want any optimization algorithm to converge
to an optimal point. In Theorem \ref{thm:criticalpointsofLDA} we
prove that the critical points of $\flda$ on $\Stiefellda$
are matrices $\matW\in\Stiefellda$ such that the columns are some $p$
FDA weight vectors. Similarly to the CCA problem, we analyze the stability of the critical points using \cite[Proposition 4.4.1]{AMS09} and \cite[Proposition 4.4.2]{AMS09}. Note that similarly to the CCA problem, solving FDA on the generalized Stiefel manifold can be reduced to a generalized eigenvalues problem, and that the discussion of the stability of critical points is independent of the choice of Riemannian metric. In \cite[Theorem 4]{brockett1991dynamical}, the stability of the global minimum of the generalized eigenvalue problem is established assuming distinct eigenvalues. In the context of Riemannian optimization, the eigenvalue problem was discussed in various references, e.g., \cite[Sections 4.6 and 4.8]{AMS09}. However, apart from \cite[Theorem 4.6.3]{AMS09} where the stability of the critical points of finding the extreme eigenvalue, i.e., $p=1$, none of the aforementioned references provides a complete analysis of the stability of the critical points in the context of Riemannian optimization. Thus, for completeness of the presentation, we present our analysis in this paper and relegate the proofs to the appendix. 

The following theorem shows that under reasonable assumptions the
critical points which solve Problem~(\ref{eq:lda_riemannian}) are
asymptotically stable and critical points which are saddle
points of Problem~(\ref{eq:lda_riemannian}) are unstable. 
The proof is analogous to the proof of Theorem \ref{thm:stabilityCCA}, and presented in Appendix \ref{subsec:proofstabilityLDA}.

\begin{theorem}
\label{thm:stabilityLDA}Consider using Riemannian optimization to
minimize $\flda$ subject to $\matW\in\Stiefellda$, and assume that
the mapping defined by the algorithm is a descent mapping. Assume
that $\rho_{1}>\rho_{2}>...>\rho_{p+1}\geq0$, then $\matW\in\Stiefellda$ such that the columns
are the $p$ leading FDA weight vectors are asymptotically stable. Furthermore, critical
points which are not a local minimum of Problem~(\ref{eq:lda_riemannian})
are unstable.
\end{theorem}

\begin{remark}
Consider the case that one or more of the $p$-dominant generalized
eigenvalues of the matrix pencil $(\mat{S_{B}},\ldaB)$ is not simple,
i.e., some $\rho_{i}=\rho_{j}$ for $1\leq i,j\leq p$. Then there
is no longer an isolated minimum of the form $\matW\in\Stiefellda$
such that the columns are the $p$ leading FDA weight vectors, because permutations and linear combinations of
the columns of $\matW$ associated with $\rho_{i}$ do not change
the optimal value. In such case, we can still guarantee that there
exists a neighborhood of the space of the $p$-dominant generalized
eigenspaces for which all the starting points converge to the space
of the $p$-dominant generalized eigenspaces. However, note that linear convergence of Riemannian gradient descent is no longer guaranteed, as the Riemannian Hessian is only positive semi-definite, but theorems such as \cite[Theorem 4.5.6]{AMS09}, require a positive definite Riemannian Hessian. \\
Indeed, any neighborhood
of these generalized eigenspaces contains a sublevel set of $\justflda$ where the only critical points of $\justflda$ in this sublevel
set belong to the space of the $p$-dominant generalized eigenspaces.
Thus, if a Riemannian optimization algorithm which induces a descent
mapping is started with an initial point within such a sublevel set,
and assuming all accumulation points are critical points of $\justflda$,
then it converges to the space of the $p$-dominant generalized eigenspaces.
\end{remark}
In the next theorem, we show that under certain assumptions,
critical points which do not solve Problem~(\ref{eq:LDAprob})
are saddle points or local maximizers. Furthermore, the algorithm is likely to converge to the desired global minimizer since under some assumptions it is
the only local minimizer (up to the signs of the columns) among the
critical points, thus making it the only asymptotically stable critical
point. The proof relies on the proof of Theorem \ref{thm:LDAHessiantheorem}, and is presented in Appendix \ref{subsec:proofmorestabilityLDA}.

\begin{theorem}\label{thm:morestabilityLDA}
Consider using Riemannian optimization
to minimize $\flda$
subject to $\mat W\in\Stiefellda$, where we use the Riemannian metric
defined by $\matM_{\matW}\in\R^{d\times d}$, which is a given preconditioning scheme. Assume that for all $i=1,...,d$ the values $\rho_{i}$
are distinct and that
$\ldaB$ is a SPD matrix. Let $\matW^{\star}$ denote the global
minimizer of $\flda$ subject to $\matW\in\Stiefellda$. Then $\matW^{\star}$ is the only local minimizer of $\flda$ on $\Stiefellda$, up to the
signs of the columns, and it is also strict, and all other critical
points are either saddle points or strict local maximizers. Thus,
$\matW^{\star}$ is the only asymptotically stable critical
point, and all other critical points are unstable.
\end{theorem}

\subsection{The Effect of Preconditioning on the Convergence}

The following theorem provides a bound on the condition number of the Riemannian
Hessian at the optimum based on how well the preconditioner approximates
a specific matrix ($\ldaB$). This theorem provides a general guideline
in designing a Riemannian preconditioner via the Riemannian metric,
and in particular motivates our proposed approach detailed in Subsection
\ref{subsec:Randomized-LDA}. The proof uses the same arguments as
the proof of Theorem \ref{thm:CCAHessiantheorem}, though they are
simplified considerably since only a single generalized Stiefel manifold
is considered. The proof appears in
Appendix \ref{subsec:proofLDAHessiantheorem}.
\begin{theorem}
\label{thm:LDAHessiantheorem} Consider using Riemannian optimization
to minimize $\flda$
subject to $\mat W\in\Stiefellda$, where we use the Riemannian metric
defined by $\matM_{\matW}\in\R^{d\times d}$, which is a given preconditioning scheme. Assume that $\rho_{1}>\rho_{2}>...>\rho_{p+1}\geq0$ and that
$\ldaB$ is an SPD matrix. Let $\matW^{\star}$ denote the global
minimizer of $\flda$ subject to $\matW\in\Stiefellda$. Then, 
\[
\kappa(\hess{\justflda(\matW^{\star}))}\leq\kappa_{{\bf FDA}}^{\star}\cdot\kappa\left(\ldaB,\matM_{\matW^{\star}}\right)
\]
where 
\[
\kappa_{{\bf FDA}}^{\star}\coloneqq\frac{\mu_{1}\left(\rho_{1}-\rho_{d}\right)}{\min\left\{ \mu_{p}\left(\rho_{p}-\rho_{p+1}\right),\min_{1\leq j<p}\frac{1}{2}\left(\mu_{j}-\mu_{j+1}\right)\left(\rho_{j}-\rho_{j+1}\right)\right\} }
\]
and $\mu_{1}>\dots>\mu_{p}>0$.
If $\matM_{\matW^{\star}}=\ldaB$ then $\kappa(\hess{\justflda(\matW^{\star}))}=\kappa_{{\bf FDA}}^{\star}$.
\end{theorem}

Note that from Theorem \ref{thm:LDAHessiantheorem}, we have a cumbersome connection between $\matN$ and the condition number. Nevertheless, from the expression for $\kappa_{{\bf FDA}}^{\star}$, one can conclude that to balance the numerator and denominator the differences between the values on the diagonal of $\matN$ should not be too small nor too large. 

The condition number for the case $p=1$ does not involve the values
in $\matN$, and thus is simple and illuminating:
\begin{corollary}\label{cor:FDA_p=1}
For $p=1$, the condition number of the Riemannian
Hessian at the optimum is at most $\frac{\rho_{1}-\rho_{d}}{\rho_{1}-\rho_{2}}\cdot\kappa\left(\ldaB,\matM_{\matW^{\star}}\right)$. If $\matM_{\matW^{\star}}=\ldaB$ then the condition number of the Riemannian
Hessian at the optimum is $\frac{\rho_{1}-\rho_{d}}{\rho_{1}-\rho_{2}}$.
\end{corollary}

The condition number bound from Corollary \ref{cor:FDA_p=1} decomposes to two components: the first
is the relative eigengap ($\rho_{1}/(\rho_{1}-\rho_{2})$) (when $\rho_{d}$ is very small or zero), which
form a natural condition number for the problem (if the first and
second generalized eigenvalues are very close, it is very hard to
distinguish between invariant subspaces corresponding to them) that
almost always appears in problems of this form, and a second component
which measures how close the preconditioner-defined metric approximates
the natural metric for the constraints. The optimal preconditioner
according to the bound, is $\matM=\mat{S_{w}}+\lambda\matI_{d}$.
However, using this preconditioner requires explicitly computing
it in $O(nd^{2})$ time. This is too expensive, since the exact discriminant
variables can be computed analytically in $O(nd^{2})$ time as well
(exact solution requires finding the inverse of $\mat{S_{w}}+\lambda\matI_{d}$
and then finding the first eigenvalue and corresponding eigenvector
of $(\mat{S_{w}}+\lambda\matI_{d})^{-1}\mat{S_{B}}$).

\subsection{\label{subsec:Randomized-LDA}Randomized Preconditioning for FDA}

The bound on the condition number in Theorem \ref{thm:LDAHessiantheorem} decomposes two components: the first, $\kappa_{{\bf FDA}}^{\star}$,
depends only on the FDA problem (and its formulation as a Riemannian
optimization problem). $\kappa_{{\bf FDA}}^{\star}$ depends on both
the gap between the $p+1$ largest generalized eigenvalues, and on
the diagonal elements of $\matN$, which are parameters
of the optimization problem. The dependence on the gap almost always
appears in problems of this form, since the more the generalized eigenvalues associated with the FDA weight vectors
we search for are distinct it is easier to distinguish between them.
The second component, $\kappa\left(\ldaB,\matM_{\matW^{\star}}\right)$, measures
how close the preconditioner, which defines the Riemannian metric, approximates $\ldaB$. A preconditioner that minimizes the
bound in Theorem \ref{thm:LDAHessiantheorem} is such that $\matM_{\matW^{\star}}=\ldaB$. However, using that preconditioner requires explicitly
computing $\ldaB$ which takes $O(nd^{2})$ arithmetic operations.
As for the CCA problem, direct methods for solving FDA require 
$\Theta(nd^{2})$ arithmetic operations as well (exact solution requires
finding the inverse of $\ldaB$ and then finding the eigenvalues and
corresponding eigenvectors of $(\ldaB)^{-1}\mat{S_{B}}$). Thus, we
want $\matM_{\matW}$ to approximate $\ldaB$ at the optimum, while allowing a cheap factorization
which is satisfied by our proposed randomized preconditioning approach.

We propose to design the preconditioner in the following way: $\matM_{\matW}\coloneqq\matM$
approximates $\ldaB=\hat{\matX}^{\T}\hat{\matX}+\lambda\matI_{d}$
via a matrix sketching procedure for $\hat{\matX}$ as described in
Section \ref{sec:rand-precond}. A full description of a randomized
preconditioned algorithm for FDA with warm-start appears in Algorithm
\ref{alg:LDAopt}. The following corollary summarize our theoretical
results regarding the proposed algorithm. Note that \noun{CountSketch}
can possibly be replaced with other sketching transforms (such as SRHT), and Riemannian
CG can be replaced with any other Riemannian optimization methods,
although the bound in the corollary might change. The proof is in Appendix \ref{subsec:proofofcorLDA}.

\begin{algorithm}[t]
\caption{\label{alg:LDAopt}Sketched Riemannian Iterative FDA with warm-start.}

\begin{algorithmic}[1]

\STATE\textbf{ Input: }$\matX\in\R^{n\times d}$, $\y\in\mathbb{N}^{n}$,
$s\geq d$, $\lambda\geq0$.

\STATE Compute matrices $\matS_{\matB}$ and $\hat{\matX}$.

\STATE \textbf{Generate random }$h:\{1,\dots,d\}\to\{1,\dots,s\}$
and $g:\{1,\dots,d\}\to\{-1,+1\}$. Let $\matS$ denote the corresponding
\noun{CountSketch} matrix.

\STATE $\matX_{\matS}\gets\matS\hat{\matX}$.

\STATE $\tilde{\matW}\gets\mathrm{exact-fda}(\matX_{\matS})$.

\STATE $\matM\gets\matX_{\matS}^{\T}\matX_{\matS}+\lambda\matI_{d}$.

\STATE \textbf{Notation: }$\matS_{\w}=\hat{\matX}^{\T}\hat{\matX}$.
Do not compute this matrix (algorithms only require taking products
with it).

\STATE \textbf{Choose: }any $\matN=\diag{\mu_{1},...,\mu_{p}}$ s.t.
$\mu_{1}>...>\mu_{p}>0$.

\STATE Using Riemannian CG, solve $\max\text{\ensuremath{\Trace{\mat W^{\T}\mat{S_{B}}\mat W\matN}}}$
s.t. $\matW\in\stiefel_{\mat{S_{w}}+\lambda\matI_{d}}(p,d)$. Use
$\matM$ for the metric. Start the iteration from $\qfm{\tilde{\matW}}{\matS_{\w}+\lambda\matI}$.

\end{algorithmic}
\end{algorithm}

\begin{corollary}
\label{cor:LDA}Consider Algorithm \ref{alg:LDAopt}. Let $\delta\in(0,1)$. If $s=\max(\left\lceil 20s_{\lambda}(\hat{\matX})^{2}/\delta\right\rceil ,d)$,
then with probability of at least $1-\delta$, the
condition number of the Riemannian Hessian at the optimum is bounded by $3\kappa_{{\bf FDA}}^{\star}$,
regardless of the condition number of $\matS_{\w}+\lambda\matI$.
Furthermore, assuming we use Riemannian CG, $n\geq d\geq p$, and all computations are
done in ambient $\R^{d\times p}$ coordinates, then the preprocessing
steps take $O(\nnz{\hat{\matX}})=O(nd)$ and $O(sd^{2})$. Assuming
a bounded number of line-search steps in each iteration then each
iteration takes $O(p\left(\nnz{\matX}+ld\right)+\nnz{\hat{\matX}}p+d^{2}p+dp^{2})$ operations.
\end{corollary}

\section{\label{sec:Numerical-Experiments}Numerical Experiments}
In the following section, we present our numerical experiments illustrating our randomized preconditioning approach.

\subsection{Experiments with Real-World Data Sets}

We report experiments with our proposed preconditioned Riemannian
optimization algorithms. The experiments are not designed to be exhaustive;
we use a prototype \noun{MATLAB} implementation. In particular we
present experiments with the preconditioned CCA and FDA algorithms
presented in Sections \ref{sec:Randomized-preconditioning-for-CCA}
and \ref{sec:Randomized-preconditioning-for-LDA}. Our aim is to assess the effectiveness
of our randomized preconditioning approach.

In addition to Algorithms \ref{alg:CCAopt} and \ref{alg:LDAopt},
we experiment with an additional preconditioning strategy based on \cite{GOS16}, which we term as \emph{Dominant Subspace Preconditioning}. This preconditioner was designed via an approximation of
the empirical correlation matrix to speed
up SVRG when solving
ridge regression problems. In our experiments, we use this method to
approximate $\Sigma_{\x\x}$ and $\Sigma_{\y\y}$ for CCA, and $\matS_{\w}+\lambda\matI_{d}$
for FDA. Specifically, the matrices $\matX,\matY$ and $\hat{\matX}$ are approximated in the following way: suppose $\matA=\matX^{\T}\matX\in\R^{d\times d}$, and let $\matX=\matU {\Lambda}^{\nicehalf}\matV^{\T}$ be an SVD decomposition of $\matX$  such that $\matA=\matU\Lambda\matU^{\T}$
is an eigendecomposition, with the diagonal entries in $\Lambda$
sorted in a descending order. Given $k$, let us denote by $\matU_{k}$ the
first $k$ columns of $\matU$, $\Lambda_{k}$ denote the leading
$k\times k$ minor of $\Lambda$, and $\lambda_{k}$ the $k$-th largest
eigenvalue of $\matA$. Then, the $k$-dominant subspace preconditioner
of $\matA+\lambda\matI_{d}$ is $\matU_{k}(\Lambda_{k}-\lambda_{k}\matI)\matU_{k}^{\T}+(\lambda_{k}+\lambda)\matI_{d}$.
The dominant subspace can be found using a sparse SVD solver (we use
\noun{MATLAB}'s svds). Moreover, its inverse can be easily computed
using the formula
\[
\matU_{k}(\Lambda_{k}+\lambda\matI)^{-1}\matU_{k}^{\T}+\frac{1}{\lambda_{k}+\lambda}(\matI_{d}-\matU_{k}\matU_{k}^{\T}).
\]

We also experiment with variants of sketched iterative CCA (Algorithm
\ref{alg:CCAopt}) and sketched iterative FDA (Algorithm \ref{alg:LDAopt})
in which Riemannian CG is replaced with Riemannian trust-region Method.

\begin{figure}[t]
\begin{centering}
\begin{tabular}{ccc}
\includegraphics[width=0.3\columnwidth]{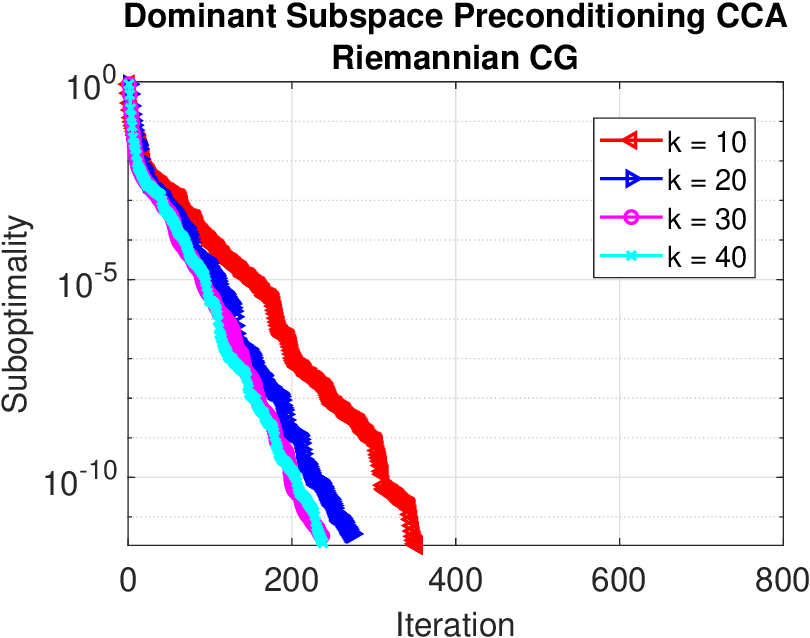} &
\includegraphics[width=0.3\columnwidth]{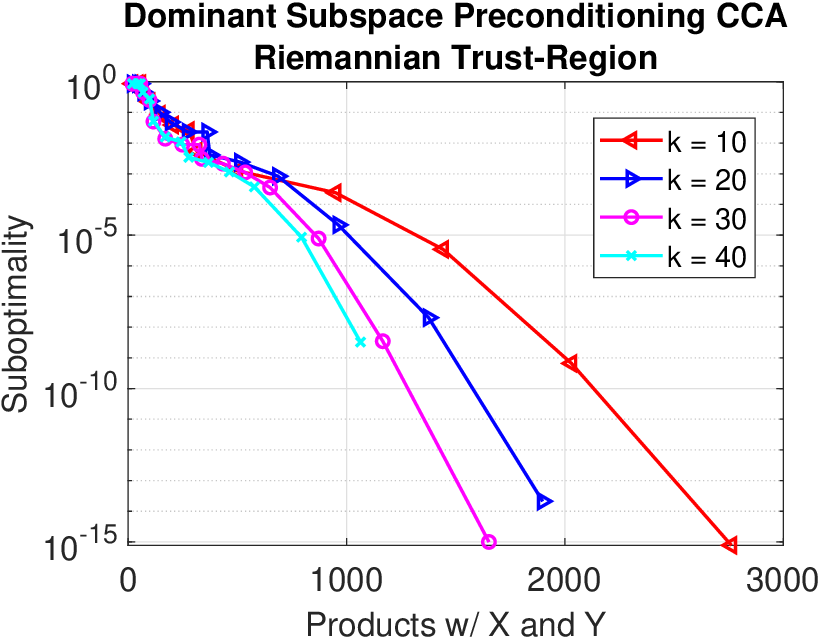} &
\includegraphics[width=0.3\columnwidth]{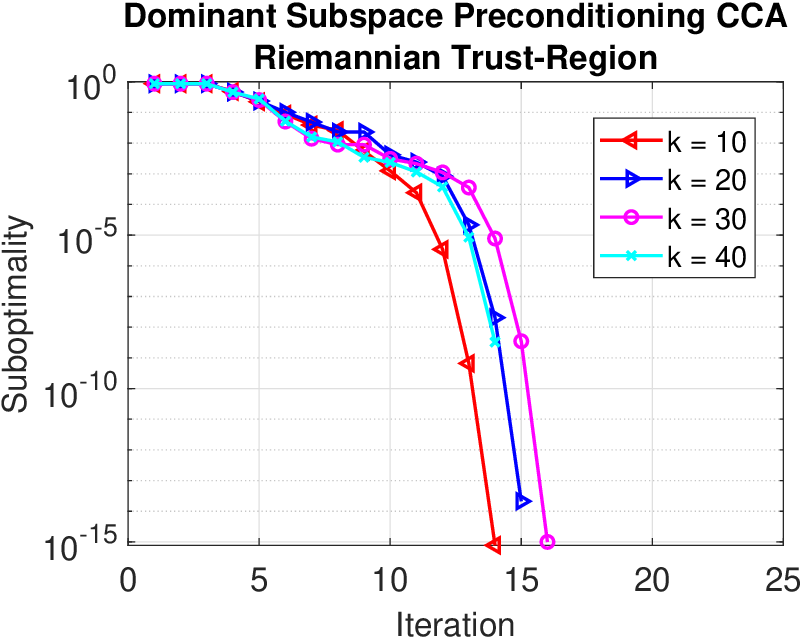}
\tabularnewline
\includegraphics[width=0.3\columnwidth]{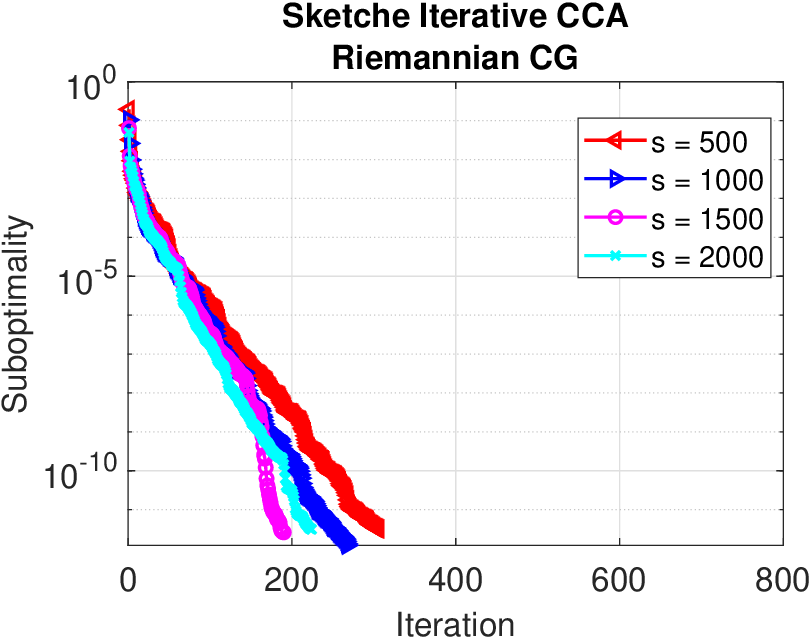} &
\includegraphics[width=0.3\columnwidth]{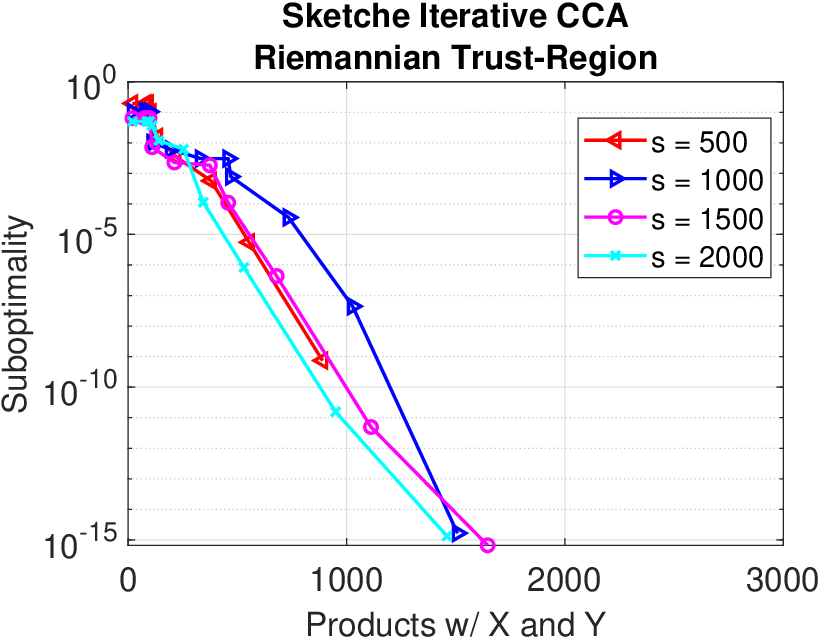} & \includegraphics[width=0.3\columnwidth]{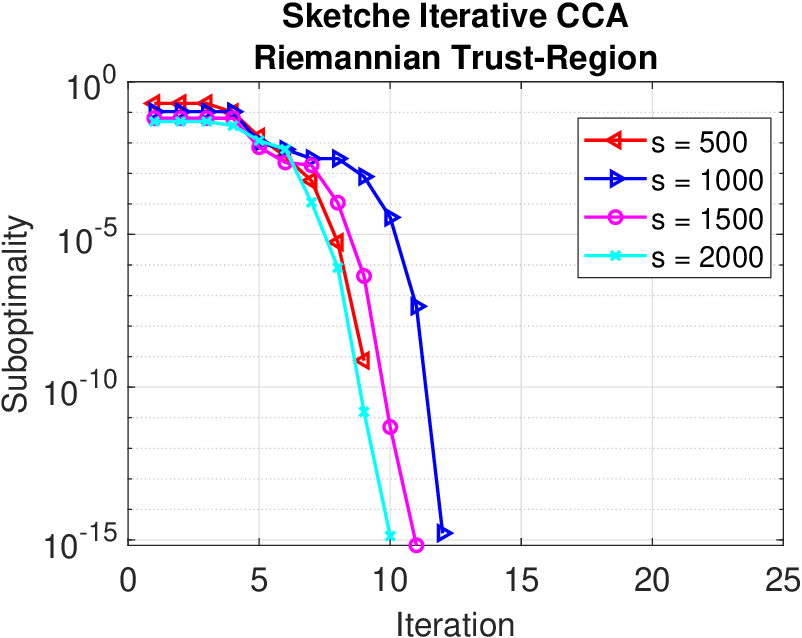}
\end{tabular}
\par\end{centering}
\caption{\label{fig:mnist-cca-p=00003D5}Results for CCA on \noun{MNIST}. For {\noun{CountSketch}} the number of rows is $s$, and $k$ is the number of singular vectors for the dominant
subspace preconditioner.}
\end{figure}

\begin{figure}[t]
\begin{centering}
\begin{tabular}{ccc}
\includegraphics[width=0.3\columnwidth]{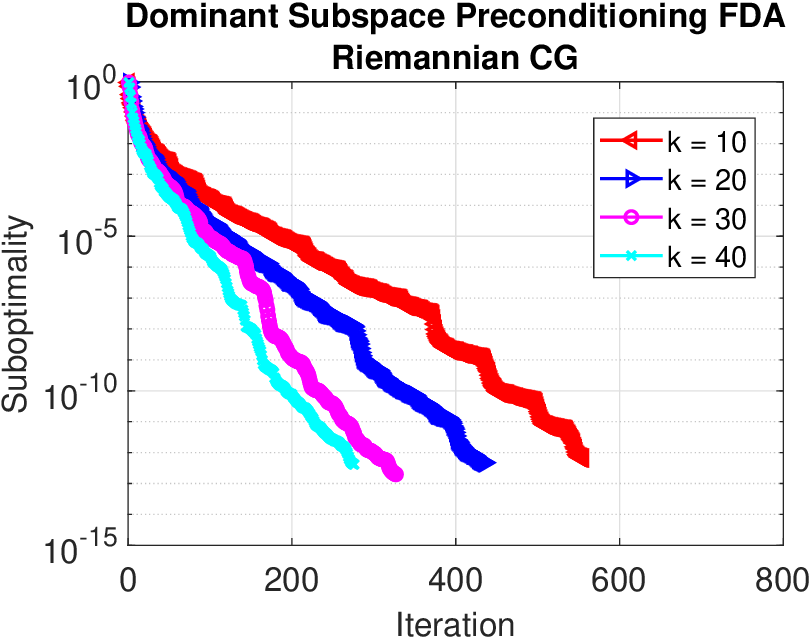} & \includegraphics[width=0.3\columnwidth]{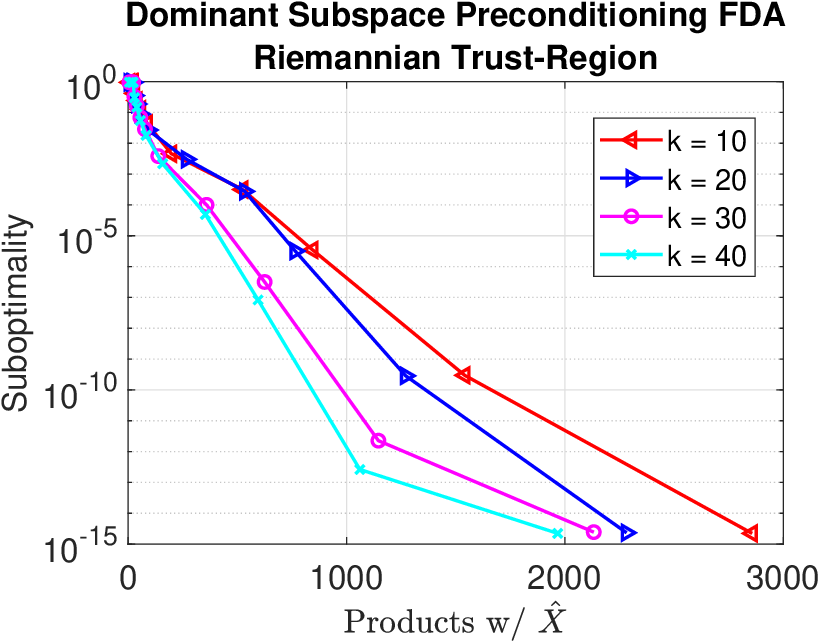} &
\includegraphics[width=0.3\columnwidth]{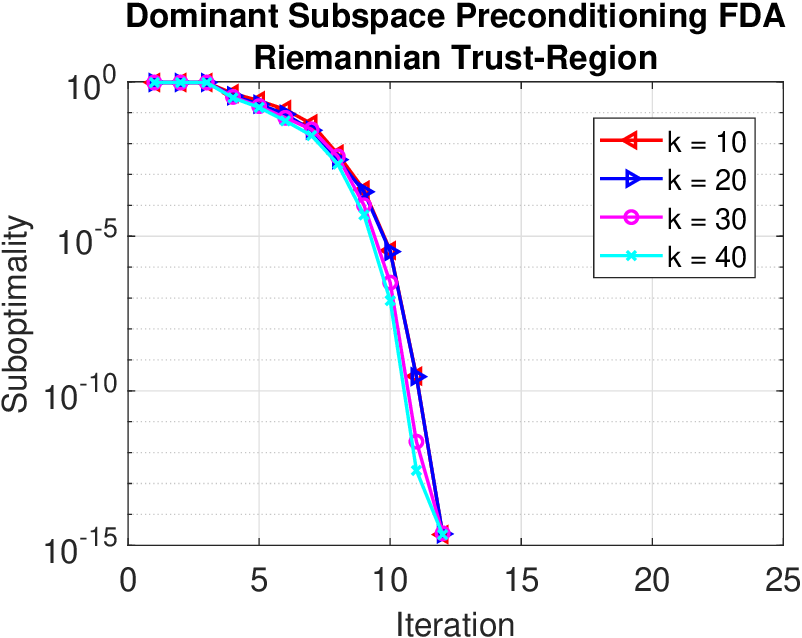}\tabularnewline
\includegraphics[width=0.3\columnwidth]{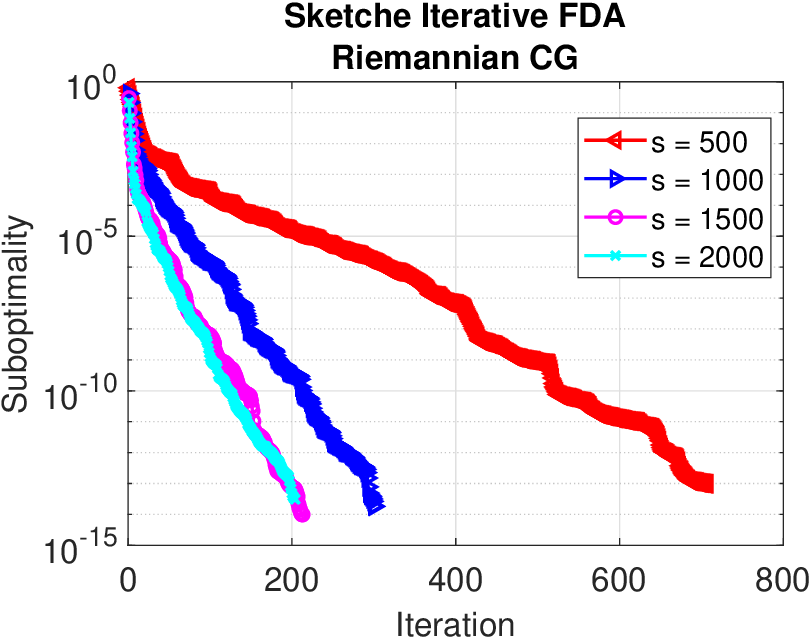} & 
\includegraphics[width=0.3\columnwidth]{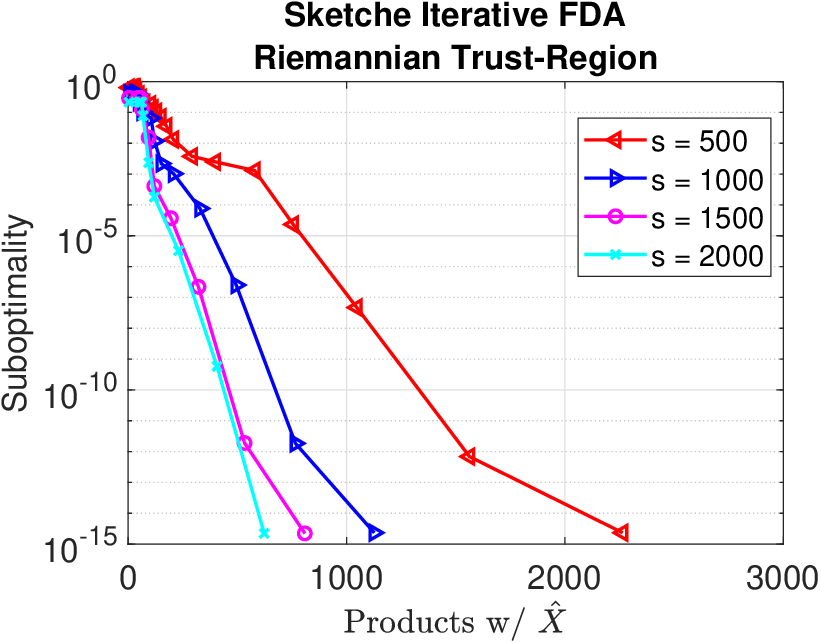} &
\includegraphics[width=0.3\columnwidth]{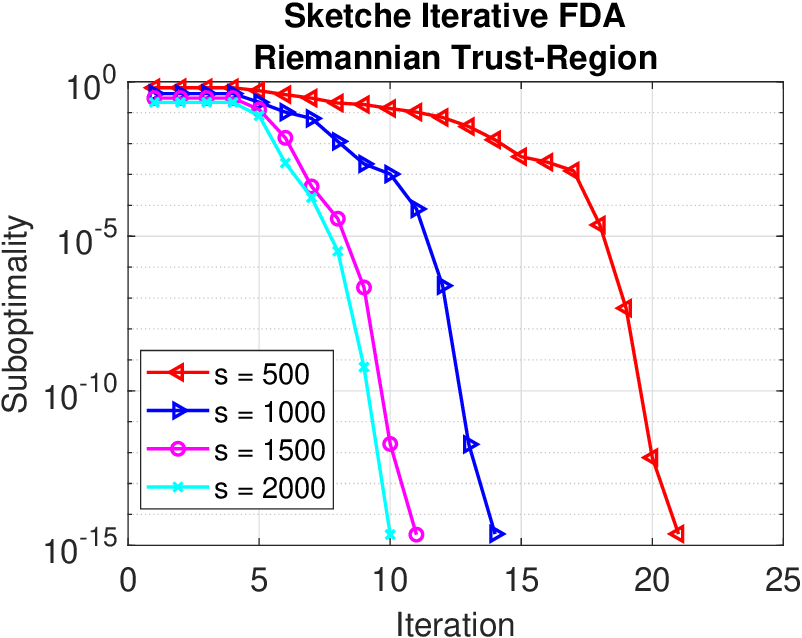} 
\end{tabular}
\par\end{centering}
\caption{\label{fig:mnist-lda-p=00003D5}Results for FDA on \noun{MNIST}. For {\noun{CountSketch}} the number of rows is $s$, and $k$ is the number of singular vectors for the dominant
subspace preconditioner.}
\end{figure}

\begin{figure}[t]
\begin{centering}
\begin{tabular}{ccc}
\includegraphics[width=0.3\columnwidth]{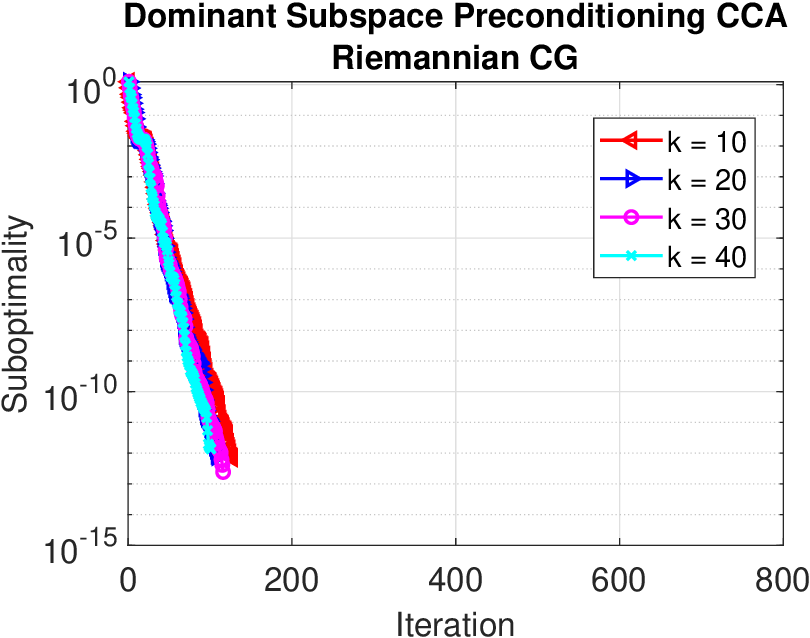} &
\includegraphics[width=0.3\columnwidth]{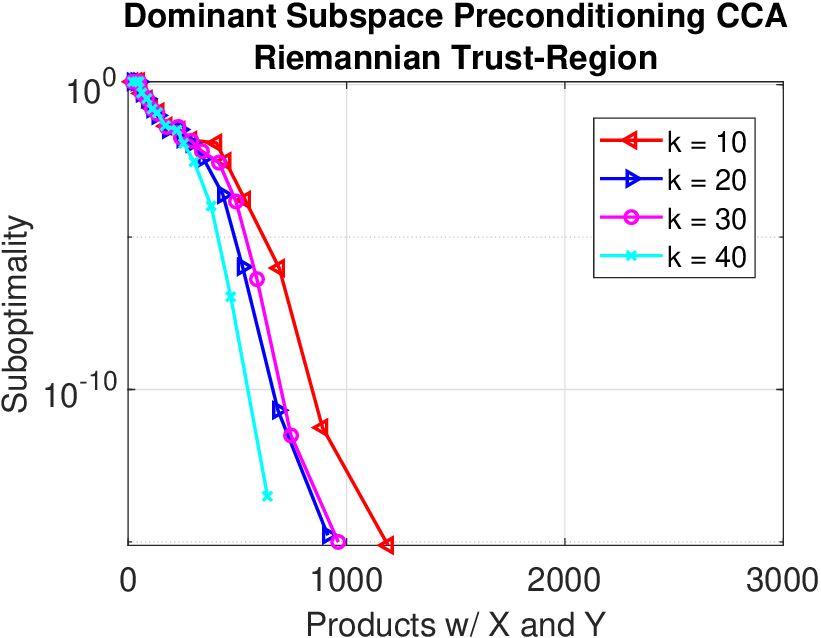} &
\includegraphics[width=0.3\columnwidth]{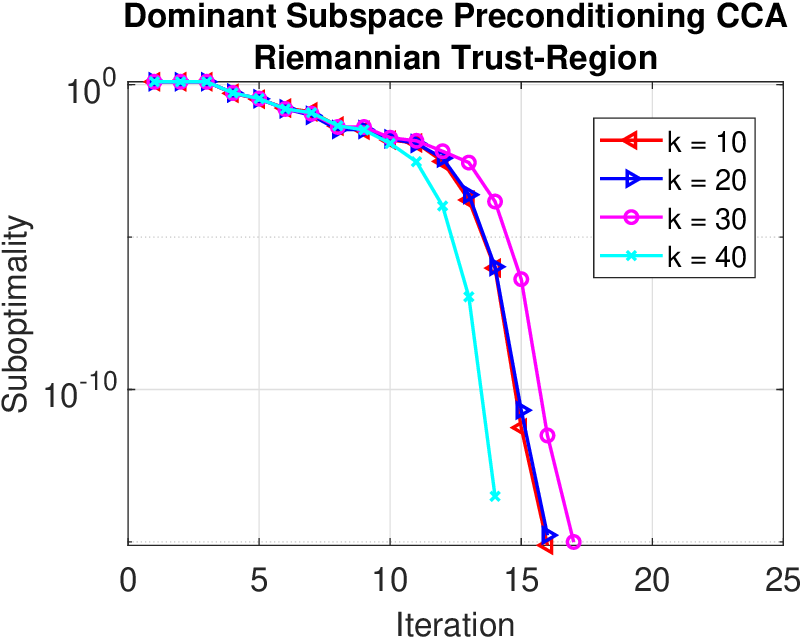}\tabularnewline
\includegraphics[width=0.3\columnwidth]{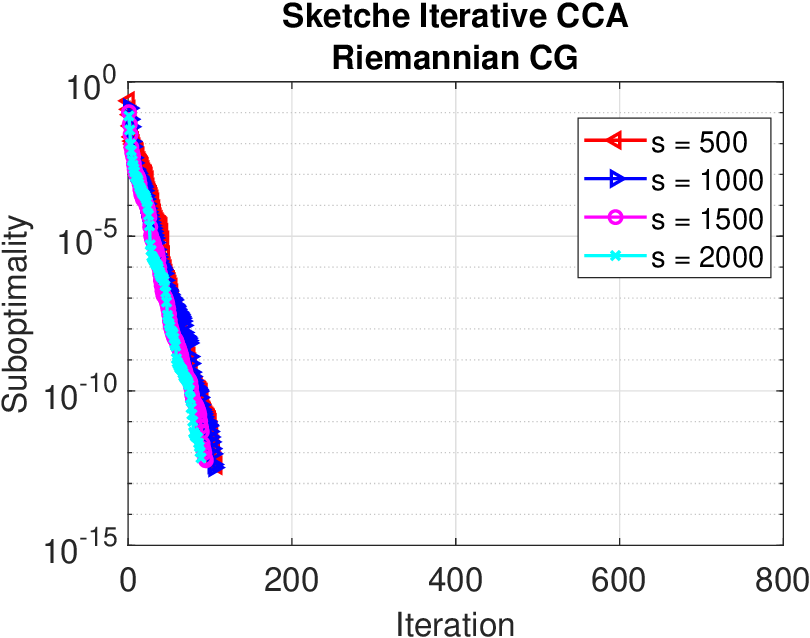} & 
\includegraphics[width=0.3\columnwidth]{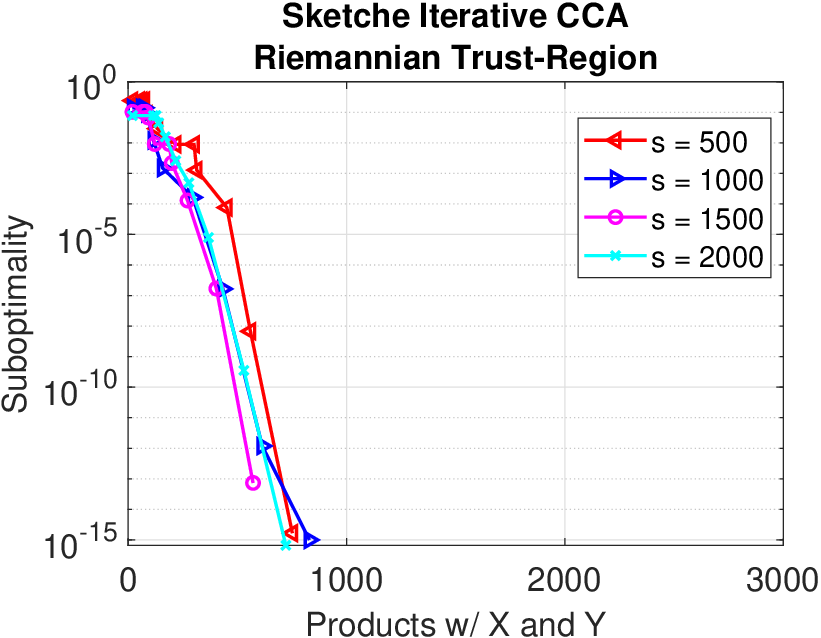} & 
\includegraphics[width=0.3\columnwidth]{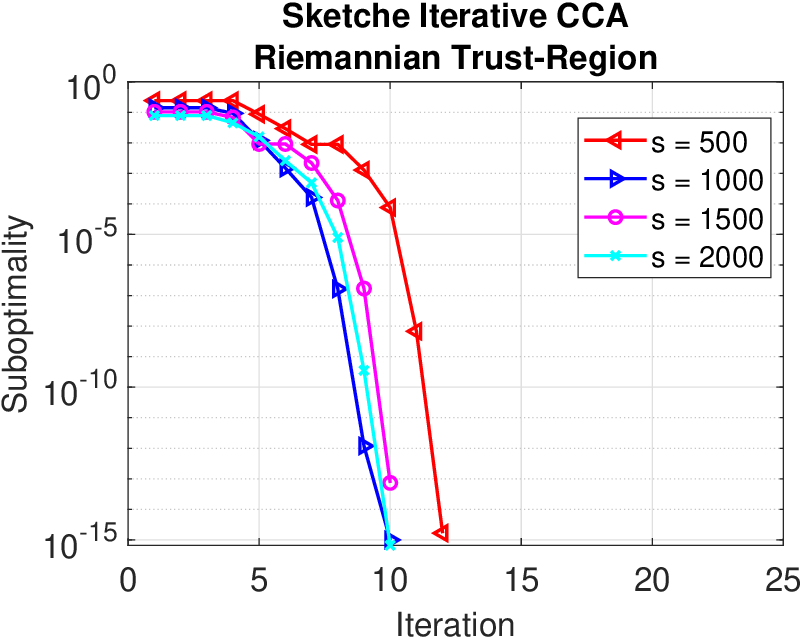}  
\end{tabular}
\par\end{centering}
\caption{\label{fig:cca-mediamill-p=00003D5}Results for CCA on \noun{MEDIAMILL}. For {\noun{CountSketch}} the number of rows is $s$, and $k$ is the number of singular vectors for the dominant
subspace preconditioner.}
\end{figure}

\begin{figure}[t]
\begin{centering}
\begin{tabular}{ccc}
\includegraphics[width=0.3\columnwidth]{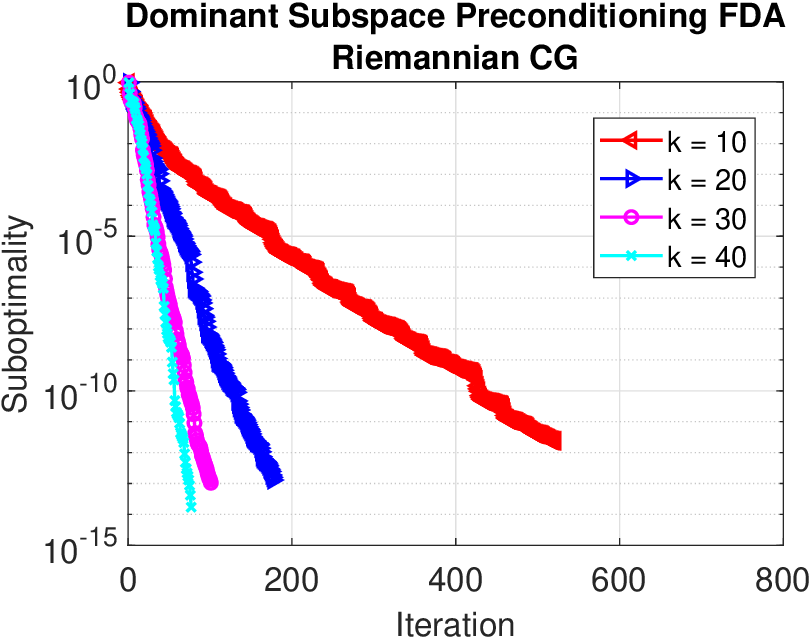} & 
\includegraphics[width=0.3\columnwidth]{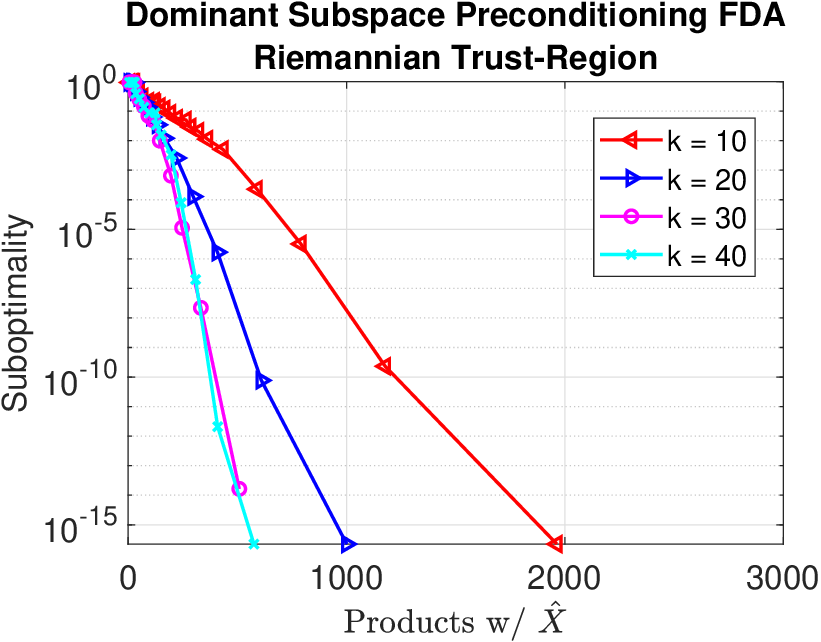} &
\includegraphics[width=0.3\columnwidth]{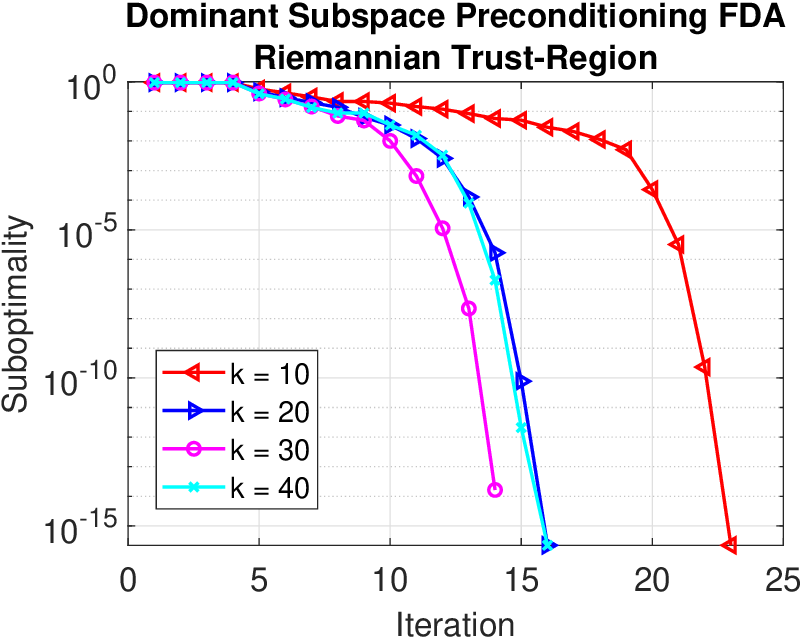}\tabularnewline
\includegraphics[width=0.3\columnwidth]{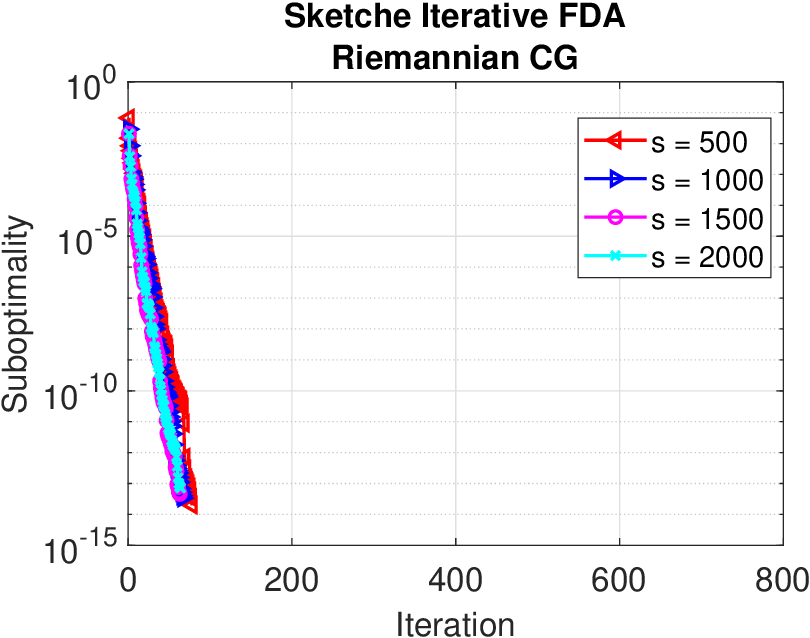} &
\includegraphics[width=0.3\columnwidth]{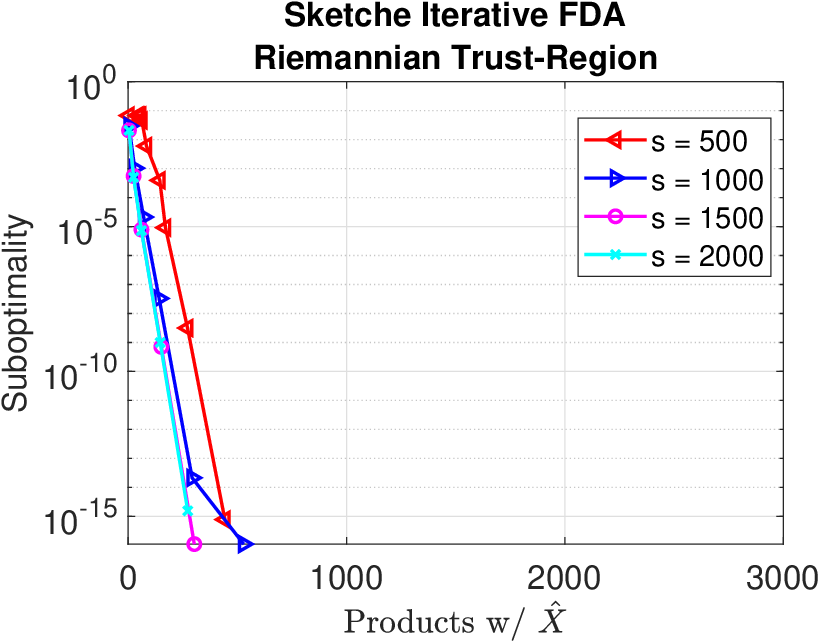} &
\includegraphics[width=0.3\columnwidth]{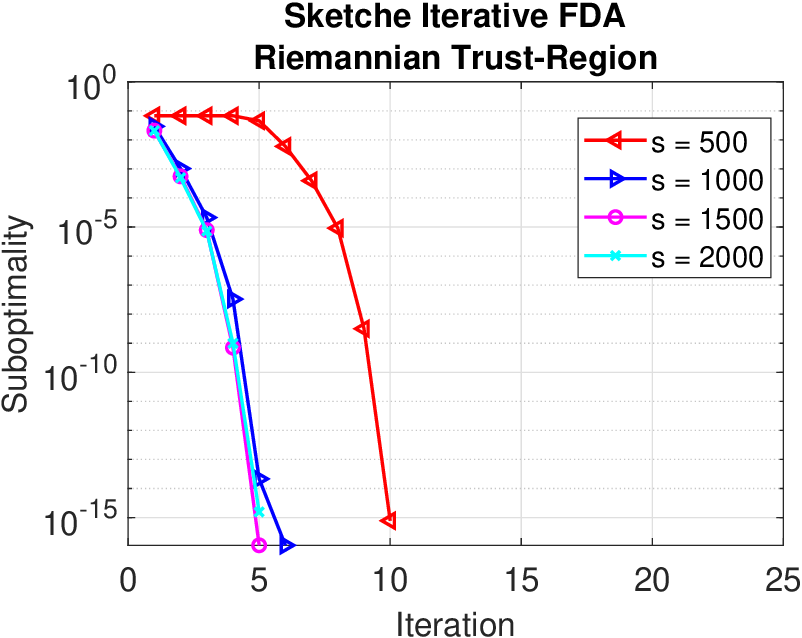}  
\end{tabular}
\par\end{centering}
\caption{\label{fig:lda-covtype-p=00003D5}Results for FDA on \noun{COVTYPE}. For {\noun{CountSketch}} the number of rows is $s$, and $k$ is the number of singular vectors for the dominant
subspace preconditioner.}
\end{figure}

We use \noun{MATLAB} for our implementations, relying on the \noun{manopt}
library \cite{manopt} for Riemannian optimization. In the \noun{manopt}
library we implemented the generalized Stiefel manifold with a non-standard
metric\footnote{\noun{manopt} has an implementation of the generalized Stiefel manifold,
but only with the standard metric $\matM=\matB$.}. The experiments we present here are with $p=3$ and $\matN=\diag{3,2.75,2}$. Recall that the only condition on the matrix $\matN$ is that it is diagonal, with strictly descending and positive diagonal elements. In the literature, typically, descending integers are chosen. We remark that this choice of $\matN$ is an arbitrary (not optimized) choice, and we achieved similar results for other choices of $\matN$ as a diagonal matrix with strictly decreasing values on the diagonal, and $p$ as well (which we exclude from this text). In the plots (Fig. \ref{fig:mnist-cca-p=00003D5}, \ref{fig:mnist-lda-p=00003D5}, \ref{fig:cca-mediamill-p=00003D5}, and \ref{fig:lda-covtype-p=00003D5}), the role of $s$ is the number of rows in the sketched data matrices, after applying a {\noun{CountSketch}} transformation on the data matrices. Similarly, $k$ is the number of singular values and vectors we use for the Dominant Subspace Preconditioning.

We did not optimize the implementation, so wall clock time is not
an appropriate metric for performance, instead, we use an alternative
metrics. In particular, the direct methods are about $10$ times faster than the iterative methods measured in wall clock time in the experiments in this subsection, however, as $n$ becomes larger, iterative methods become faster than the direct methods (see Subsection \ref{subsec:synexp}). In Riemannian trust-region Method, different iterations do a
variable amount of passes over the data, so we measure passes directly,
i.e., products with the data matrices, as this is the dominant cost
of our algorithm. As for Riemannian CG, the \noun{manopt} solver restricts
the number of line-search steps in each iteration to $25$, otherwise
the step is rejected. In practice, in our experiments the number of
products in different iterations is between $6$ to $26$ for CCA
and between $3$ to $13$ for FDA, so we report the number of iterations.
For every iteration or pass over the data, $t$, we plot the suboptimality of the current iterate: $|\sum_{i=1}^{p}\sigma_{i}\mu_{i}+\justfcca(\matZ_{t})|/\sum_{i=1}^{p}\sigma_{i}\mu_{i}$
for CCA where $\matZ_{t}$ represents the parameters at iteration $t$, and $|\sum_{i=1}^{p}\rho_{i}\mu_{i}+\justflda(\matW_{t})|/\sum_{i=1}^{p}\rho_{i}\mu_{i}$
for FDA where $\matW_{t}$ represents the parameters at iteration $t$. We use \noun{manopt}'s default stopping criteria: the optimization
process terminates if the norm of the Riemannian gradient drops below
$10^{-6}$. We cap the number of iterations by 1000.

We use in our experiments three popular data sets: \noun{MNIST} (Figures \ref{fig:mnist-cca-p=00003D5} and \ref{fig:mnist-lda-p=00003D5}), \noun{MEDIANILL} (Figure \ref{fig:cca-mediamill-p=00003D5})
and \noun{COVTYPE} (Figure \ref{fig:lda-covtype-p=00003D5})\footnote{data sets were downloaded for libsvm's website: https://www.csie.ntu.edu.tw/\textasciitilde cjlin/libsvmtools/datasets/.}.
\noun{MNIST} is used for testing CCA and FDA, where for CCA we try
to correlate the left side of the image to the right side of the image.
\noun{MEDIANILL} ($43,907$ examples) is a multilabel data set, so we use
it to test CCA. \noun{COVTYPE} is a large ($581,012$ examples) labeled
data set, and we use it to test FDA.

Consider Figure \ref{fig:mnist-cca-p=00003D5} (CCA on \noun{MNIST}). For Riemannian CG the
number of products per iteration is never bigger than $20$. As a
reference, for CCA the number of iterations required for CG and Trust-Region
with $\matM=\Sigma$ is $218$  and $15$  correspondingly, whereas
without a preconditioner ($\matM=\matI_{d}$) the CG does not converge
even after $1000$ iterations and the Trust-Region required $21$
iterations to converge. We clearly see the direct correspondence between
sketch quality (as measured by the sketch size $s$) and number of
iterations. Furthermore, the number of iterations is close to optimal
after sketching to only $s=2000$ rows (there are 60,000 examples in the
original data set) or using only 40 singular vectors for the dominant
subspace preconditioner (there are 784 features in the data set)\footnote{Interestingly, with $s\geq500$ the subspsace embedding preconditioner
uses less iterations than the optimal preconditioner. This is because
of the use of sketching based warm-start.}.

Consider Figure \ref{fig:mnist-lda-p=00003D5} (FDA on \noun{MNIST}). For Riemannian
CG the number of products per iteration is never bigger than $9$.
As a reference, for FDA the number of iterations required for CG and
Trust-Region with $\matM=\matS_{\w}+\lambda\matI_{d}$ is $95$  and $11$  correspondingly,
whereas without a preconditioner ($\matM=\matI_{d}$) the CG does not
converge even after $1000$ iterations and the Trust-Region required
$20$ iterations to converge, so again we see that the sketching reduces
the number of iterations.

Consider Figure \ref{fig:cca-mediamill-p=00003D5} (CCA on \noun{MEDIANILL}). For Riemannian CG the number of products per
iteration is never bigger than $20$. As a reference, the number of
iterations required for CG and Trust-Region with $\matM=\Sigma$
is $107$  and $15$  correspondingly, whereas without a preconditioner
($\matM=\matI_{d}$) the number of iterations for CG is $550$ and
for Trust-Region is $23$. The data set has 30,993 examples and 221
features, so again we see that we can sketch to relatively small size
($s=2000$ or $k=40$) and get an effective preconditioner.

Consider Figure \ref{fig:lda-covtype-p=00003D5} (FDA on \noun{COVTYPE}). For Riemannian CG the number of products per
iteration is never bigger than $10$. As a reference, the number of
iterations required for CG and Trust-Region with $\matM=\matS_{\w}+\lambda\matI_{d}$
is $71$  and $19$  correspondingly, whereas without a preconditioner
($\matM=\matI_{d}$) the CG does not converge even after $1000$ iterations
and the Trust-Region required $44$ iterations to converge. Considering
that the data set has over half a million examples, subspace embedding
preconditioning is highly effective, as it sketches the data to a
comparatively very small size. The data set has only 50 features,
so dominant subspace preconditioning is less effective for this data
set.

\subsection{\label{subsec:synexp}Synthetic Experiments and Comparison with ALS}

\begin{figure}[t]
\begin{centering}
\begin{tabular}{cc}
\includegraphics[width=0.5\columnwidth]{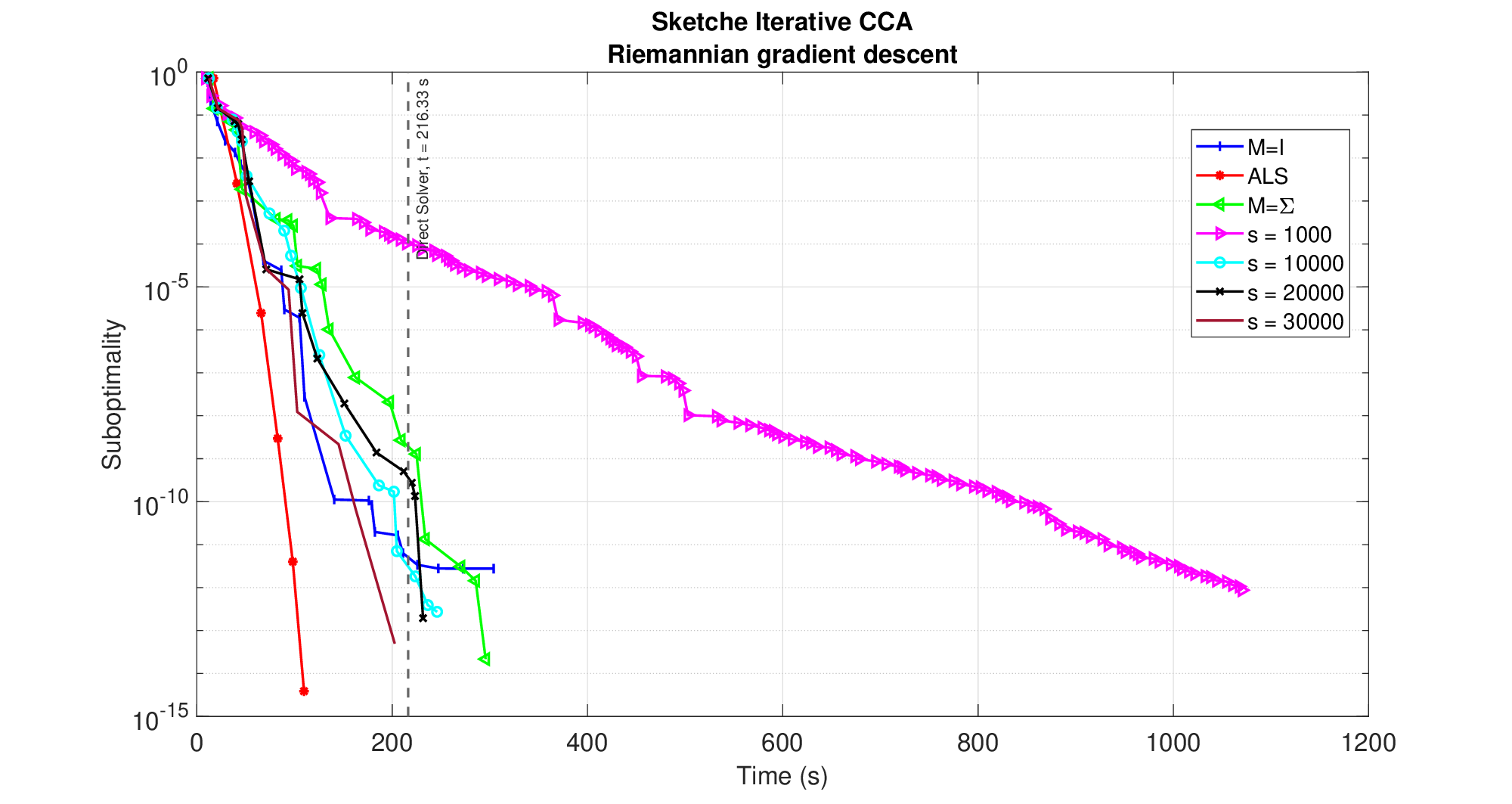} &
\includegraphics[width=0.5\columnwidth]{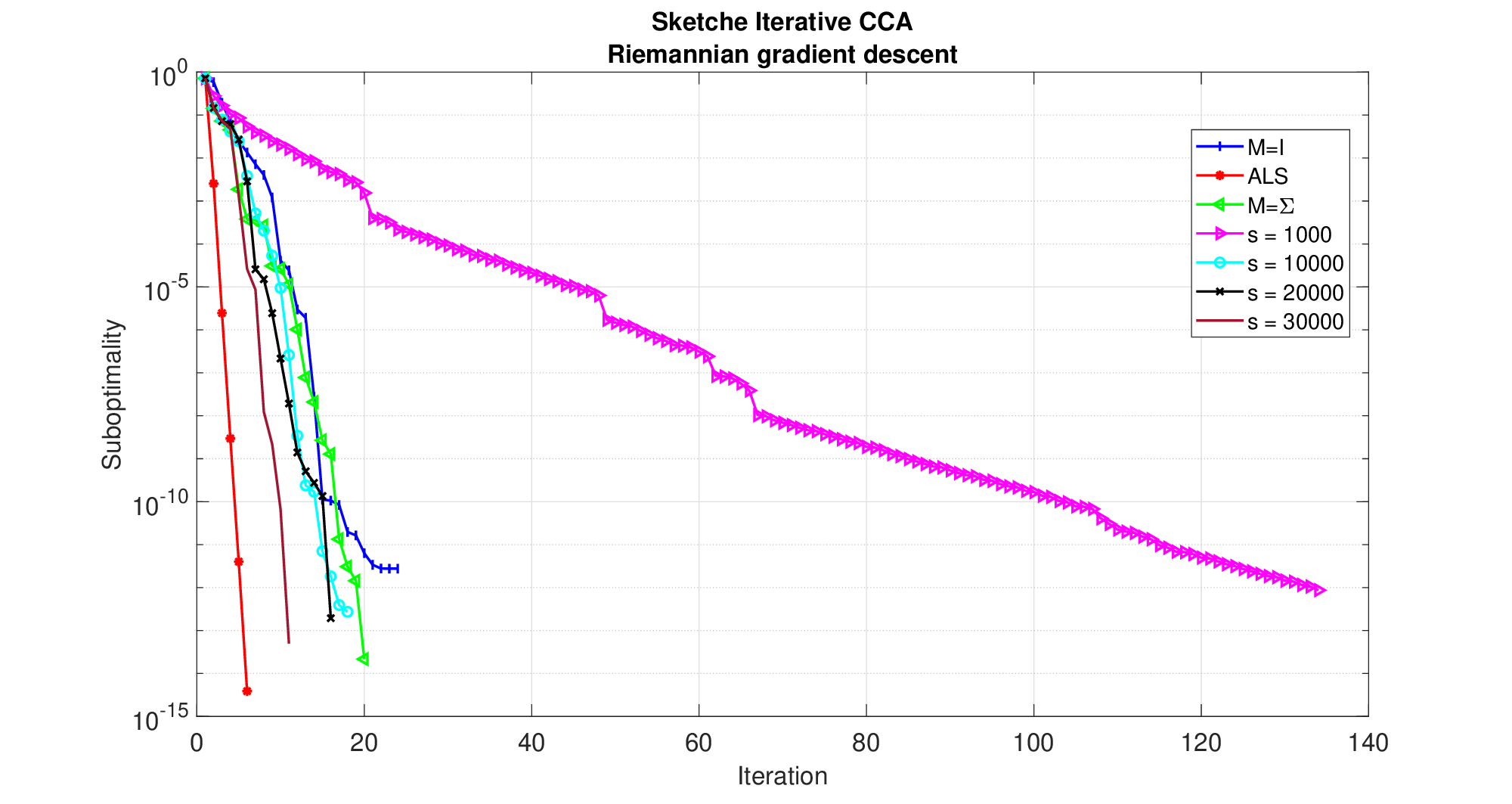} 
\end{tabular}
\par\end{centering}
\caption{\label{fig:CCA_reviews}Results for CCA on a synthetic experiment. Left - suboptimality vs. time in seconds. Right - suboptimality vs. iteration count. we use \noun{CountSketch} as the sketching transform, where $s$ is the number of rows after sketching is applied to the data matrices.}
\end{figure}

\begin{figure}[t]
\begin{centering}
\begin{tabular}{cc}
\includegraphics[width=0.5\columnwidth]{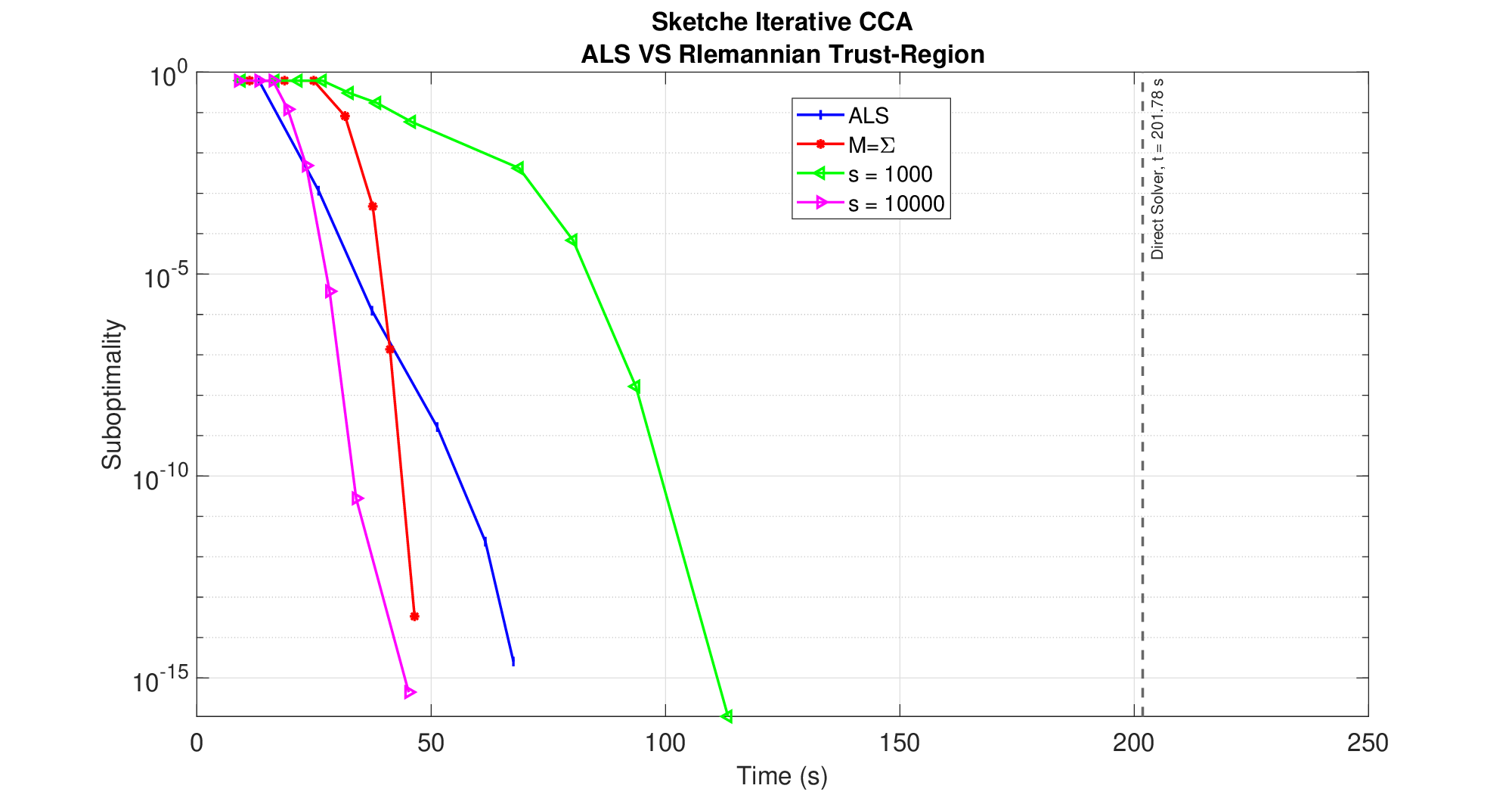} &
\includegraphics[width=0.5\columnwidth]{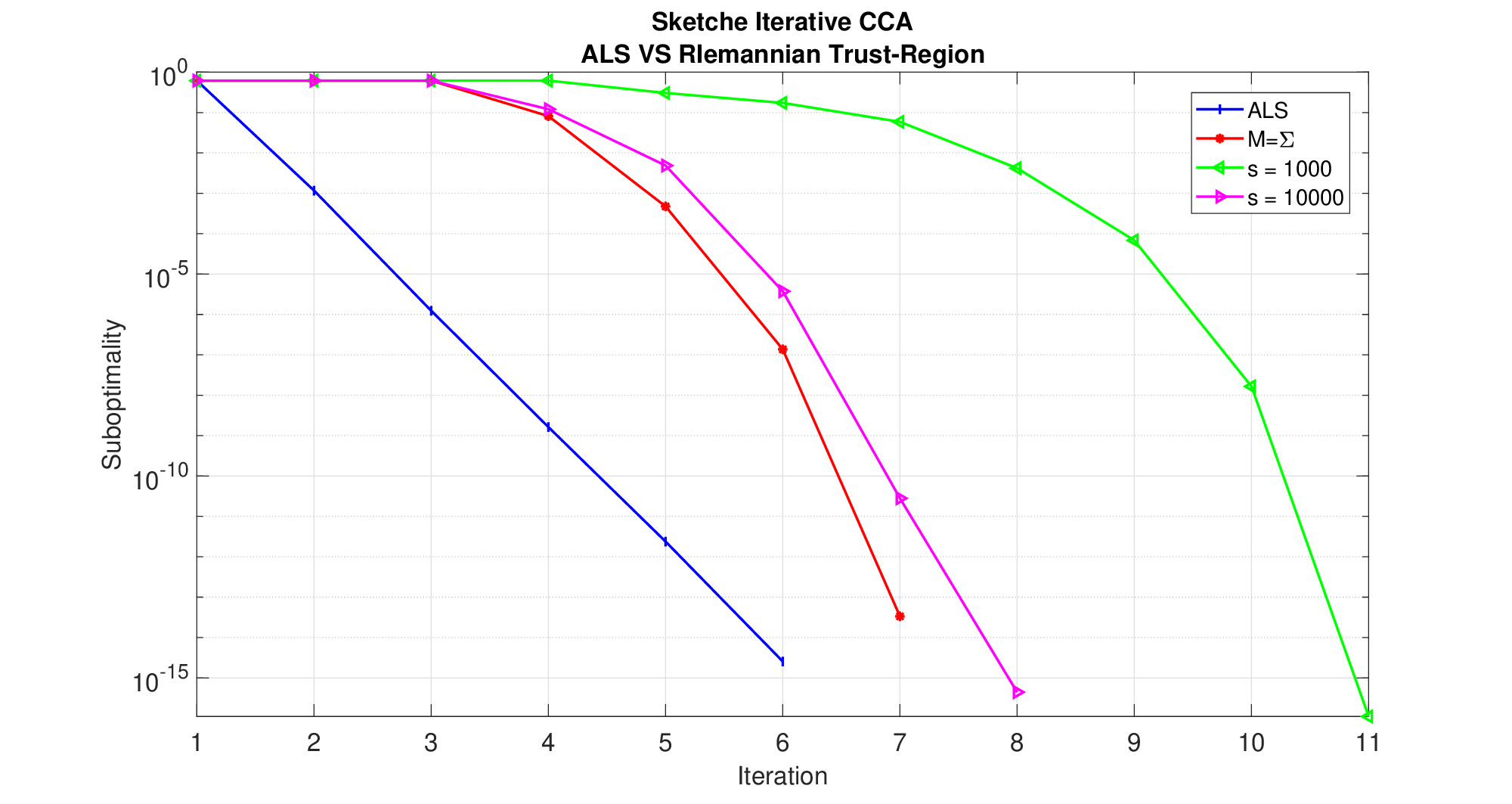} 
\end{tabular}
\par\end{centering}
\caption{\label{fig:CCA_reviews_trust}Results for CCA with Riemannian trust-region on a synthetic experiment. Left - suboptimality vs. time in seconds. Right - suboptimality vs. iteration count. we use \noun{CountSketch} as the sketching transform, where $s$ is the number of rows after sketching is applied to the data matrices.}
\end{figure}

In addition, we report a set of synthetic experiments for CCA with comparisons to other methods. To perform comparisons to existing methods, we conducted a set of synthetic experiments in which we compared the running times of our algorithms for CCA to ALS algorithm for CCA \cite[Algorithm 1]{GeEtAl16}, that is \cite[Algorithm 5.2]{golub1995canonical}. Note that the aforementioned algorithm is equivalent to a Riemannian gradient descent with a standard metric and a step-size of $\alpha_t=(\u_t^{\T}\Sigma_{\x\y}\v_t)^{-1}$ at each iteration $t$ (this is valid whenever $\alpha_0=(\u_0^{\T}\Sigma_{\x\y}\v_0)^{-1} > 0$), for $p=1$. Thus, we add in Fig. \ref{fig:CCA_reviews} experiments for $p=1$ of Riemannian gradient descent with a step-size of $\alpha_t=(\u_t^{\T}\Sigma_{\x\y}\v_t)^{-1}$. We compare using the standard metric with a prescribed step-size \cite[Algorithm 5.2]{golub1995canonical} vs. line search and our sketched metrics. In addition, we add the running times for direct methods. For FDA, one can perform similar experiments, comparing our algorithms to the power method for $(\ldaB)^{-\nicehalf} \mat{S_{B}} (\ldaB)^{-\nicehalf}$, as it is equivalent to Riemannian gradient descent with the standard metric and a step-size of $\alpha_t = (\w_t^{\T}\mat{S_{B}}\w_t)^{-1}$ for $p=1$ (this is valid whenever $\alpha_0 = (\w_0^{\T}\mat{S_{B}}\w_0)^{-1} > 0$). We do not include these experiments in the paper as the results are similar to the experiments for CCA. Moreover, we also added experiments with Riemannian trust-region in Fig. \ref{fig:CCA_reviews_trust}, to illustrate an advantage Riemannian optimization framework provides, e.g., second-order methods. 

Our method is most beneficial in the regime $n\gg d$ and ill-conditioned data matrices. Thus, the experiments presented in Fig. \ref{fig:CCA_reviews} and Fig. \ref{fig:CCA_reviews_trust} are performed as follows. A synthetic data set is generated with two ill-conditioned random data matrices $
\matX\in\R^{10^6 \times 500}$ and $\matY\in \R^{10^6 \times 500}$ with a condition number of the order of $10^{14}$ correspondingly. We use a small regularization of $10^{-6}$ multiplied by the average eigenvalue of the Gram matrices of $\matX$ and $\matY$ correspondingly. We measure suboptimality vs. time in seconds (including pre-processing time of creating the sketched matrices) and iteration count. 

In Fig. \ref{fig:CCA_reviews}, we present results for ALS algorithm and gradient descent with line search with $\matM = \matI$, i.e., no preconditioning, $\matM = \Sigma$, i.e., optimal preconditioning, and four sketched matrices with \noun{CountSketch} of sizes $s= 10^3, 10^4, 2*10^4, 3*10^4$. 

ALS is the fastest ($109.7$ seconds). Solving with $\matM = \Sigma$ takes $295.72$ seconds including pre-processing, for $s=10^4$ ($0.1$\% of the data) it takes $245.99$ seconds including pre-processing, for $s=2*10^4$ ($0.2$\% of the data) it takes $231.62$ seconds including pre-processing, and for $s=3*10^4$ ($0.3$\% of the data) it takes $202.64$ seconds including pre-processing. In particular, excluding pre-processing time, solving with $\matM = \Sigma$ takes $289.5$ seconds, whereas, for $s=10^4$ it takes $241.77$ seconds, for $s=2*10^4$ it takes $227.74$ seconds, for $s=3*10^4$ it takes $197.74$ seconds. Note that the pre-processing time is almost identical and faster for the sketched matrices since storing and computing the full data matrices takes considerable time. For comparison, solving directly using \cite[Section 5]{bjorck1973numerical} takes $216.33$ seconds. Note that for $s=3*10^4$ the iterative method is faster, due to the larger $n$. Note that the condition number of the Riemannian Hessian at the optimum changes according to the preconditioning, as the ill-conditioning of the problem arises from the data matrices. For $\matM=\matI$ the condition number computed via \noun{Manopt} is $2.66*10^5$, for $s=10^3$ it is $34.1$, however for $s=10^4$ it is reduced to $2.48$, for $s=2*10^4$ it becomes $1.89$, for $s=3*10^4$ it becomes $1.68$, and for $\matM=\Sigma$ it becomes $1.09$. Indeed, we can observe that our randomized preconditioning scheme behaves as theory predicts, having better conditioning of the Riemannian Hessian at the optimum, as our sketched matrices approximate better and better the standard metric. Moreover, as seen in these results, for very large matrices, sketching indeed accelerates the algorithms, and results in faster methods. 

In Fig. \ref{fig:CCA_reviews_trust}, we present Riemannian trust-region compared to ALS, with the metric $\matM=\Sigma$, and the metrics achieved by sketching with $s=10^3$ and $s=10^4$ (results for $\matM=\matI$ were very similar to $s=10^3$). All the iterative methods were faster than solving directly using \cite[Section 5]{bjorck1973numerical} ($201.78$ seconds). While $s=10^3$ achieved the slowest convergence  ($113.32$ seconds including preprocessing) among the iterative methods, due to its ill-conditioning, both $\matM=\Sigma$ and $s=10^4$ were faster than the ALS, and $s=10^4$ was the fastest ($45.18$ seconds including preprocessing, and similar runtime without preprocessing of the algorithm to $\matM=\Sigma$). 

\section{Conclusions}

In this paper we propose faster randomized methods for orthogonality
constrained problems. Our method is specifically designed for typical
structure of orthogonality constraints in machine learning, i.e., the
constraints are defined by a Gram matrix of a data matrix where one
dimension is much larger than the other. We use the framework of Riemannian
optimization as the underlying iterative methods which we precondition
by incorporating a randomized preconditioner via the Riemannian metric
(Riemannian preconditioning). The aforementioned technique can be
used to precondition any core Riemannian method. Our method can also
be applied to constraints which are described by the product of two
or more generalized Stiefel manifolds.

We demonstrated our method and proposed a preconditioning strategy
for two problems: CCA and FDA. For both of these examples we evaluate
the computational costs, and bound the condition number of the Riemanninan
Hessian at the optimum. This in turn, allows us to reason about the
effect of the proposed randomized preconditioner.

As a future research direction we believe our method can be extended
beyond orthogonality constrained problems for other equality constraints
by identifying similar problem structure.

\section*{Acknowledgements}
The authors thank Bart Vandereycken for useful discussions. This
research was supported by the Israel Science Foundation (grant no.
1272/17).

\appendix

\section{\label{sec:analogue_of_lem_SRHT} Analogue of Lemma \ref{lem:sketching} for SRHT}
In this section, we formulate an analogue for Lemma \ref{lem:sketching}, which was formulated for \noun{CountSketch} sketching transform to SRHT. First, we recall the definition of SRHT from \cite[Section 1.1]{tropp2011improved} and \cite[Definition 1.2]{boutsidis2013improved}.
\begin{definition}\label{def:SRHT}
Fix integers $s$ and $n=2^m$ where $s<n$, and $m=1,2,3,...$. An SRHT matrix is an
$s\times n$ matrix of the form
\begin{equation*}
\matS = \sqrt{\frac{s}{n}} \matR \mat{H} \matD,
\end{equation*}
where $\matD\in \R^{n\times n}$ is a random diagonal matrix whose entries are independent random signs, i.e., random variables uniformly distributed on ${\pm 1}$, $\mat{H} \in \R^{n \times n}$ is a normalized Walsh-Hadamard matrix, and $\matR\in \R^{s\times n}$ is a subset of $s$ rows from the $n\times n$ identity matrix, where the rows are chosen uniformly at random and without replacement. 
\end{definition}
Next, recall that according to \cite[Theorem 2.1]{ailon2009fast} and \cite[Lemma 1.3]{boutsidis2013improved}, constructing $\matS$ and computing $\matS x$ given $\x\in \R^n$, $s<n$, can be done in $O(n\log{(s)})$ operations. Thus, computing $\matS \matZ$ for $\matZ\in \R^{n\times d}$ takes $O(nd\log{(s)})$ operations (see Table \ref{tab:costs}).

Finally, we prove an analogue of Lemma \ref{lem:sketching} for SRHT using the approximate matrix multiplication property for SRHT \cite[Theorem 9]{https://doi.org/10.48550/arxiv.1507.02268}.

\begin{lemma}\label{lem:sketching_SRHT}
Assume that $\lambda>0$ or that $\matZ\in\R^{n\times d}$
has full column rank.
Suppose that $\matS\in \R^{s \times n}$ is an SRHT matrix with $$s = \Omega(4 s_{\lambda}(\matZ)^2(s_{\lambda}(\matZ) / \TNorm{\matQ}^2 +\log{(2 s_{\lambda}(\matZ)/\delta)}\log{(s_{\lambda}(\matZ) / (\TNorm{\matQ}^2 \delta)}))),$$ where $\matQ$ is the matrix $\matQ$ in the $\lambda$-QR factorization of $\matZ$ \cite[Definition 28]{ACW16b},
for some $\delta\in(0,1)$. Then with probability of at least $1-\delta$
we have that all the generalized eigenvalues of the pencil $(\matZ^{\T}\matZ+\lambda\matI,\matZ^{\T}\matS^{\T}\matS\matZ+\lambda\matI)$
are contained in the interval $[1/2,3/2]$ and $\kappa(\matZ^{\T}\matZ+\lambda\matI,\matZ^{\T}\matS^{\T}\matS\matZ+\lambda\matI)\leq3$.
\end{lemma}

\begin{proof}
The proof is similar to the proof of Lemma \ref{lem:sketching}. To show that the generalized eigenvalues of the pencil
$(\matZ^{\T}\matZ+\lambda\matI,\matZ^{\T}\matS^{\T}\matS\matZ+\lambda\matI)$
are contained in the interval $[1/2,3/2]$ and that $\kappa(\matZ^{\T}\matZ+\lambda\matI,\matZ^{\T}\matS^{\T}\matS\matZ+\lambda\matI)\leq3$
to hold, it is enough to show that 
\[
\frac{1}{2}(\matZ^{\T}\matZ+\lambda\matI)\preceq\matZ^{\T}\matS^{\T}\matS\matZ+\lambda\matI\preceq\frac{3}{2}(\matZ^{\T}\matZ+\lambda\matI)\,.
\]
Let $\matZ=\matQ\matR$ be a $\lambda$-QR factorization of $\matZ$ \cite[Definition 28]{ACW16b}. Left-multiplying by $\matR^{-\T}$ and right-multiplying by $\matR^{-1}$
on both sizes, we find it suffices to show that with probability of
at least $1-\delta$ we have
\[
\frac{1}{2}\matI_{d}\preceq\matQ^{\T}\matS^{\T}\matS\matQ+\lambda\matR^{-\T}\matR^{-1}\preceq\frac{3}{2}\matI_{d}
\]
 or, equivalently, 
\[
\TNorm{\matQ^{\T}\matS^{\T}\matS\matQ+\lambda\matR^{-\T}\matR^{-1}-\matI_{d}}\leq\frac{1}{2}\,.
\]
Since $\matQ^{\T}\matS^{\T}\matS\matQ+\lambda\matR^{-\T}\matR^{-1}-\matI_{d}=\matQ^{\T}\matS^{\T}\matS\matQ-\matQ^{\T}\matQ$
\cite[Fact 29]{ACW16b}, it is enough to show that 
\[
\TNorm{\matQ^{\T}\matS^{\T}\matS\matQ-\matQ^{\T}\matQ}\leq\frac{1}{2}\,.
\]
From \cite[Theorem 9]{https://doi.org/10.48550/arxiv.1507.02268}, we have that for any fixed matrix $\matA$, and an SRHT matrix $\matS_{0}$ with $m = \Omega ((\FNorm{\matA}^2 / \TNorm{\matA}^2 +\log{(1/(\varepsilon \delta))}\log{(\FNorm{\matA}^2 / (\TNorm{\matA}^2 \delta)})) / \varepsilon^2)$
rows, taking $k \coloneqq \FNorm{\matA}^2 / \TNorm{\matA}^2$ in \cite[Theorem 9]{https://doi.org/10.48550/arxiv.1507.02268}, we have that with probability of at least $1-\delta$, 
\[
\TNorm{\matA^{\T}\matS_{0}^{\T}\matS_{0}\matA-\matA^{\T}\matA}\leq\epsilon\cdot\TNorm{\matA}^2\leq \epsilon\cdot\FNorm{\matA}^2\,.
\]
Since $\FNormS{\matQ}=s_{\lambda}(\matZ)$ \cite[Fact 30]{ACW16b}, then with
$$s = \Omega (4 s_{\lambda}(\matZ)^2(s_{\lambda}(\matZ) / \TNorm{\matQ}^2 +\log{(2 s_{\lambda}(\matZ)/\delta)}\log{(s_{\lambda}(\matZ) / (\TNorm{\matQ}^2 \delta)}))),$$ we have 
\[
\TNorm{\matQ^{\T}\matS^{\T}\matS\matQ-\matQ^{\T}\matQ}\leq\frac{1}{2}
\]
with probability of at least $1-\delta$.
\end{proof}

\section{Omitted Proofs From Section~\ref{sec:Randomized-preconditioning-for-CCA}}
In this section we give omitted proofs from Section~\ref{sec:Randomized-preconditioning-for-CCA}.

\subsection{Proof of Theorem \ref{thm:criticalpointofCCA}}\label{subsec:proofcriticalpointofCCA}

\begin{proof}
Recall that the critical points are defined to be the points where the Riemannian gradient
is zero, and as such, whether a point is a critical point or not does
not depend on the choice of Riemannian metric (see, \cite[Eq. 3.31]{AMS09}).
Thus, for the sake of identifying the critical points, we can use
$\matM_{\matZ}\coloneqq \Sigma$ for all $\matZ\in\elpCCA$ and get a simplified form for the Riemannian gradient (Eq. \eqref{eq:gradCCA}):
\begin{eqnarray*}
\grad{\fcca} & = & -\left[\begin{array}{c}
\Sigma_{\x\x}^{-1}\Sigma_{\x\y}\matV\matN-\matU\sym{\left(\matU^{\T}\Sigma_{\x\x}\Sigma_{\x\x}^{-1}\Sigma_{\x\y}\matV\matN\right)}\\
\Sigma_{\y\y}^{-1}\Sigma_{\x\y}^{\T}\matU\matN-\matV\sym{\left(\matV^{\T}\Sigma_{\y\y}\Sigma_{\y\y}^{-1}\Sigma_{\x\y}^{\T}\matU\matN\right)}
\end{array}\right]\\
 & = & -\left[\begin{array}{c}
\left(\matI_{d_{\x}}-\matU\matU^{\T}\Sigma_{\x\x}\right)\Sigma_{\x\x}^{-1}\Sigma_{\x\y}\matV\matN+\matU\skew{(\matU^{\T}\Sigma_{\x\y}\matV\matN)}\\
\left(\matI_{d_{\y}}-\matV\matV^{\T}\Sigma_{\y\y}\right)\Sigma_{\y\y}^{-1}\Sigma_{\x\y}^{\T}\matU\matN+\matV\skew{(\matV^{\T}\Sigma_{\x\y}^{\T}\matU\matN)}
\end{array}\right].
\end{eqnarray*}

There is a connection between the canonical correlations and the singular value decomposition of $\matR\coloneqq\Sigma_{\x\x}^{-\nicehalf}\Sigma_{\x\y}\Sigma_{\y\y}^{-\nicehalf}$ \cite{WangEtAl16}: a pair of vectors $\u$ and $\v$ are canonical correlation vectors corresponding to the same canonical correlation if and only if $\tilde{\u}=\Sigma_{\x\x}^{\nicehalf}\u$ and  $\tilde{\v}=\Sigma_{\y\y}^{\nicehalf}\v$ are left and right singular vectors of $\matR$ corresponding to the same singular value.

canonical correlation vectors are the columns of $\matU$ and $\matV$ it and only if the columns of $(\tilde{\matU},\tilde{\matV})=(\Sigma_{\x\x}^{\nicehalf}\matU,\Sigma_{\y\y}^{\nicehalf}\matV)$ are $p$ left and right coordinated singular vectors of the matrix $\matR$. Thus, we use this relation to prove this theorem.

To show that all $\matZ\in\elpCCA$ such that the columns of $\matU$ and $\matV$ are left and right coordinated canonical correlation vectors not necessarily on the same phase are critical points, let $\alpha_{1},\dots,\alpha_{p}$ be some singular values of $\matR$,
and let $\tilde{\u}_{1},\dots,\tilde{\u}_{p},\tilde{\v}_{1},\dots,\tilde{\v}_{p}$
be the corresponding left and right singular vectors not necessarily
on the same phase. Writing $\tilde{\u}_{1},\dots,\tilde{\u}_{p}$
as the columns of $\tilde{\matU}$ and and $\tilde{\v}_{1},\dots,\tilde{\v}_{p}$
as the columns of $\tilde{\matV}$, and defining $\matU=\Sigma_{\x\x}^{-\nicehalf}\tilde{\matU},\matV=\Sigma_{\x\x}^{-\nicehalf}\tilde{\matV}$,
the following two equations hold:
\[
\Sigma_{\x\y}\matV=\Sigma_{\x\x}\matU\matA,\quad  \Sigma_{\x\y}^{\T}\matU=\Sigma_{\y\y}\matV\matA
\]
where $\matA\coloneqq\diag{\beta_{1},\dots,\beta_{p}}$ such that
$|\beta_{i}|=\alpha_{i}$ for $i=1,...,p$. Letting $\matZ=(\matU,\matV)\in\elpCCA$,
we have 
$$
\grad{\fcca}  =  -\left[\begin{array}{c}
\matU\matA\matN-\matU\sym{\left(\matU^{\T}\Sigma_{\x\x}\matU\matA\matN\right)}\\
\matV\matA\matN-\matV\sym{\left(\matV^{\T}\Sigma_{\y\y}\matV\matA\matN\right)}
\end{array}\right]=\left[\begin{array}{c}
\mat 0_{d_{\x}\times p}\\
\mat 0_{d_{\y}\times p}
\end{array}\right]
$$
where the last line  is true since $\matA\matN$ is diagonal.

To show the other side, note that if the Riemannian gradient nullifies, then 
$$
\left[\begin{array}{c}
\left(\matI_{d_{\x}}-\matU\matU^{\T}\Sigma_{\x\x}\right)\Sigma_{\x\x}^{-1}\Sigma_{\x\y}\matV\matN+\matU\skew{(\matU^{\T}\Sigma_{\x\y}\matV\matN)}\\
\left(\matI_{d_{\y}}-\matV\matV^{\T}\Sigma_{\y\y}\right)\Sigma_{\y\y}^{-1}\Sigma_{\x\y}^{\T}\matU\matN+\matV\skew{(\matV^{\T}\Sigma_{\x\y}^{\T}\matU\matN)}
\end{array}\right] = \left[\begin{array}{c}
\mat 0_{d_{\x}\times p}\\
\mat 0_{d_{\y}\times p}
\end{array}\right].
$$
By using similar reasoning as in \cite[Subsection 4.8.2]{AMS09},
$$(\matI_{d_{\x}}-\matU\matU^{\T}\Sigma_{\x\x})\Sigma_{\x\x}^{-1}\Sigma_{\x\y}\matV\matN,$$
belongs to the orthogonal compliment of the column space of $\matU$
(with respect to the matrix $\Sigma_{\x\x}$), and $\matU\skew{(\matU^{\T}\Sigma_{\x\y}\matV\matN)}$
belongs to the column space of $\matU$. The same is true for $(\matI_{d_{\y}}-\matV\matV^{\T}\Sigma_{\y\y})\Sigma_{\y\y}^{-1}\Sigma_{\x\y}^{\T}\matU\matN$ and $\matV\skew{(\matV^{\T}\Sigma_{\x\y}^{\T}\matU\matN)}$ with respect to $\matV$ and $\Sigma_{\y\y}$. Thus, we get that the gradient vanishes if and only if the aforementioned four factors vanish.

From $\left(\matI_{d_{\x}}-\matU\matU^{\T}\Sigma_{\x\x}\right)\Sigma_{\x\x}^{-1}\Sigma_{\x\y}\matV\matN=\mat 0_{d_{\x}\times p}$ we get $\Sigma_{\x\x}^{-1}\Sigma_{\x\y}\matV=\matU\left(\matU^{\T}\Sigma_{\x\y}\matV\right)$,
since $\matN$ is an invertible matrix. Also, since $\matU\in\stiefel_{\Sigma_{\x\x}}(p,d_{\x})$
it is a full (column) rank matrix then $\matU\skew{(\matU^{\T}\Sigma_{\x\y}\matV\matN)}$ vanishes if and only if $\skew{(\matU^{\T}\Sigma_{\x\y}\matV\matN)}$ vanishes, which leads to $\left(\matU^{\T}\Sigma_{\x\y}\matV\right)\matN=\matN\left(\matU^{\T}\Sigma_{\x\y}\matV\right)$.
This implies that $\matU^{\T}\Sigma_{\x\y}\matV$ is diagonal because
any rectangular matrix that commutes with a diagonal matrix with distinct
entries is diagonal. Thus, we have 
\begin{equation}
\Sigma_{\x\x}^{-1}\Sigma_{\x\y}\matV=\matU\mat D\ ,\label{eq:diag1}
\end{equation}
where $\mat D=\matU^{\T}\Sigma_{\x\y}\matV$ is diagonal. Similarly,
we have
\begin{equation}
\Sigma_{\y\y}^{-1}\Sigma_{\x\y}^{\T}\matU=\matV\mat D\ ,\label{eq:diag2}
\end{equation}
where $\mat D=\matV^{\T}\Sigma_{\x\y}^{\T}\matU=\matU^{\T}\Sigma_{\x\y}\matV$
is diagonal. From Eq.~(\ref{eq:diag1}) and Eq.~(\ref{eq:diag2})
we get
\[
\matR\Sigma_{\y\y}^{\nicehalf}\matV=\Sigma_{\x\x}^{\nicehalf}\matU\mat D,\ \matR^{\T}\Sigma_{\x\x}^{\nicehalf}\matU=\Sigma_{\y\y}^{\nicehalf}\matV\mat D.
\]
This implies that the columns of $(\tilde{\matU},\tilde{\matV})=(\Sigma_{\x\x}^{\nicehalf}\matU,\Sigma_{\y\y}^{\nicehalf}\matV)$
are some $p$ left and right singular vectors of $\matR$ not necessarily
on the same phase, but corresponding to the same singular values.

Finally, to identify the optimal solutions, note that at the critical points the objective function is the sum of the canonical correlations multiplied by a diagonal element of $\matN$, and a sign corresponding to the canonical correlation vectors in the columns of $\matU$ and $\matV$ and the correspondence of their phase. Thus, the optimal solutions that minimize $\fcca$ on $\elpCCA$
are $\matZ=(\matU,\matV)\in\elpCCA$ such that the columns of
$\matU$ and $\matV$ correspond to the top $p$-canonical correlations on the \emph{same phase}. Otherwise, we can increase the value of the objective function either by flipping the sign of one of the vectors or by
replacing a canonical vector with another that corresponds to a smaller canonical correlation. Moreover, if we assume that $\sigma_{1}>\sigma_{2}>...>\sigma_{p+1}\geq0$,
then for the aforementioned $\matZ=(\matU,\matV)\in\elpCCA$, the
columns of $(\tilde{\matU},\tilde{\matV})$ belong each to a one dimensional
singular left and right space and so do the columns of the corresponding
$(\matU,\matV)$, i.e., unique solution up the the signs of the columns
of $(\matU,\matV)$. In the case where some $\sigma_{i}=\sigma_{j}$
for $1\leq i,j\leq p$, permutations of the columns of $\matU$
and $\matV$ associated with $\sigma_{i}$ keep the solution optimal making it non-unique.
\end{proof}

\subsection{Proof of Theorem \ref{thm:stabilityCCA}}\label{subsec:proofstabilityCCA}

\begin{proof}
To prove the asymptotic stability of $\matZ=(\matU,\matV)$
that minimize $\fcca$ on $\elpCCA$, we use \cite[Proposition 4.4.2]{AMS09}. Recall from Theorem \ref{thm:criticalpointofCCA} that 
$\matZ$ that solve Problem~(\ref{eq:costriemannianccaBrockett})
are unique up the the signs of the columns of $\matU$ and $\matV$,
making these points isolated global (and consequently local) minimizers
of $\fcca$ on $\elpCCA$. According to \cite[Proposition 4.4.2]{AMS09},
such points $\matZ$ are asymptotically stable.

Suppose $\matZ$ is a critical point of $\fcca$
on $\elpCCA$ which is not a local minimum. Then, there exists compact
neighborhoods with either no other critical points, if there are no
multiplicities of the canonical correlations, or where all other critical
point achieve the same value for the cost function, if there are multiplicities. Thus, according to \cite[Proposition 4.4.1]{AMS09},
such $\matZ$ are unstable. 
\end{proof}

\subsection{Proof of Theorem~\ref{thm:CCAHessiantheorem}}\label{subsec:proofCCAHessiantheorem}

\begin{proof}
In order to bound the condition number of the Riemannian Hessian at
$\matZ^{\star}\in\elpCCA$, we need to bound its maximal and minimal
eigenvalues. Thus, to prove the
theorem we analyze the eigenvalues of the Riemannian Hessian at some
critical point $\matZ\in\elpCCA$ (in particular at $\matZ^{\star}\in\elpCCA$)
using the Courant-Fischer Theorem (also called the minimax principle,
see \cite[Chapter 1, Section 6.10]{kato2013perturbation}) for the
compact self-adjoint linear operator $\hess{\fcca[\cdot]}:T_{\matZ}\elpCCA\to T_{\matZ}\elpCCA$
over the finite dimensional vector space $T_{\matZ}\elpCCA$:
\begin{eqnarray}
\lambda_{k}(\hess{\fcca}) & = & \min_{U,\dim(U)=k-1}\max_{\mat 0\neq\xi_{\matZ}\in U^{\perp}}R(\xi_{\matZ}),\label{eq:CCAopt1a}
\end{eqnarray}
\begin{eqnarray}
\lambda_{k}(\hess{\fcca}) & = & \max_{U,\dim(U)=k}\min_{\mat 0\neq\xi_{\matZ}\in U}R(\xi_{\matZ}),\label{eq:CCAopt2a}
\end{eqnarray}
where 
\[
R(\xi_{\matZ})\coloneqq\frac{g_{\matZ}(\xi_{\matZ},\hess{\fcca[\xi_{\matZ}]})}{g_{\matZ}(\xi_{\matZ},\xi_{\matZ})},
\]
In the above, $\lambda_{k}(\hess{\fcca})$ is the $k$th largest eigenvalue
(i.e., eigenvalues are ordered in a descending order) of $\hess{\fcca}$,
and $U$ is a linear subspace of $T_{\matZ}\elpCCA$. In particular,
the maximal and minimal eigenvalues are given by the formulas
\begin{equation}
\lambda_{\max}(\hess{\fcca})=\max_{\mat 0\neq\xi_{\matZ}\in T_{\matZ}\elpCCA}R(\xi_{\matZ}),\label{eq:CCAopt1}
\end{equation}
 
\begin{equation}
\lambda_{\min}(\hess{\fcca})=\min_{\mat 0\neq\xi_{\matZ}\in T_{\matZ}\elpCCA}R(\xi_{\matZ}).\label{eq:CCAopt2}
\end{equation}

We begin by simplifying the quotient $R(\xi_{\matZ})$. Recall that
any critical point of $\justfcca$ is a matrix $\matZ=(\matU,\matV)\in\elpCCA$
such that the columns of $(\tilde{\matU},\tilde{\matV})=(\Sigma_{\x\x}^{\nicehalf}\matU,\Sigma_{\y\y}^{\nicehalf}\matV)$
are $p$ left and right singular vectors (not necessarily on the same
phase) of the matrix $\matR\coloneqq\Sigma_{\x\x}^{-\nicehalf}\Sigma_{\x\y}\Sigma_{\y\y}^{-\nicehalf}$ (see Theorem \ref{thm:criticalpointofCCA}).
Let $\alpha_{1},\dots,\alpha_{p}$ be some singular values of the
matrix $\matR$, and let $\tilde{\u}_{1},\dots,\tilde{\u}_{p},\tilde{\v}_{1},\dots,\tilde{\v}_{p}$
be the corresponding left and right singular vectors (not necessarily
on the same phase). Writing $\tilde{\u}_{1},\dots,\tilde{\u}_{p}$
as the columns of $\tilde{\matU}$ and and $\tilde{\v}_{1},\dots,\tilde{\v}_{p}$
as the columns of $\tilde{\matV}$, and defining $\matU=\Sigma_{\x\x}^{-\frac{1}{2}}\tilde{\matU},\matV=\Sigma_{\x\x}^{-\frac{1}{2}}\tilde{\matV}$,
the following two equations hold:
\begin{equation}
\Sigma_{\x\y}\matV=\Sigma_{\x\x}\matU\matA\label{eq:CCAproofform1}
\end{equation}
\begin{equation}
\Sigma_{\x\y}^{\T}\matU=\Sigma_{\y\y}\matV\matA\label{eq:CCAproofform2}
\end{equation}
where $\matA\coloneqq\diag{\beta_{1},\dots,\beta_{p}}$ such that
$|\beta_{i}|=\alpha_{i}$ for $i=1,...,p$. Letting $\matZ=(\matU,\matV)\in\elpCCA$,
plugging in the ambient coordinates formula for the Riemannian Hessian
(Eq. (\ref{eq:HessCCA})), the Riemannian gradient nullifies (see
Theorem \ref{thm:criticalpointofCCA}) and we have {\footnotesize{}
\begin{multline}
\hess{\fcca}[\xi_{\matZ}] = \mat{\Pi}_{\matZ}(\matM^{-1}_{\matZ}[-\sigmacca\xi_{\matZ}\matN+\Sigma\left[\begin{array}{c}
\xi_{\matU}\matU^{\T}\Sigma_{\x\y}\matV\matN\\
\xi_{\matV}\matV^{\T}\Sigma_{\x\y}^{\T}\matU\matN
\end{array}\right]\\+\Sigma\left[\begin{array}{cc}
\xi_{\matU}\matU^{\T}\\
 & \xi_{\matV}\matV^{\T}
\end{array}\right]\matM_{\matZ}\grad{\fcca}]) 
  =  \mat{\Pi}_{\matZ}\left(\matM^{-1}_{\matZ}\left[-\sigmacca\xi_{\matZ}+\Sigma\xi_{\matZ}\matA\right]\matN\right).\label{eq:hessatcritical}
\end{multline}
}{\footnotesize\par}

Plugging in the formula for the Riemannian Hessian at a critical point
(Eq. (\ref{eq:hessatcritical})), the quotient $R(\xi_{\matZ})$ is reduced to
to 
\[
R(\xi_{\matZ})=\frac{\Trace{\xi_{\matZ}^{\T}\matM_{\matZ}\mat{\Pi}_{\matZ}\left(\matM^{-1}_{\matZ}\left[-\sigmacca\xi_{\matZ}+\Sigma\xi_{\matZ}\matA\right]\matN\right)}}{\Trace{\xi_{\matZ}^{\T}\matM_{\matZ}\xi_{\matZ}}}.
\]
Now, using the fact that the projection to the tangent space is self-adjoint
with respect to the Riemannian metric and that for any $\xi_{\matZ}\in T_{\matZ}\elpCCA$
we have $\mat{\Pi}_{\matZ}\left(\xi_{\matZ}\right)=\xi_{\matZ}$,
we further see that 
\[
\frac{\Trace{\xi_{\matZ}^{\T}\matM_{\matZ}\mat{\Pi}_{\matZ}\left(\matM^{-1}_{\matZ}\left[-\sigmacca\xi_{\matZ}+\Sigma\xi_{\matZ}\matA\right]\matN\right)}}{\Trace{\xi_{\matZ}^{\T}\matM_{\matZ}\xi_{\matZ}}}=\frac{\Trace{\xi_{\matZ}^{\T}\left(-\sigmacca\xi_{\matZ}+\Sigma\xi_{\matZ}\matA\right)\matN}}{\Trace{\xi_{\matZ}^{\T}\matM_{\matZ}\xi_{\matZ}}}.
\]
Obviously, we can also write 
\begin{multline}
\frac{\Trace{\xi_{\matZ}^{\T}\left(-\sigmacca\xi_{\matZ}+\Sigma\xi_{\matZ}\matA\right)\matN}}{\Trace{\xi_{\matZ}^{\T}\matM_{\matZ}\xi_{\matZ}}} \\ = \frac{\Trace{\xi_{\matZ}^{\T}\left(-\sigmacca\xi_{\matZ}+\Sigma\xi_{\matZ}\matA\right)\matN}}{\Trace{\xi_{\matZ}^{\T}\Sigma\xi_{\matZ}}}\cdot\frac{\Trace{\xi_{\matZ}^{\T}\Sigma\xi_{\matZ}}}{\Trace{\xi_{\matZ}^{\T}\matM_{\matZ}\xi_{\matZ}}}.\label{eq:CCAproofRayleigh}
\end{multline}
Using Eq. \eqref{eq:CCAproofRayleigh}, a simplified
form of the quotient $R(\xi_{\matZ})$, we can estimate upper and lower bounds on $R(\xi_{\matZ})$ where $\mat 0\neq\xi_{\matZ}\in T_{\matZ}\elpCCA$
in order to bound the condition number of the Riemannian Hessian at
$\matZ^{\star}\in\elpCCA$. Since for $\xi_{\matZ}\neq\mat 0$
the term $\Trace{\xi_{\matZ}^{\T}\Sigma\xi_{\matZ}}/\Trace{\xi_{\matZ}^{\T}\matM_{\matZ}\xi_{\matZ}}>0$,
then the upper and lower bounds on
\begin{equation}
\frac{\Trace{\xi_{\matZ}^{\T}\left(-\sigmacca\xi_{\matZ}+\Sigma\xi_{\matZ}\matA\right)\matN}}{\Trace{\xi_{\matZ}^{\T}\Sigma\xi_{\matZ}}},\label{eq:ccarayleighforMeqB}
\end{equation}
together with the upper and lower bounds of $\Trace{\xi_{\matZ}^{\T}\Sigma\xi_{\matZ}}/\Trace{\xi_{\matZ}^{\T}\matM_{\matZ}\xi_{\matZ}}>0$
bound the condition number of the Riemannian Hessian at $\matZ^{\star}\in\elpCCA$.

We begin by estimating the term $\Trace{\xi_{\matZ}^{\T}\Sigma\xi_{\matZ}}/\Trace{\xi_{\matZ}^{\T}\matM_{\matZ}\xi_{\matZ}}>0$.
We use the vectorization operator, and the Kronecker product to rewrite
it in the following form
\begin{equation}
\frac{\Trace{\xi_{\matZ}^{\T}\Sigma\xi_{\matZ}}}{\Trace{\xi_{\matZ}^{\T}\matM_{\matZ}\xi_{\matZ}}}=\frac{\vectorization{\xi_{\matZ}}^{\T}\left(\matI_{p}\otimes\Sigma\right)\vectorization{\xi_{\matZ}}}{\vectorization{\xi_{\matZ}}^{\T}\left(\matI_{p}\otimes\matM_{\matZ}\right)\vectorization{\xi_{\matZ}}}.\label{eq:CCAhessiansecondRayleighquotient}
\end{equation}
Eq. (\ref{eq:CCAhessiansecondRayleighquotient}) is the generalized
Rayleigh quotient for the matrix pencil $\left(\matI_{p}\otimes\Sigma,\matI_{p}\otimes\matM_{\matZ}\right)$.
Note that $\matI_{p}\otimes\Sigma$ and $\matI_{p}\otimes\matM_{\matZ}$
are both SPD matrices, thus the generalized eigenvalues of the matrix
pencil $\left(\matI_{p}\otimes\Sigma,\matI_{p}\otimes\matM_{\matZ}\right)$
are equivalent to the eigenvalues of the matrix $\left(\matI_{p}\otimes\matM_{\matZ}\right)^{-1}\left(\matI_{p}\otimes\Sigma\right)=\matI_{p}\otimes\matM^{-1}_{\matZ}\Sigma$.
According to \cite[Section 2]{minka2000old} the eigenvalues $\matI_{p}\otimes\matM^{-1}_{\matZ}\Sigma$
are $p$ copies of each of the eigenvalues of $\matM^{-1}_{\matZ}\Sigma$.
Thus, the maximal and minimal eigenvalues of the matrix pencil $\left(\matI_{p}\otimes\Sigma,\matI_{p}\otimes\matM_{\matZ}\right)$
denoted by $\tilde{\lambda}_{\text{max}}$ and $\tilde{\lambda}_{\min}$
are equivalent to the maximal and minimal generalized eigenvalues
of the matrix pencil $\left(\Sigma,\matM_{\matZ}\right)$, and so is the corresponding
condition number 
\[
\kappa\left(\matI_{p}\otimes\Sigma,\matI_{p}\otimes\matM_{\matZ}\right)=\frac{\tilde{\lambda}_{\text{max}}}{\tilde{\lambda}_{\min}}=\kappa\left(\Sigma,\matM_{\matZ}\right).
\]

Recall the definition of the generalized eigenvalues. The generalized eigenvalues of the matrix pencil $(\matA,\matB)$, where $\matA\in\R^{d\times d}$ and
$\matB\in\R^{d\times d}$ is a symmetric positive semi-definite matrix
such that $\ker(\matB)\subseteq\ker(\matA)$, are defined as follows:
if for $\lambda\in\R$ and $\v\notin\ker(\matB)$ it holds that $\matA\v=\lambda\matB\v$
then $\lambda$ is a generalized eigenvalue and $\v$ is a generalized
eigenvector of the matrix pencil $(\matA,\matB)$. The generalized
eigenvalues are denoted by $\lambda_{1}(\matA,\matB)\geq\lambda_{2}(\matA,\matB)\geq\dots\geq\lambda_{\rank{\matB}}(\matA,\matB)$. Therefore, using the Courant-Fischer Theorem for the matrix pencil
$\left(\matI_{p}\otimes\Sigma,\matI_{p}\otimes\matM_{\matZ}\right)$
we have
\begin{eqnarray*}
\tilde{\lambda}_{\text{max}} & \coloneqq & \lambda_{\max}(\matI_{p}\otimes\Sigma,\matI_{p}\otimes\matM_{\matZ})\\
 & = & \max_{0\neq\xi_{\matZ}\in\R^{d\times p}}\frac{\vectorization{\xi_{\matZ}}^{\T}\left(\matI_{p}\otimes\Sigma\right)\vectorization{\xi_{\matZ}}}{\vectorization{\xi_{\matZ}}^{\T}\left(\matI_{p}\otimes\matM_{\matZ}\right)\vectorization{\xi_{\matZ}}}\\
 & \geq & \max_{0\neq\xi_{\matZ}\in T_{\matZ}\elpCCA}\frac{\vectorization{\xi_{\matZ}}^{\T}\left(\matI_{p}\otimes\Sigma\right)\vectorization{\xi_{\matZ}}}{\vectorization{\xi_{\matZ}}^{\T}\left(\matI_{p}\otimes\matM_{\matZ}\right)\vectorization{\xi_{\matZ}}}\\
 & = & \max_{0\neq\xi_{\matZ}\in T_{\matZ}\elpCCA}\frac{\Trace{\xi_{\matZ}^{\T}\Sigma\xi_{\matZ}}}{\Trace{\xi_{\matZ}^{\T}\matM_{\matZ}\xi_{\matZ}}}
\end{eqnarray*}
and
\begin{eqnarray*}
\tilde{\lambda}_{\min} & \coloneqq & \lambda_{\min}(\matI_{p}\otimes\Sigma,\matI_{p}\otimes\matM_{\matZ})\\
 & = & \min_{0\neq\xi_{\matZ}\in\R^{d\times p}}\frac{\vectorization{\xi_{\matZ}}^{\T}\left(\matI_{p}\otimes\Sigma\right)\vectorization{\xi_{\matZ}}}{\vectorization{\xi_{\matZ}}^{\T}\left(\matI_{p}\otimes\matM_{\matZ}\right)\vectorization{\xi_{\matZ}}}\\
 & \leq & \min_{0\neq\xi_{\matZ}\in T_{\matZ}\elpCCA}\frac{\vectorization{\xi_{\matZ}}^{\T}\left(\matI_{p}\otimes\Sigma\right)\vectorization{\xi_{\matZ}}}{\vectorization{\xi_{\matZ}}^{\T}\left(\matI_{p}\otimes\matM_{\matZ}\right)\vectorization{\xi_{\matZ}}}\\
 & = & \min_{0\neq\xi_{\matZ}\in T_{\matZ}\elpCCA}\frac{\Trace{\xi_{\matZ}^{\T}\Sigma\xi_{\matZ}}}{\Trace{\xi_{\matZ}^{\T}\matM_{\matZ}\xi_{\matZ}}}
\end{eqnarray*}

Next, we analyze Eq. (\ref{eq:ccarayleighforMeqB}). Recall that $\xi_{\matZ}=(\xi_{\matU},\xi_{\matV})\in T_{\matZ}\elpCCA$
where we have $\xi_{\matU}\in T_{\matU}\stiefel_{\Sigma_{\x\x}}(p,d_{\x})$
and $\xi_{\matV}\in T_{\matV}\stiefel_{\Sigma_{\y\y}}(p,d_{\y})$,
thus we can rewrite the tangent vectors in the following form:
\[
\xi_{\matU}=\matU\mat{\Omega}_{\matU}+\matU_{\Sigma_{\x\x}\perp}\mat K_{\matU},\ \xi_{\matV}=\matV\mat{\Omega}_{\matV}+\matV_{\Sigma_{\y\y}\perp}\mat K_{\matV},
\]
where $\matU_{\Sigma_{\x\x}\perp}$ is $\Sigma_{\x\x}$-orthogonal
to $\matU$ so that the union of the columns of $\matU$ and $\matU_{\Sigma_{\x\x}\perp}$
is a basis to $\R^{d_{\x}}$, and similarly $\matV_{\Sigma_{\y\y}\perp}$
is $\Sigma_{\y\y}$-orthogonal to $\matV$ so that the union of the
columns of $\matV$ and $\matV_{\Sigma_{\y\y}\perp}$ is a basis to
$\R^{d_{\y}}$, $\mat{\Omega}_{\matU}=-\mat{\Omega}_{\matU}^{\T}\in\R^{p\times p}$,
$\mat{\Omega}_{\matV}=-\mat{\Omega}_{\matV}^{\T}\in\R^{p\times p}$,
$\mat K_{\matU}\in\R^{\left(d_{\x}-p\right)\times p}$ and $\mat K_{\matV}\in\R^{\left(d_{\y}-p\right)\times p}$.
Note that we can always make the choice of the columns of $\matU_{\Sigma_{\x\x}\perp}$
and $\matV_{\Sigma_{\y\y}\perp}$ to be such that $(\Sigma_{\x\x}^{\nicehalf}\matU_{\Sigma_{\x\x}\perp},\Sigma_{\y\y}^{\nicehalf}\matV_{\Sigma_{\y\y}\perp})$
are some $\min\left\{ d_{\x}-p,d_{\y}-p\right\} $ left and right
singular vectors not necessarily on the same phase of the matrix $\matR$
belonging to the same singular values. Without loss of generality
suppose $d_{\x}\geq d_{\y}$. With this choice we have
\[
\Sigma_{\x\y}\matV_{\Sigma_{\y\y}\perp}=\Sigma_{\x\x}\matU_{\Sigma_{\x\x}\perp}\tilde{\matA},\ \Sigma_{\x\y}^{\T}\matU_{\Sigma_{\x\x}\perp}=\Sigma_{\y\y}\matV_{\Sigma_{\y\y}\perp}\tilde{\matA}^{\T},
\]
and
\[
\matU_{\Sigma_{\x\x}\perp}^{\T}\Sigma_{\x\y}\matV_{\Sigma_{\y\y}\perp}=\tilde{\matA}\in\R^{\left(d_{\x}-p\right)\times\left(d_{\y}-p\right)},\ \matV_{\Sigma_{\y\y}\perp}^{\T}\Sigma_{\x\y}^{\T}\matU_{\Sigma_{\x\x}\perp}=\tilde{\matA}^{\T}\in\R^{\left(d_{\y}-p\right)\times\left(d_{\x}-p\right)},
\]
where $\tilde{\matA}$ is a diagonal matrix (but not necessarily a
square matrix), with the corresponding values on the diagonal $\beta_{p+1},...,\beta_{d_{\y}}$,
which satisfy $|\beta_{i}|=\alpha_{i}$ for $i=p+1,...,d_{\y}$.

Now, 
\[
\xi_{\matZ}^{\T}\Sigma\xi_{\matZ}=\left[\begin{array}{cc}
\mat{\Omega}_{\matU}^{\T} & \mat{\Omega}_{\matV}^{\T}\end{array}\right]\left[\begin{array}{c}
\mat{\Omega}_{\matU}\\
\mat{\Omega}_{\matV}
\end{array}\right]+\left[\begin{array}{cc}
\mat K_{\matU}^{\T} & \mat K_{\matV}^{\T}\end{array}\right]\left[\begin{array}{c}
\mat K_{\matU}\\
\mat K_{\matV}
\end{array}\right]
\]
and 
\[
\xi_{\matZ}^{\T}\sigmacca\xi_{\matZ}=\left[\begin{array}{cc}
\mat{\Omega}_{\matU}^{\T} & \mat{\Omega}_{\matV}^{\T}\end{array}\right]\left[\begin{array}{cc}
 & \matA\\
\matA
\end{array}\right]\left[\begin{array}{c}
\mat{\Omega}_{\matU}\\
\mat{\Omega}_{\matV}
\end{array}\right]+\left[\begin{array}{cc}
\mat K_{\matU}^{\T} & \mat K_{\matV}^{\T}\end{array}\right]\left[\begin{array}{cc}
 & \tilde{\matA}\\
\tilde{\matA}^{\T}
\end{array}\right]\left[\begin{array}{c}
\mat K_{\matU}\\
\mat K_{\matV}
\end{array}\right].
\]
Let 
\[
\m_{\matZ}\coloneqq\vectorization{\left[\begin{array}{c}
\mat{\Omega}_{\matU}\\
\mat{\Omega}_{\matV}
\end{array}\right]},\quad\k_{\matZ}\coloneqq\vectorization{\left[\begin{array}{c}
\mat K_{\matU}\\
\mat K_{\matV}
\end{array}\right]}
\]
Then,
\[
\Trace{\xi_{\matZ}^{\T}\Sigma\xi_{\matZ}}=\m_{\matZ}^{\T}\m_{\matZ}+\k_{\matZ}^{\T}\k_{\matZ},
\]
{\footnotesize{}
\begin{multline*}
\Trace{-\left[\begin{array}{cc}
\mat{\Omega}_{\matU}^{\T} & \mat{\Omega}_{\matV}^{\T}\end{array}\right]\left[\begin{array}{cc}
 & \matA\\
\matA
\end{array}\right]\left[\begin{array}{c}
\mat{\Omega}_{\matU}\\
\mat{\Omega}_{\matV}
\end{array}\right]\matN+\left[\begin{array}{cc}
\mat{\Omega}_{\matU}^{\T} & \mat{\Omega}_{\matV}^{\T}\end{array}\right]\left[\begin{array}{c}
\mat{\Omega}_{\matU}\\
\mat{\Omega}_{\matV}
\end{array}\right]\matA\matN}=\\ \m_{\matZ}^{\T}\left(\matA\matN\otimes\matI_{2p}-\matN\otimes\left[\begin{array}{cc}
 & \matA\\
\matA
\end{array}\right]\right)\m_{\matZ}
\end{multline*}
\begin{multline*}
\Trace{-\left[\begin{array}{cc}
\mat K_{\matU}^{\T} & \mat K_{\matV}^{\T}\end{array}\right]\left[\begin{array}{cc}
 & \tilde{\matA}\\
\tilde{\matA}^{\T}
\end{array}\right]\left[\begin{array}{c}
\mat K_{\matU}\\
\mat K_{\matV}
\end{array}\right]\matN+\left[\begin{array}{cc}
\mat K_{\matU}^{\T} & \mat K_{\matV}^{\T}\end{array}\right]\left[\begin{array}{c}
\mat K_{\matU}\\
\mat K_{\matV}
\end{array}\right]\matA\matN}=\\ \k_{\matZ}^{\T}\left(\matA\matN\otimes\matI_{2p}-\matN\otimes\left[\begin{array}{cc}
 & \tilde{\matA}\\
\tilde{\matA}^{\T}
\end{array}\right]\right)\k_{\matZ}
\end{multline*}
}thus,{\footnotesize{}
\begin{multline*}
\frac{\Trace{\xi_{\matZ}^{\T}\left(-\sigmacca\xi_{\matZ}+\Sigma\xi_{\matZ}\matA\right)\matN}}{\Trace{\xi_{\matZ}^{\T}\Sigma\xi_{\matZ}}}\\=\frac{\m_{\matZ}^{\T}\left(\matA\matN\otimes\matI_{2p}-\matN\otimes\left[\begin{array}{cc}
 & \matA\\
\matA
\end{array}\right]\right)\m_{\matZ}+\k_{\matZ}^{\T}\left(\matA\matN\otimes\matI_{2p}-\matN\otimes\left[\begin{array}{cc}
 & \tilde{\matA}\\
\tilde{\matA}^{\T}
\end{array}\right]\right)\k_{\matZ}}{\m_{\matZ}^{\T}\m_{\matZ}+\k_{\matZ}^{\T}\k_{\matZ}}.
\end{multline*}
}{\footnotesize\par}

Recalling that $\mat{\Omega}_{\matU}=-\mat{\Omega}_{\matU}^{\T}$
and $\mat{\Omega}_{\matV}=-\mat{\Omega}_{\matV}^{\T}$, and both are
real matrices (so the elements of the main diagonals are $0$), we
have 
\begin{equation}
\m_{\matZ}^{\T}\m_{\matZ}=2\left(\sum_{1\leq j<i\leq p}\left(\mat{\Omega}_{\matU}\right)_{ij}^{2}+\sum_{1\leq j<i\leq p}\left(\mat{\Omega}_{\matV}\right)_{ij}^{2}\right),\label{eq:mzmz}
\end{equation}
 and
\begin{multline}
\m_{\matZ}^{\T}\left(\matA\matN\otimes\matI_{2p}-\matN\otimes\left[\begin{array}{cc}
 & \matA\\
\matA
\end{array}\right]\right)\m_{\matZ}\\=\sum_{1\leq j<i\leq p}\left[\left(\mu_{i}\beta_{i}+\mu_{j}\beta_{j}\right)\left(\left(\mat{\Omega}_{\matU}\right)_{ij}^{2}+\left(\mat{\Omega}_{\matV}\right)_{ij}^{2}\right)-2\left(\beta_{i}\mu_{j}+\beta_{j}\mu_{i}\right)\left(\mat{\Omega}_{\matU}\right)_{ij}\left(\mat{\Omega}_{\matV}\right)_{ij}\right].\label{eq:mzAmz}
\end{multline}
Thus, only the $p(p-1)/2$ entries below the diagonal of $\mat{\Omega}_{\matU}$
and the $p(p-1)/2$ entries below the diagonal of $\mat{\Omega}_{\matV}$
determine the values of Eq. (\ref{eq:mzmz}) and Eq. (\ref{eq:mzAmz}).
Let us now denote by $\tilde{\m}_{\matZ}$ the column stack of $(\mat{\Omega}_{\matU},\mat{\Omega}_{\matV})$,
but only with the the subdiagonal entries of $\mat{\Omega}_{\matU}$
and of $\mat{\Omega}_{\matV}$ (i.e., $\m_{\matZ}$ ``purged'' of
the superdiagonal elements). We then have, 
\[
\m_{\matZ}^{\T}\m_{\matZ}=2\tilde{\m}_{\matZ}^{\T}\tilde{\m}_{\matZ}
\]
and
\[
\m_{\matZ}^{\T}\left(\matA\matN\otimes\matI_{2p}-\matN\otimes\left[\begin{array}{cc}
 & \matA\\
\matA
\end{array}\right]\right)\m_{\matZ}=\tilde{\m}_{\matZ}^{\T}\Psi\tilde{\m}_{\matZ}
\]
where $\mat{\Psi}\in\R^{p(p-1)\times p(p-1)}$ is a matrix defined
as follows: $\Psi$ is a block diagonal matrix, where the blocks are
of descending order: $2\left(p-1\right),2\left(p-2\right),...,2=2$.
The $j$th ($1\leq j\leq p-1$) block, denoted by $\mat{\Psi}_{j}$,
is of the order $2\left(p-j\right)$ and has following form:
\[
\mat{\Psi}_{j}\coloneqq\left[\begin{array}{cc}
\matD_{j} & -\mat T_{j}\\
-\mat T_{j} & \matD_{j}
\end{array}\right],
\]
where 
\[
\matD_{j}=\diag{\mu_{j+1}\beta_{j+1}+\mu_{j}\beta_{j},\mu_{j+2}\beta_{j+2}+\mu_{j}\beta_{j},...,\mu_{p}\beta_{p}+\mu_{j}\beta_{j}},
\]
and 
\[
\mat T_{j}=\diag{\beta_{j+1}\mu_{j}+\beta_{j}\mu_{j+1},\beta_{j+2}\mu_{j}+\beta_{j}\mu_{j+2},...,\beta_{p}\mu_{j}+\beta_{j}\mu_{p}}.
\]
We make the following change of variables: $\tilde{\d}_{\matZ}\coloneqq\sqrt{2}\tilde{\m}_{\matZ}$. Finally, Eq. (\ref{eq:ccarayleighforMeqB}) is rewritten in the
following way 
\begin{multline}
\frac{\Trace{\xi_{\matZ}^{\T}\left(-\sigmacca\xi_{\matZ}+\Sigma\xi_{\matZ}\matA\right)\matN}}{\Trace{\xi_{\matZ}^{\T}\Sigma\xi_{\matZ}}}=\\\frac{\left[\begin{array}{cc}
\tilde{\d}_{\matZ}^{\T} & \k_{\matZ}^{\T}\end{array}\right]\blockdiag{\frac{1}{2}\mat{\Psi},\matA\matN\otimes\matI_{d-2p}-\matN\otimes\left[\begin{array}{cc}
 & \tilde{\matA}\\
\tilde{\matA}^{\T}
\end{array}\right]}\left[\begin{array}{c}
\tilde{\d}_{\matZ}\\
\k_{\matZ}
\end{array}\right]}{\left[\begin{array}{cc}
\tilde{\d}_{\matZ}^{\T} & \k_{\matZ}^{\T}\end{array}\right]\left[\begin{array}{c}
\tilde{\d}_{\matZ}\\
\k_{\matZ}
\end{array}\right]}.\label{eq:RayleighCCA}
\end{multline}

Note that the mapping $\varphi(\cdot):T_{\matZ}\elpCCA\to\R^{pd-p(p+1)}$
defined by 
\[
\varphi(\xi_{\matZ})\coloneqq\left[\begin{array}{c}
\tilde{\d}_{\matZ}\\
\k_{\matZ}
\end{array}\right],
\]
is a coordinate chart of the elements of $T_{\matZ}\elpCCA$, since
$\varphi(\cdot)$ is a bijection (one-to-one correspondence) of the
elements of $T_{\matZ}\elpCCA$ onto $\R^{pd-p(p+1)}$. Indeed, $\k_{\matZ}$
is a column stack of $(\mat K_{\matU},\mat K_{\matV})$, thus we can
retract the matrices $\mat K_{\matU}$ and $\mat K_{\matV}$. Similarly
$\tilde{\d}_{\matZ}$ is proportional to $\tilde{\m}_{\matZ}$ which
is a column stack of $(\mat{\Omega}_{\matU},\mat{\Omega}_{\matV})$,
but only with the the subdiagonal entries of $\mat{\Omega}_{\matU}$
and of $\mat{\Omega}_{\matV}$. Since $\mat{\Omega}_{\matU}$ and
$\mat{\Omega}_{\matV}$ are skew-symmetric matrices, we can retract
$\mat{\Omega}_{\matU}$ and $\mat{\Omega}_{\matV}$. With the matrices
$\mat K_{\matU},\mat K_{\matV},\mat{\Omega}_{\matU}$ and $\mat{\Omega}_{\matV}$
at hand, we can fully retract $\xi_{\matZ}=(\xi_{\matU},\xi_{\matV})$.

The eigenvalues and corresponding eigenvectors of any linear operator
over a finite dimensional vector space do not depend on the choice
of coordinate chart and basis, thus the eigenvalues and eigenvectors
of $\hess{\fcca[\cdot]}:T_{\matZ}\elpCCA\to T_{\matZ}\elpCCA$ which
are computed using the Courant-Fischer Theorem for compact self-adjoint
linear operators (Eq. (\ref{eq:CCAopt1a}) and Eq. (\ref{eq:CCAopt2a})),
can be also computed by the Courant-Fischer Theorem for symmetric
matrices \cite[Theorem 4.2.6]{horn2012matrix} after applying $\varphi(\cdot)$.
In particular, Eq. (\ref{eq:RayleighCCA}) determines the signs of
the eigenvalues of the Riemannian Hessian at any $\matZ\in\elpCCA$
(in the special case $\matM_{\matZ}\coloneqq\Sigma$, the eigenvalues of Eq. (\ref{eq:CCAHessianoptimal})
are the eigenvalues of the Riemannian Hessian at a critical point $\matZ\in\elpCCA$),
and the bounds of Eq. (\ref{eq:RayleighCCA}) together with the bounds
of the term $\Trace{\xi_{\matZ}^{\T}\Sigma\xi_{\matZ}}/\Trace{\xi_{\matZ}^{\T}\matM_{\matZ}\xi_{\matZ}}$
bound the condition number of the Riemannian Hessian at $\matZ^{\star}\in\elpCCA$.

To that end, we perform the following computation. The righthand side
of Eq. (\ref{eq:RayleighCCA}) is a Rayleigh quotient, so according
to the Courant-Fischer Theorem for symmetric matrices the eigenvalues
of the $pd-p(p+1)\times pd-p(p+1)$ symmetric matrix 
\begin{equation}
\blockdiag{\frac{1}{2}\mat{\Psi},\matA\matN\otimes\matI_{d-2p}-\matN\otimes\left[\begin{array}{cc}
 & \tilde{\matA}\\
\tilde{\matA}^{\T}
\end{array}\right]},\label{eq:CCAHessianoptimal}
\end{equation}
are determined by critical values of Eq. (\ref{eq:RayleighCCA}).
The set of eigenvalues of the matrix in Eq.~(\ref{eq:CCAHessianoptimal})
is equal to the union of the set of eigenvalues of $\frac{1}{2}\mat{\Psi}$
and 
\[
\Phi\coloneqq\matA\matN\otimes\matI_{d-2p}-\matN\otimes\left[\begin{array}{cc}
 & \tilde{\matA}\\
\tilde{\matA}^{\T}
\end{array}\right].
\]
The matrix $\Phi$ is a $p(d-2p)\times p(d-2p)$ block diagonal matrix,
where all the blocks are $p\times p$, and the $j$th block is
\[
\mu_{j}\beta_{j}\matI_{d-2p}-\mu_{j}\left[\begin{array}{cc}
 & \tilde{\matA}\\
\tilde{\matA}^{\T}
\end{array}\right].
\]
Thus, the eigenvalues of $\Phi$ are $\mu_{j}\left(\beta_{j}\pm\beta_{i}\right)$
and $\mu_{j}\beta_{j}$ for $j=1,...,p$ and $i=p+1,...,d_{\y}$.
In summary, we have $p(d_{\y}-p)$ eigenvalues of the form $\mu_{j}\left(\beta_{j}+\beta_{i}\right)$,
similarly $p(d_{\y}-p)$ eigenvalues of the form $\mu_{j}\left(\beta_{j}-\beta_{i}\right)$,
and $p(d_{\x}-d_{\y})$ eigenvalues of the form $\mu_{j}\beta_{j}$.
From the definition of $\Psi$, we see that the eigenvalues of $\frac{1}{2}\mat{\Psi}$
are: $\frac{1}{2}\left[\left(\mu_{i}\beta_{i}+\mu_{j}\beta_{j}\right)\pm\left(\beta_{i}\mu_{j}+\beta_{j}\mu_{i}\right)\right]$
for $1\leq j<i\leq p$. These eigenvalues can also be rewritten as:
$\frac{1}{2}\left(\mu_{j}+\mu_{i}\right)\left(\beta_{j}+\beta_{i}\right)$
and $\frac{1}{2}\left(\mu_{j}-\mu_{i}\right)\left(\beta_{j}-\beta_{i}\right)$.

Now, we have all the eigenvalues of the matrix in Eq.~(\ref{eq:CCAHessianoptimal}):
$p(d_{\y}-p)$ eigenvalues of the form $\mu_{j}\left(\beta_{j}+\beta_{i}\right)$
where $j=1,...,p$ and $i=p+1,...,d_{\y}$, $p(d_{\y}-p)$ eigenvalues
of the form $\mu_{j}\left(\beta_{j}-\beta_{i}\right)$ where $j=1,...,p$
and $i=p+1,...,d_{\y}$, $p(d_{\x}-d_{\y})$ eigenvalues of the form
$\mu_{j}\beta_{j}$ where $j=1,...,p$, $p(p-1)/2$ eigenvalues of
the form $\frac{1}{2}\left(\mu_{j}+\mu_{i}\right)\left(\beta_{j}+\beta_{i}\right)$
where $1\leq j<i\leq p$, and $p(p-1)/2$ eigenvalues of the form
$\frac{1}{2}\left(\mu_{j}-\mu_{i}\right)\left(\beta_{j}-\beta_{i}\right)$
where $1\leq j<i\leq p$. 

Finally, we bound the condition number of the Riemannian Hessian at
$\matZ^{\star}=(\matU^{\star},\matV^{\star})\in\elpCCA$. In such
case, $\beta_{1}=\sigma_{1}>...>\beta_{p}=\sigma_{p}$. Without loss
of generality, we can always choose $\matU_{\Sigma_{\x\x}\perp}$
and $\matV_{\Sigma_{\y\y}\perp}$ such that $\beta_{p+1}=\sigma_{p+1}\geq...\geq\beta_{d_{\y}}=\sigma_{d_{\y}}$.
Then, we have that Eq. (\ref{eq:ccarayleighforMeqB}) is bounded by
the minimal and maximal eigenvalues of Eq. (\ref{eq:CCAHessianoptimal}).
Thus,
\begin{multline*}
0<\max_{\mat 0\neq\xi_{\matZ^{\star}}\in T_{\matZ^{\star}}\elpCCA}\frac{\Trace{\xi_{\matZ^{\star}}^{\T}\left(-\sigmacca\xi_{\matZ^{\star}}+\Sigma\xi_{\matZ^{\star}}\matA\right)\matN}}{\Trace{\xi_{\matZ^{\star}}^{\T}\Sigma\xi_{\matZ^{\star}}}}=\\ \max\left\{ \mu_{1}(\sigma_{1}+\sigma_{p+1}),\frac{1}{2}(\mu_{1}+\mu_{2})(\sigma_{1}+\sigma_{2})\right\} ,
\end{multline*}
and 

\begin{multline*}
\min_{\mat 0\neq\xi_{\matZ^{\star}}\in T_{\matZ^{\star}}\elpCCA}\frac{\Trace{\xi_{\matZ^{\star}}^{\T}\left(-\sigmacca\xi_{\matZ^{\star}}+\Sigma\xi_{\matZ^{\star}}\matA\right)\matN}}{\Trace{\xi_{\matZ^{\star}}^{\T}\Sigma\xi_{\matZ^{\star}}}}=\\ \min\left\{ \mu_{p}(\sigma_{p}-\sigma_{p+1}),\min_{1\leq j<i\leq p}\frac{1}{2}\left(\mu_{j}-\mu_{i}\right)\left(\sigma_{j}-\sigma_{i}\right)\right\} >0.
\end{multline*}
We use Eq.~(\ref{eq:CCAopt1}) and Eq.~(\ref{eq:CCAopt2}) to bound
the condition number of the Riemannian Hessian at $\matZ^{\star}=(\matU^{\star},\matV^{\star})\in\elpCCA$:
\begin{eqnarray*}
\lambda_{\max}(\hess{\justfcca(\matZ^{\star})}) & = & \max_{\mat 0\neq\xi_{\matZ^{\star}}\in T_{\matZ^{\star}}\elpCCA}\frac{g_{\matZ^{\star}}(\xi_{\matZ^{\star}},\hess{\justfcca(\matZ^{\star})[\xi_{\matZ^{\star}}]})}{g_{\matZ^{\star}}(\xi_{\matZ^{\star}},\xi_{\matZ^{\star}})}\\
 & \leq & \max\left\{ \mu_{1}(\sigma_{1}+\sigma_{p+1}),\frac{1}{2}(\mu_{1}+\mu_{2})(\sigma_{1}+\sigma_{2})\right\} \cdot\tilde{\lambda}_{\text{max}}\ ,
\end{eqnarray*}
and 
\begin{eqnarray*}
\lambda_{\min}(\hess{\justfcca(\matZ^{\star})}) & = & \min_{\mat 0\neq\xi_{\matZ^{\star}}\in T_{\matZ^{\star}}\elpCCA}\frac{g_{\matZ^{\star}}(\xi_{\matZ^{\star}},\hess{\justfcca(\matZ^{\star})[\xi_{\matZ^{\star}}]})}{g_{\matZ^{\star}}(\xi_{\matZ^{\star}},\xi_{\matZ^{\star}})}\\
 & \geq & \min\left\{ \mu_{p}(\sigma_{p}-\sigma_{p+1}),\min_{1\leq j<p}\frac{1}{2}\left(\mu_{j}-\mu_{j+1}\right)\left(\sigma_{j}-\sigma_{j+1}\right)\right\} \cdot\tilde{\lambda}_{\min}\ .
\end{eqnarray*}
Finally,
\[
\kappa(\hess{\justfcca(\matZ^{\star})})=\frac{\lambda_{\max}(\hess{\justfcca(\matZ^{\star})})}{\lambda_{\min}(\hess{\justfcca(\matZ^{\star})})}\leq\kappa_{{\bf CCA}}^{\star}\cdot\kappa\left(\Sigma,\matM_{\matZ}\right),
\]
where 
\[
\kappa_{{\bf CCA}}^{\star}=\frac{\max\left\{ \mu_{1}(\sigma_{1}+\sigma_{p+1}),\frac{1}{2}(\mu_{1}+\mu_{2})(\sigma_{1}+\sigma_{2})\right\} }{\min\left\{ \mu_{p}(\sigma_{p}-\sigma_{p+1}),\min_{1\leq j<p}\frac{1}{2}\left(\mu_{j}-\mu_{j+1}\right)\left(\sigma_{j}-\sigma_{j+1}\right)\right\} }.
\]
In the special case $\matM_{\matZ^{\star}}=\Sigma$, the bound on the condition number
of the Riemannian Hessian at $\matZ^{\star}=(\matU^{\star},\matV^{\star})\in\elpCCA$
is reduced to an equality 
\[
\kappa(\hess{\justfcca(\matZ^{\star})})=\kappa_{{\bf CCA}}^{\star}.
\]
\end{proof}

\subsection{Proof of Theorem \ref{thm:morestabilityCCA}}\label{subsec:proofmorestabilityCCA}

\begin{proof}
We prove that the global minimizer of
$\fcca$ subject to $\matZ\in\elpCCA$, denoted by $\matZ^{\star}$, is the
only local minimum of $\fcca$ (and it is also strict) and all other critical points are either
saddle points or strict local maximizers.

Eq. \eqref{eq:CCAproofRayleigh} helps to determine the signs of the eigenvalues of the Riemannian Hessian at
any critical point $\mat Z\in\elpCCA$ and in particular at $\matZ^{\star}$:
the matrices $\Sigma$ and $\matM_{\matZ}$ are both SPD matrices, therefore
for $\xi_{\matZ}\neq\mat 0$ the term $\Trace{\xi_{\matZ}^{\T}\Sigma\xi_{\matZ}}/\Trace{\xi_{\matZ}^{\T}\matM_{\matZ}\xi_{\matZ}}>0$,
thus only Eq. (\ref{eq:ccarayleighforMeqB}),
where $\mat 0\neq\xi_{\matZ}\in T_{\matZ}\elpCCA$ determines the
signs. In addition, at a critical point $\matZ\in\elpCCA$, Eq. (\ref{eq:ccarayleighforMeqB})
equals to the quotient $R(\xi_{\matZ})$ for the choice $\matM_{\matZ}\coloneqq\Sigma$,
since in \cite[Proposition 5.5.6 and Eq. 5.35]{AMS09} it is shown that at
a critical point the term $g_{\matZ}(\xi_{\matZ},\hess{\fcca[\xi_{\matZ}]})$,
which is the numerator of $R(\xi_{\matZ})$, does not depend on the
choice of Riemannian metric. Thus, the optimal values of Eq. (\ref{eq:ccarayleighforMeqB})
satisfying Eq. (\ref{eq:CCAopt1a}) or Eq. (\ref{eq:CCAopt2a}) are
the eigenvalues of the Riemannian Hessian at $\matZ\in\elpCCA$ with
the choice $\matM_{\matZ}\coloneqq\Sigma$. Obviously, classification of the critical
points does not depend on the Riemannian metric. Therefore, we can
classify the critical points using the signs of the eigenvalues of
the Riemannian Hessian at any critical point $\matZ\in\elpCCA$ with the choice $\matM_{\matZ}\coloneqq\Sigma$.

Recall from the proof of Theorem \ref{thm:CCAHessiantheorem} that in the special case $\matM_{\matZ}\coloneqq=\Sigma$, the eigenvalues of Eq. (\ref{eq:CCAHessianoptimal})
are also the eigenvalues of the Riemannian Hessian at a critical point $\matZ\in\elpCCA$. The eigenvalues are:
$p(d_{\y}-p)$ eigenvalues of the form $\mu_{j}\left(\beta_{j}+\beta_{i}\right)$
where $j=1,...,p$ and $i=p+1,...,d_{\y}$, $p(d_{\y}-p)$ eigenvalues
of the form $\mu_{j}\left(\beta_{j}-\beta_{i}\right)$ where $j=1,...,p$
and $i=p+1,...,d_{\y}$, $p(d_{\x}-d_{\y})$ eigenvalues of the form
$\mu_{j}\beta_{j}$ where $j=1,...,p$, $p(p-1)/2$ eigenvalues of
the form $\frac{1}{2}\left(\mu_{j}+\mu_{i}\right)\left(\beta_{j}+\beta_{i}\right)$
where $1\leq j<i\leq p$, and $p(p-1)/2$ eigenvalues of the form
$\frac{1}{2}\left(\mu_{j}-\mu_{i}\right)\left(\beta_{j}-\beta_{i}\right)$
where $1\leq j<i\leq p$.  $\mu_{i}>0$ for $i=1,...,p$, and $\mu_{j}-\mu_{i}>0$
for $j<i$. Also $|\beta_{i}|=\alpha_{i}$ for $i=1,...,d_{\y}$. Thus, the signs of the eigenvalues
of the matrix in Eq.~(\ref{eq:CCAHessianoptimal}) are only determined
by $\beta_{j}$, $\beta_{j}+\beta_{i}$ and $\beta_{j}-\beta_{i}$,
where $\beta_{j}$ are up to their sign the $p$ singular values of
the matrix $\matR\coloneqq\Sigma_{\x\x}^{-\nicehalf}\Sigma_{\x\y}\Sigma_{\y\y}^{-\nicehalf}$ corresponding to left and right singular vectors
which are the columns of $(\tilde{\matU},\tilde{\matV})=(\Sigma_{\x\x}^{\nicehalf}\matU,\Sigma_{\y\y}^{\nicehalf}\matV)$.
In particular, for an optimal $\matZ^{\star}=(\matU^{\star},\matV^{\star})\in\elpCCA$
such that the columns of $(\tilde{\matU},\tilde{\matV})=(\Sigma_{\x\x}^{\nicehalf}\matU^{\star},\Sigma_{\y\y}^{\nicehalf}\matV^{\star})$
are ordered left and right $p$-dominant singular vectors of the matrix
$\matR$ on the same phase, then we have that $\beta_{1}=\sigma_{1}>...>\beta_{p}=\sigma_{p}$.
We conclude that $\beta_{i}$ such that $i=p+1,...,d_{\y}$ satisfies
$|\beta_{i}|=\sigma_{i}$, which leads to $\beta_{j}>0$, $\beta_{j}+\beta_{i}>0$
and $\beta_{j}-\beta_{i}>0$ where $j=1,...,p$ and $i=p+1,...,d_{\y}$
or $1\leq j<i\leq p$. Therefore, in this case all the eigenvalues
of the matrix Eq.~(\ref{eq:CCAHessianoptimal}) are strictly positive,
the matrix Eq.~(\ref{eq:CCAHessianoptimal}) is SPD, and consequently
the eigenvalues of the Riemannian Hessian at $\matZ^{\star}\in\elpCCA$
are all strictly positive. This proves that $\matZ^{\star}$ is a
strict local minimum of $\fcca$ on $\elpCCA$, see \cite[Proposition 6.5.]{boumal2022intromanifolds}.

If we prove that $\matZ^{\star}$ is the only local minimum
(up to the signs of the columns of $\matU^{\star}$ and $\matV^{\star}$),
then $\matZ^{\star}$ is the only asymptotically stable critical point
following Theorem \ref{thm:stabilityCCA}. In order to prove it we
further assume that for all $i=1,...,q$ the values $\sigma_{i}$
are distinct, then we can conclude the following. Suppose $\matZ=(\matU,\matV)$
is any other critical point differs from $\matZ^{\star}$ at the optimal
value, i.e., such that the columns of $(\tilde{\matU},\tilde{\matV})=(\Sigma_{\x\x}^{\nicehalf}\matU,\Sigma_{\y\y}^{\nicehalf}\matV)$
are left and right singular vectors corresponding to some $p$ singular
values of the matrix $\matR$ not necessarily on the same phase, so
that Eq. (\ref{eq:CCAproofform1}) and Eq. (\ref{eq:CCAproofform2})
hold, and there exists at least one $1\leq j\leq p$ for which $\beta_{j}\neq\sigma_{j}$.
We consider the different cases:
\begin{enumerate}
\item Suppose $\beta_{1},\dots,\beta_{p}$ are not ordered in any particular
order (possible for $p\geq3$), then there exists $j$ such that $\beta_{j}$
is larger than some $\beta_{k}$ and smaller than $\beta_{m}$ where
$j<k,m\leq p$, then $\beta_{j}-\beta_{k}>0$ and $\beta_{j}-\beta_{m}<0$.
In this case there are both strictly positive and strictly negative
eigenvalues of the Riemannian Hessian at $\matZ$ for the choice $\matM_{\matZ}\coloneqq\Sigma$,
therefore, $\matZ$ is a saddle point.
\item Suppose $\beta_{1},\dots,\beta_{p}$ are ordered in a descending order.
Since $\matZ$ is not an optimal solution of Problem~(\ref{eq:costriemannianccaBrockett}),
then there exists at least one $1\leq j\leq p$ for which $\beta_{j}\neq\sigma_{j}$.
Thus, on the one hand $\beta_{j}-\beta_{i}>0$ where $1\leq j<i\leq p$,
but on the other hand there exists at least one pair of indexes $j=1,...,p$
and $i=p+1,...,d_{\y}$ such that $\beta_{j}+\beta_{i}<0$ or $\beta_{j}-\beta_{i}<0$,
otherwise, it contradicts the assumption that there exists at least
one $1\leq j\leq p$ for which $\beta_{j}\neq\sigma_{j}$. In this
case there are both strictly positive and strictly negative eigenvalues
of the Riemannian Hessian at $\matZ$ for the choice $\matM_{\matZ}\coloneqq\Sigma$,
therefore, $\matZ$ is a saddle point.
\item Suppose $\beta_{1},\dots,\beta_{p}$ are ordered in an ascending order.
Then, $\beta_{j}-\beta_{i}<0$ where $1\leq j<i\leq p$. Now, we consider
two sub-cases:
\begin{enumerate}
\item There exists at least one $1\leq j\leq p$ for which $\beta_{j}\neq-\sigma_{j}$.
Then, there exists at least one pair of indexes $j=1,...,p$ and $i=p+1,...,d_{\y}$
such that $\beta_{j}+\beta_{i}>0$ or $\beta_{j}-\beta_{i}>0$, otherwise
it contradicts the assumption that there exists at least one $1\leq j\leq p$
for which $\beta_{j}\neq-\sigma_{j}$. In this case there are both
strictly positive and strictly negative eigenvalues of the Riemannian
Hessian at $\matZ$ for the choice $\matM_{\matZ}\coloneqq=\Sigma$, therefore, $\matZ$
is a saddle point.
\item Consider the case $\beta_{1}=-\sigma_{1}<...<\beta_{p}=-\sigma_{p}$.
Then, $\beta_{j}<0$, $\beta_{j}+\beta_{i}<0$ and $\beta_{j}-\beta_{i}<0$
where $j=1,...,p$ and $i=p+1,...,d_{\y}$ or $1\leq j<i\leq p$.
Therefore, in this case all the eigenvalues of the Riemannian Hessian
at $\matZ$ are strictly negative for the choice $\matM_{\matZ}\coloneqq\Sigma$,
therefore, $\matZ$ is a strict local maximizer. Since it is the only
local maximizer up to the signs of the columns of $\matU$ and $\matV$,
and it is also a global maximizer, it is the unique maximizer.
\end{enumerate}
\end{enumerate}
In all cases, $\matZ$ is not a local minimizer. Thus, $\matZ^{\star}$ is the only local minimum
(up to the signs of the columns of $\matU^{\star}$ and $\matV^{\star}$). According to Theorem
\ref{thm:stabilityCCA} all the other critical points are unstable.
\end{proof}

\subsection{Proof of Corollary \ref{cor:CCA}}\label{subsec:proofofcorCCA}

\begin{proof}
The condition number bound follows from Lemma \ref{lem:sketching}
and Theorem \ref{thm:CCAHessiantheorem} with one additional argument.
From Lemma \ref{lem:sketching} we know that with probability of at
least $1-\delta$ we have that all the generalized eigenvalues of
the pencil $(\Sigma_{\x\x},\matM^{(\x\x)})$ are contained in the
interval $[1/2,3/2]$, and the same is true for the pencil $(\Sigma_{\y\y},\matM^{(\y\y)})$.
Recall that all the generalized eigenvalues of the pencil $(\Sigma_{\x\x},\matM^{(\x\x)})$
are also generalized eigenvalues of the pencil $\left(\Sigma,\matM\right)$,
and the same is true for the generalized eigenvalues of $(\Sigma_{\y\y},\matM^{(\y\y)})$.
Indeed, after an appropriate padding with zeros each generalized eigenvector of $(\Sigma_{\x\x},\matM^{(\x\x)})$ or $(\Sigma_{\y\y},\matM^{(\y\y)})$
is also a generalized eigenvector of $\left(\Sigma,\matM\right)$. Thus, since $(\Sigma_{\x\x},\matM^{(\x\x)})$ and $(\Sigma_{\y\y},\matM^{(\y\y)})$
have $d_{\x}$ and $d_{\y}$ generalized eigenvalues and corresponding eigenvectors, they characterize all the generalized eigenvalues and corresponding eigenvectors of $\left(\Sigma,\matM\right)$.
Thus, also the eigenvalues of the pencil $\left(\Sigma,\matM\right)$
are contained in the interval $[1/2,3/2]$, and subsequently $\kappa\left(\Sigma,\matM\right)\leq3.$
Using Theorem \ref{thm:CCAHessiantheorem} we get the require bound
for the condition number.

The costs are evident from Table \ref{tab:costs}, once we observe
that none of the operations require forming $\Sigma_{\x\x},\Sigma_{\y\y}$
or $\Sigma_{\x\y}$, but instead require taking product of these matrices
with vectors. These products can be computed in cost proportional
to the number of non-zeros in $\matX$ and/or $\matY$ by iterated
products. In addition, we use the fact that $\matS\matX$ and $\matS\matY$
can be computed in $O(\nnz{\matX})=O(nd_{\x})$ and $O(\nnz{\matY})=O(nd_{\y})$
operations. The required preprocessing is to factorize $\matM^{(\x\x)}$
and $\matM^{(\y\y)}$, so we can efficiently take products with $\left(\matM^{(\x\x)}\right)^{-1}$
and $\left(\matM^{(\y\y)}\right)^{-1}$, which can be done in $O(sd^{2}+d^{2})$
as explained in Section \ref{sec:rand-precond}. Assuming a bounded
number of line-search steps in each iteration of Riemannian CG, each
iteration requires a bounded number of computations of each of the
following: objective function evaluation costs $O\left(p\left(\nnz{\matX}+\nnz{\matY}\right)\right)=O\left(pnd\right)$,
retraction costs $O\left(p\left(\nnz{\matX}+\nnz{\matY}\right)+dp^{2}\right)$,
vector transport and Riemannian gradient computation take $O\left(p\left(\nnz{\matX}+\nnz{\matY}\right)+dp^{2}+d^{2}p\right)$.
\end{proof}

\section{Omitted Proofs From Section~\ref{sec:Randomized-preconditioning-for-LDA}}
In this section we give omitted proofs from Section~\ref{sec:Randomized-preconditioning-for-LDA}.

\subsection{Proof of Theorem~\ref{thm:criticalpointsofLDA}}\label{subsec:proofcriticalpointsofLDA}
\begin{proof}
Recall that the critical points are defined to be the points where the Riemannian gradient is zero, and as such,
whether a point is a critical point or not does not depend on the
choice of Riemannian metric (see \cite[Eq. 3.31]{AMS09}). Thus, for
the sake of identifying the critical points, we can assume that $\matM_{\matW}\coloneqq\ldaB$
and use a simplified form for the Riemannian gradient: 
\begin{align*}
\grad{\flda} & =  \mat{\Pi}_{\matW}\left((\ldaB)^{-1}\nabla\bar{f}_{{\bf FDA}}(\mat W)\right)\\
 & =  -\mat{\Pi}_{\matW}\left((\ldaB)^{-1}\mat{S_{B}}\mat W\matN\right)\\
 & =  -\left[(\ldaB)^{-1}\mat{S_{B}}\mat W\matN-\matW\sym{\left(\mat W^{\T}\mat{S_{B}}\mat W\matN\right)}\right]\\
 & =  -\left[\left(\matI_{d}-\mat W\mat W^{\T}(\ldaB)\right)(\ldaB)^{-1}\mat{S_{B}}\mat W\matN \right. \\ & \left.   + \ \mat W\skew{\left(\mat W^{\T}\mat{S_{B}}\mat W\matN\right)}\right].
\end{align*}
Let $\alpha_{1},\dots,\alpha_{p}$ be some generalized eigenvalues
of the matrix pencil $(\mat{S_{B}},\ldaB)$, and let $\w_{1},...,\mathbf{w}_{p}$
be the corresponding generalized eigenvectors (making them some $p$ FDA weight vectors). Writing $\w_{1},...,\mathbf{w}_{p}$
as the columns of $\matW$, the following equation holds:
\[
\mat{S_{B}}\mat W=(\ldaB)\mat W\matA
\]
where $\matA=\diag{\alpha_{1},\dots,\alpha_{p}}$. Letting $\matW\in\Stiefellda$,
we have 
\begin{eqnarray*}
\grad{\flda} & = & -\left[(\ldaB)^{-1}\mat{S_{B}}\mat W\matN-\matW\sym{\left(\mat W^{\T}\mat{S_{B}}\mat W\matN\right)}\right]\\
 & = & -\left[\mat W\matA\matN-\matW\sym{\left(\matA\matN\right)}\right] = \mat 0_{d\times p}.
\end{eqnarray*}

To show the other side, note that if the Riemannian gradient nullifies,
then 
\begin{eqnarray*}
\left[\left(\matI_{d}-\mat W\mat W^{\T}(\ldaB)\right)(\ldaB)^{-1}\mat{S_{B}}\mat W\matN+\mat W\skew{\left(\mat W^{\T}\mat{S_{B}}\mat W\matN\right)}\right] & = & \mat 0_{d\times p}\ .
\end{eqnarray*}
By using similar reasoning as in \cite[Subsection 4.8.2]{AMS09},
$$\left(\matI_{d}-\mat W\mat W^{\T}(\ldaB)\right)(\ldaB)^{-1}\mat{S_{B}}\mat W\matN,$$
belongs to the orthogonal compliment of the column space of $\mat W$
(with respect to the matrix $(\ldaB)$), and $\mat W\skew{\left(\mat W^{\T}\mat{S_{B}}\mat W\matN\right)}$
belongs to the column space of $\mat W$. Thus, we get that the gradient
vanishes if and only if the following two formulas hold:
\begin{equation}
\left(\matI_{d}-\mat W\mat W^{\T}(\ldaB)\right)(\ldaB)^{-1}\mat{S_{B}}\mat W\matN=\mat 0_{d\times p}\ ,\label{eq:critical1-2}
\end{equation}
and 
\begin{equation}
\mat W\skew{\left(\mat W^{\T}\mat{S_{B}}\mat W\matN\right)}=\mat 0_{d\times p}\ .\label{eq:critical2-2}
\end{equation}
From Eq.~(\ref{eq:critical1-2}) we get
\[
(\ldaB)^{-1}\mat{S_{B}}\mat W=\mat W\left(\mat W^{\T}\mat{S_{B}}\mat W\right),
\]
since $\matN$ is an invertible matrix. Also, since $\mat W\in\Stiefellda$
it is a full (column) rank matrix then Eq.~(\ref{eq:critical2-2})
vanishes if and only if 
\[
\skew{\left(\mat W^{\T}\mat{S_{B}}\mat W\matN\right)}=\mat 0_{p}\ ,
\]
which leads to $\left(\mat W^{\T}\mat{S_{B}}\mat W\right)\matN=\matN\left(\mat W^{\T}\mat{S_{B}}\mat W\right)$.
This implies that $\left(\mat W^{\T}\mat{S_{B}}\mat W\right)$ is
diagonal because any rectangle matrix that commutes with a diagonal
matrix with distinct entries is diagonal. Finally we get 
\[
(\ldaB)^{-1}\mat{S_{B}}\mat W=\mat W\mat D\ ,
\]
where $\mat D=\mat W^{\T}\mat{S_{B}}\mat W$ is a diagonal matrix.
This implies that the columns of $\mat W$ correspond to some $p$
generalized eigenvectors the matrix pencil $(\mat{S_{B}},\ldaB)$, thus, making the columns some $p$ FDA weight vectors.

Finally, to identify the optimal solutions, note that at the critical points the objective function is a sum of the generalized eigenvalues corresponding to the columns of $\matW$ multiplied by a diagonal element of $\matN$. Thus, the optimal solutions that minimize $\flda$ on $\Stiefellda$
are $\matW\in\Stiefellda$ such that the columns correspond to the $p$ leading FDA weight vectors. Otherwise, we can increase the value of the objective function by
replacing a weight vector with another that corresponds to a smaller generalized eigenvalue. Moreover, if we assume that $\rho_{1}>\rho_{2}>...>\rho_{p+1}\geq0$,
then for the aforementioned $\matW$, the
columns belong each to a one dimensional
generalized eigenspace, i.e., unique solution up the the signs of the columns. In the case where some $\rho_{i}=\rho_{j}$
for $1\leq i,j\leq p$, permutations of the columns of $\matW$ associated with $\rho_{i}$ keep the solution optimal making it non-unique.
\end{proof}

\subsection{Proof of Theorem \ref{thm:stabilityLDA}}\label{subsec:proofstabilityLDA}

\begin{proof}
To prove the asymptotic stability of a $\matW\in\Stiefellda$ such that the columns are the $p$ leading FDA weight vectors
we use \cite[Proposition 4.4.2]{AMS09}. Recall from Theorem \ref{thm:criticalpointsofLDA} that $\matW$ that solve Problem~(\ref{eq:lda_riemannian})
are unique up the the signs of the columns of $\matW$, making these
points isolated global (and consequently local) minimizers of $\justflda$
on $\Stiefellda$. According to \cite[Proposition 4.4.2]{AMS09},
such points $\matW$ are asymptotically stable.

Suppose $\matW$ is a critical point of $\flda$ on
$\Stiefellda$ which is not a local minimum. Then, there exists compact
neighborhoods with either no other critical points, if there are no
multiplicities of the generalized eigenspaces, or where all other
critical point achieve the same value for the cost function, if there
are multiplicities of the generalized eigenspaces. Thus, according
to \cite[Proposition 4.4.1]{AMS09}, such $\matW$ are unstable.
\end{proof}

\subsection{Proof of Theorem~\ref{thm:LDAHessiantheorem}}\label{subsec:proofLDAHessiantheorem}
\begin{proof} The proof is similar to the proof of Theorem \ref{thm:CCAHessiantheorem}. In order to
bound the condition number of the Riemannian Hessian at $\matW^{\star}$,
we need to bound its maximal and minimal eigenvalues. Thus,
to prove the theorem we analyze the eigenvalues of the Riemannian
Hessian at some critical point $\matW\in\Stiefellda$ (in particular
at $\matW^{\star}$) using the Courant-Fischer Theorem
for the compact self-adjoint linear operator $$\hess{\flda[\cdot]}:T_{\matW}\Stiefellda\to T_{\matW}\Stiefellda,$$
over the finite-dimensional vector space $T_{\matW}\Stiefellda$, see
\cite[Chapter 1, Section 6.10]{kato2013perturbation}:
\begin{eqnarray}
\lambda_{k}(\hess{\flda}) & = & \min_{U,\dim(U)=k-1}\max_{\mat 0\neq\xi_{\matW}\in U^{\perp}}R(\xi_{\matW}),\label{eq:LDAopt1a-1}
\end{eqnarray}
\begin{eqnarray}
\lambda_{k}(\hess{\flda}) & = & \max_{U,\dim(U)=k}\min_{\mat 0\neq\xi_{\matW}\in U}R(\xi_{\matW}),\label{eq:LDAopt2a-1}
\end{eqnarray}
where 
\[
R(\xi_{\matW})\coloneqq\frac{g_{\matW}(\xi_{\matW},\hess{\flda[\xi_{\matW}]})}{g_{\matW}(\xi_{\matW},\xi_{\matW})},
\]
$\lambda_{k}(\hess{\flda})=\rho_{k}$ is the $k$th largest eigenvalue
(descending order) of $\hess{\flda}$, and $U$ is a linear subspace
of $T_{\matW}\Stiefellda$. In particular, the maximal and minimal
eigenvalues are given by the formulas:
\begin{equation}
\lambda_{\max}(\hess{\flda})=\max_{\mat 0\neq\xi_{\matW}\in T_{\matW}\Stiefellda}R(\xi_{\matW}),\label{eq:LDAopt1-1}
\end{equation}
 
\begin{equation}
\lambda_{\min}(\hess{\flda})=\min_{\mat 0\neq\xi_{\matW}\in T_{\matW}\Stiefellda}R(\xi_{\matW}).\label{eq:LDAopt2-1}
\end{equation}

We begin by simplifying the quotient $R(\xi_{\matW})$. Recall that
any critical point of $\justflda(\cdot)$ is a matrix $\matW\in\Stiefellda$
such that the columns are some $p$ generalized eigenvectors of the
matrix pencil $(\mat{S_{B}},\ldaB)$ (see Theorem \ref{thm:criticalpointsofLDA}).
Let $\alpha_{1},\dots,\alpha_{p}$ be some generalized eigenvalues of the
matrix pencil $(\mat{S_{B}},\ldaB)$, and let $\w_{1},...,\mathbf{w}_{p}$
be the corresponding generalized eigenvectors. Writing $\w_{1},...,\mathbf{w}_{p}$
as the columns of $\matW$, the following equation holds:
\[
\mat{S_{B}}\mat W=(\ldaB)\mat W\matA
\]
where $\matA=\diag{\alpha_{1},\dots,\alpha_{p}}$. Letting $\matW\in\Stiefellda$,
plugging in the ambient coordinates formula for the Riemannian Hessian
(Eq. (\ref{eq:LDAhess})), the Riemannian gradient nullifies (see
Theorem \ref{thm:criticalpointsofLDA}) and we have{\footnotesize{}
\begin{eqnarray}
\text{Hess}\flda[\xi_{\mat W}] & = & \mat{\Pi}_{\matW}\left(\matM^{-1}_{\matW}\left[-\mat{S_{B}}\xi_{\mat W}\matN+(\ldaB)\xi_{\mat W}\left(\mat W^{\T}\mat{S_{B}}\mat W\matN+\grad{f(\mat W)}\right)\right]\right)\nonumber \\
 & = & \mat{\Pi}_{\matW}\left(\matM^{-1}_{\matW}\left[-\mat{S_{B}}\xi_{\mat W}\matN+(\ldaB)\xi_{\mat W}\left(\mat W^{\T}(\ldaB)\mat W\matA\matN\right)\right]\right)\nonumber \\
 & = & \mat{\Pi}_{\matW}\left(\matM^{-1}_{\matW}\left[-\mat{S_{B}}\xi_{\mat W}+(\ldaB)\xi_{\mat W}\matA\right]\matN\right).\label{eq:hessatcritical-1-1}
\end{eqnarray}
}{\footnotesize\par}

Plugging in the formula for the Riemannian Hessian at a critical point
(Eq. (\ref{eq:hessatcritical-1-1})), the quotient $R(\xi_{\matW})$
reduces to 
\[
R(\xi_{\matW})=\frac{\Trace{\xi_{\matW}^{\T}\matM_{\matW}\mat{\Pi}_{\matW}\left(\matM^{-1}_{\matW}\left[-\mat{S_{B}}\xi_{\mat W}+(\ldaB)\xi_{\mat W}\matA\right]\matN\right)}}{\Trace{\xi_{\matW}^{\T}\matM_{\matW}\xi_{\matW}}}.
\]
Now using the fact that the projection to the tangent space is self-adjoint
with respect to the Riemannian metric,
and that for any $\xi_{\matW}\in T_{\matW}\Stiefellda$ we have $\mat{\Pi}_{\matW}\left(\xi_{\matW}\right)=\xi_{\matW}$,
we further see that 
\begin{multline*}
\frac{\Trace{\xi_{\matW}^{\T}\matM_{\matW}\mat{\Pi}_{\matW}\left(\matM^{-1}_{\matW}\left[-\mat{S_{B}}\xi_{\mat W}+(\ldaB)\xi_{\mat W}\matA\right]\matN\right)}}{\Trace{\xi_{\matW}^{\T}\matM_{\matW}\xi_{\matW}}}=\\ \frac{\Trace{\xi_{\matW}^{\T}\left[-\mat{S_{B}}\xi_{\mat W}+(\ldaB)\xi_{\mat W}\matA\right]\matN}}{\Trace{\xi_{\matW}^{\T}\matM_{\matW}\xi_{\matW}}}
\end{multline*}
Obviously, we can also write 
\begin{multline}
\frac{\Trace{\xi_{\matW}^{\T}\left[-\mat{S_{B}}\xi_{\mat W}+(\ldaB)\xi_{\mat W}\matA\right]\matN}}{\Trace{\xi_{\matW}^{\T}\matM_{\matW}\xi_{\matW}}} \\ =  \frac{\Trace{\xi_{\matW}^{\T}\left[-\mat{S_{B}}\xi_{\mat W}+(\ldaB)\xi_{\mat W}\matA\right]\matN}}{\Trace{\xi_{\matW}^{\T}(\ldaB)\xi_{\matW}}}\cdot\frac{\Trace{\xi_{\matW}^{\T}(\ldaB)\xi_{\matW}}}{\Trace{\xi_{\matW}^{\T}\matM_{\matW}\xi_{\matW}}}.
\label{eq:LDAproofRayleigh}
\end{multline}
Using Eq. \eqref{eq:LDAproofRayleigh}, a simplified
form of the quotient $R(\xi_{\matW})$, we can estimate upper and lower bounds on $R(\xi_{\matW})$ where $\mat 0\neq\xi_{\matW}\in T_{\matW}\Stiefellda$
in order to bound the condition number of the Riemannian Hessian at
$\matW^{\star}$. Since for $\xi_{\matW}\neq\mat 0$
the term $$\Trace{\xi_{\matW}^{\T}(\ldaB)\xi_{\matW}}/\Trace{\xi_{\matW}^{\T}\matM_{\matW}\xi_{\matW}}>0,$$
the upper and lower bounds on
\begin{equation}
\frac{\Trace{\xi_{\matW}^{\T}\left[-\mat{S_{B}}\xi_{\mat W}+(\ldaB)\xi_{\mat W}\matA\right]\matN}}{\Trace{\xi_{\matW}^{\T}(\ldaB)\xi_{\matW}}},\label{eq:ldarayleighforMeqB-1}
\end{equation}
together with the upper and lower bounds of $\Trace{\xi_{\matW}^{\T}(\ldaB)\xi_{\matW}}/\Trace{\xi_{\matW}^{\T}\matM_{\matW}\xi_{\matW}}$
bound the condition number of the Riemannian Hessian at $\matW^{\star}$. 

We begin by estimating the term $\Trace{\xi_{\matW}^{\T}(\ldaB)\xi_{\matW}}/\Trace{\xi_{\matW}^{\T}\matM_{\matW}\xi_{\matW}}$. We use the vectorization operator and the Kronecker Product to
rewrite it in the following form
\begin{equation}
\frac{\Trace{\xi_{\matW}^{\T}(\ldaB)\xi_{\matW}}}{\Trace{\xi_{\matW}^{\T}\matM_{\matW}\xi_{\matW}}}=\frac{\vectorization{\xi_{\matW}}^{\T}\left(\matI_{p}\otimes(\ldaB)\right)\vectorization{\xi_{\matW}}}{\vectorization{\xi_{\matW}}^{\T}\left(\matI_{p}\otimes\matM_{\matW}\right)\vectorization{\xi_{\matW}}}.\label{eq:LDAhessiansecondRayleighquotient-1}
\end{equation}
The righthand side of Eq.~(\ref{eq:LDAhessiansecondRayleighquotient-1})
is the generalized Rayleigh quotient for the matrix pencil $$\left(\matI_{p}\otimes(\ldaB),\matI_{p}\otimes\matM_{\matW}\right).$$
Note that $\matI_{p}\otimes\left(\mat{S_{w}}+\lambda\matI_{d}\right)$
and $\matI_{p}\otimes\matM_{\matW}$ are both SPD matrices, thus the eigenvalues
of the matrix pencil $$\left(\matI_{p}\otimes(\ldaB),\matI_{p}\otimes\matM_{\matW}\right)$$
are equivalent to the eigenvalues of the matrix $\left(\matI_{p}\otimes\matM_{\matW}\right)^{-1}\left(\matI_{p}\otimes(\ldaB)\right)=\matI_{p}\otimes\matM^{-1}_{\matW}(\ldaB)$
and all positive. According to \cite[Section 2]{minka2000old} the
eigenvalues $\matI_{p}\otimes\matM^{-1}_{\matW}(\ldaB)$ are $p$ copies
of each of the eigenvalues of $\matM^{-1}_{\matW}(\ldaB)$. Thus, the maximal
and minimal eigenvalues of the matrix pencil $\left(\matI_{p}\otimes(\ldaB),\matI_{p}\otimes\matM_{\matW}\right)$
denoted by $\tilde{\lambda}_{\text{max}}$ and $\tilde{\lambda}_{\min}$
are equivalent to the maximal and minimal eigenvalues of the matrix
pencil $\left(\ldaB,\matM_{\matW}\right)$, and so is the corresponding condition
number 
\[
\kappa\left(\matI_{p}\otimes\left(\mat{S_{w}}+\lambda\matI_{d}\right),\matI_{p}\otimes\matM_{\matW}\right)=\frac{\tilde{\lambda}_{\text{max}}}{\tilde{\lambda}_{\min}}=\kappa\left(\left(\mat{S_{w}}+\lambda\matI_{d}\right),\matM_{\matW}\right).
\]
Therefore, using the Courant-Fischer Theorem 
\begin{eqnarray*}
\tilde{\lambda}_{\text{max}} & \eqqcolon & \lambda_{\max}(\matI_{p}\otimes(\ldaB),\matI_{p}\otimes\matM_{\matW})\\
 & = & \max_{\mat 0\neq\xi_{\matW}\in\R^{d\times p}}\frac{\vectorization{\xi_{\matW}}^{\T}\left(\matI_{p}\otimes(\ldaB)\right)\vectorization{\xi_{\matW}}}{\vectorization{\xi_{\matW}}^{\T}\left(\matI_{p}\otimes\matM_{\matW}\right)\vectorization{\xi_{\matW}}}\\
 & \geq & \max_{\mat 0\neq\xi_{\matW}\in T_{\matW}\Stiefellda}\frac{\vectorization{\xi_{\matW}}^{\T}\left(\matI_{p}\otimes(\ldaB)\right)\vectorization{\xi_{\matW}}}{\vectorization{\xi_{\matW}}^{\T}\left(\matI_{p}\otimes\matM_{\matW}\right)\vectorization{\xi_{\matW}}}
\end{eqnarray*}
and
\begin{eqnarray*}
\tilde{\lambda}_{\min} & \eqqcolon & \lambda_{\min}(\matI_{p}\otimes(\ldaB),\matI_{p}\otimes\matM_{\matW})\\
 & = & \min_{\mat 0\neq\xi_{\matW}\in\R^{d\times p}}\frac{\vectorization{\xi_{\matW}}^{\T}\left(\matI_{p}\otimes(\ldaB)\right)\vectorization{\xi_{\matW}}}{\vectorization{\xi_{\matW}}^{\T}\left(\matI_{p}\otimes\matM_{\matW}\right)\vectorization{\xi_{\matW}}}\\
 & \leq & \min_{\mat 0\neq\xi_{\matW}\in T_{\matW}\Stiefellda}\frac{\vectorization{\xi_{\matW}}^{\T}\left(\matI_{p}\otimes(\ldaB)\right)\vectorization{\xi_{\matW}}}{\vectorization{\xi_{\matW}}^{\T}\left(\matI_{p}\otimes\matM_{\matW}\right)\vectorization{\xi_{\matW}}}
\end{eqnarray*}

Next, we consider Eq. (\ref{eq:ldarayleighforMeqB-1}). Recall that
we can rewrite any tangent vector, $\xi_{\matW}\in T_{\matW}\Stiefellda$,
in the following form: 
\[
\xi_{\matW}=\matW\mat{\Omega}_{\matW}+\matW_{(\ldaB)\perp}\mat K_{\matW},
\]
where $\matW_{(\ldaB)\perp}$ is $(\ldaB)$-orthogonal to $\matW$
so that the union of the columns of $\matW$ and $\matW_{\left(\mat{S_{w}}+\lambda\matI_{d}\right)\perp}$
is a basis to $\R^{d}$, and $\mat{\Omega}_{\matW}=-\mat{\Omega}_{\matW}^{\T}\in\R^{p\times p}$.
Note that we can always make the choice of the columns of $\matW_{(\ldaB)\perp}$
to be some $d-p$ generalized eigenvalues of the matrix pencil $(\mat{S_{B}},\ldaB)$.
With this choice we have
\[
\mat{S_{B}}\matW_{(\ldaB)\perp}=(\ldaB)\matW_{(\ldaB)\perp}\tilde{\matA},\ \matW_{(\ldaB)\perp}^{\T}\mat{S_{B}}\matW_{(\ldaB)\perp}=\tilde{\matA},
\]
where $\tilde{\matA}\in\R^{\left(d-p\right)\times\left(d-p\right)}$
is a diagonal matrix with the corresponding generalized eigenvalues
on the diagonal $\alpha_{p+1},...,\alpha_{d}$. 

Now, we have 
\[
\xi_{\matW}^{\T}(\ldaB)\xi_{\matW}=\mat{\Omega}_{\matW}^{\T}\mat{\Omega}_{\matW}+\mat K_{\matW}^{\T}\mat K_{\matW}
\]
and
\[
\xi_{\matW}^{\T}\mat{S_{B}}\xi_{\matW}=\mat{\Omega}_{\matW}^{\T}\matA\mat{\Omega}_{\matW}+\mat K_{\matW}^{\T}\tilde{\matA}\mat K_{\matW}\ .
\]

Substitute these computations into Eq. (\ref{eq:ldarayleighforMeqB-1})
\begin{scriptsize}
\begin{eqnarray*}
\frac{\Trace{\xi_{\matW}^{\T}\left[-\mat{S_{B}}\xi_{\mat W}+(\ldaB)\xi_{\mat W}\matA\right]\matN}}{\Trace{\xi_{\matW}^{\T}(\ldaB)\xi_{\matW}}} & = & \frac{\Trace{-\left(\mat{\Omega}_{\matW}^{\T}\matA\mat{\Omega}_{\matW}+\mat K_{\matW}^{\T}\tilde{\matA}\mat K_{\matW}\right)\matN+\left(\mat{\Omega}_{\matW}^{\T}\mat{\Omega}_{\matW}+\mat K_{\matW}^{\T}\mat K_{\matW}\right)\matA\matN}}{\Trace{\mat{\Omega}_{\matW}^{\T}\mat{\Omega}_{\matW}+\mat K_{\matW}^{\T}\mat K_{\matW}}}\\
 & = & \frac{\vectorization{\mat{\Omega}_{\matW}}^{\T}\matD_{1}\vectorization{\mat{\Omega}_{\matW}}+\vectorization{\mat K_{\matW}}^{\T}\matD_{2}\vectorization{\mat K_{\matW}}}{\vectorization{\mat{\Omega}_{\matW}}^{\T}\vectorization{\mat{\Omega}_{\matW}}+\vectorization{\mat K_{\matW}}^{\T}\vectorization{\mat K_{\matW}}},
\end{eqnarray*}
\end{scriptsize}
where 
\[
\matD_{1}\coloneqq\matA\matN\otimes\matI_{p}-\matN\otimes\matA,\:\mat D_{2}\coloneqq\matA\matN\otimes\matI_{d-p}-\matN\otimes\tilde{\matA}.
\]

RRecall that $\mat{\Omega}_{\matW}=-\mat{\Omega}_{\matW}^{\T}\in\R^{p\times p}$,
and both are real matrices (so the elements of the main diagonals
are 0), we have 
\begin{equation}
\vectorization{\mat{\Omega}_{\matW}}^{\T}\vectorization{\mat{\Omega}_{\matW}}=2\sum_{p\geq i>j}\left(\mat{\Omega}_{\matW}\right)_{ij}^{2}\label{eq:omegawomegaw-1}
\end{equation}
and 
\begin{equation}
\vectorization{\mat{\Omega}_{\matW}}^{\T}\left(\matA\matN\otimes\matI_{p}-\matN\otimes\matA\right)\vectorization{\mat{\Omega}_{\matW}}=\sum_{1\leq j<i\leq p}\left(\mat{\Omega}_{\matW}\right)_{ij}^{2}\left(\mu_{j}-\mu_{i}\right)\left(\alpha_{j}-\alpha_{i}\right).\label{eq:omegawAomegaw-1}
\end{equation}
Thus, only the $p(p-1)/2$ entries below the diagonal of $\mat{\Omega}_{\matW}$
determine the values of Eq. (\ref{eq:omegawomegaw-1}) and Eq. (\ref{eq:omegawAomegaw-1}).
Let us now denote by $\m_{\matW}$ the column stack of $\mat{\Omega}_{\matW}$,
but only with the the subdiagonal entries of $\mat{\Omega}_{\matW}$
(i.e., $\vectorization{\mat{\Omega}_{\matW}}$ purged of the superdiagonal
elements).  We then have 
\[
\vectorization{\mat{\Omega}_{\matW}}^{\T}\vectorization{\mat{\Omega}_{\matW}}=2\m_{\matW}^{\T}\m_{\matW}
\]
and 
\[
\vectorization{\mat{\Omega}_{\matW}}^{\T}\left(\matA\matN\otimes\matI_{p}-\matN\otimes\matA\right)\vectorization{\mat{\Omega}_{\matW}}=\m_{\matW}^{\T}\mat{\Psi}\m_{\matW}
\]
where $\mat{\Psi}\in\R^{p(p-1)/2\times p(p-1)/2}$ is a block diagonal
matrix, where the blocks are of descending order from $p-1,p-2,...,1$,
and each block is a diagonal matrix as well. The $j$th ($1\leq j\leq p-1$)
block, denoted by $\mat{\Psi}_{j}$, is of the order $p-j$ and has
the following form
\[
\mat{\Psi}_{j}=\diag{\left(\mu_{j}-\mu_{j+1}\right)\left(\alpha_{j}-\alpha_{j+1}\right),\left(\mu_{j}-\mu_{j+2}\right)\left(\alpha_{j}-\alpha_{j+2}\right),...,\left(\mu_{j}-\mu_{p}\right)\left(\alpha_{j}-\alpha_{p}\right)}
\]
Now, we make the following change of variables: $\d_{\matW}\coloneqq\sqrt{2}\m_{\matW}$,
$\k_{\matW}\coloneqq\vectorization{\mat K_{\matW}}$. Finally, Eq.
(\ref{eq:ldarayleighforMeqB-1}) is rewritten in the following way{\footnotesize{}
\begin{multline}
\frac{\Trace{\xi_{\matW}^{\T}\left[-\mat{S_{B}}\xi_{\mat W}+(\ldaB)\xi_{\mat W}\matA\right]\matN}}{\Trace{\xi_{\matW}^{\T}(\ldaB)\xi_{\matW}}} \\ =  \frac{\left[\begin{array}{cc}
\d_{\matW}^{\T} & \k_{\matW}^{\T}\end{array}\right]\blockdiag{\frac{1}{2}\mat{\Psi},\matA\matN\otimes\matI_{d-p}-\matN\otimes\tilde{\matA}}\left[\begin{array}{c}
\d_{\matW}\\
\k_{\matW}
\end{array}\right]}{\left[\begin{array}{cc}
\d_{\matW}^{\T} & \k_{\matW}^{\T}\end{array}\right]\left[\begin{array}{c}
\d_{\matW}\\
\k_{\matW}
\end{array}\right]}.\label{eq:RayleighLDA-1}
\end{multline}
}{\footnotesize\par}

Note that the mapping $\varphi(\cdot):T_{\matW}\Stiefellda\to\R^{pd-p(p+1)/2}$
defined by 
\[
\varphi(\xi_{\matZ})\coloneqq\left[\begin{array}{c}
\d_{\matW}\\
\k_{\matW}
\end{array}\right],
\]
is a coordinate chart of the elements of $T_{\matW}\Stiefellda$,
since $\varphi(\cdot)$ is a bijection (one-to-one correspondence)
of the elements of $T_{\matW}\Stiefellda$ onto $\R^{pd-p(p+1)/2}$.
Indeed, $\k_{\matW}$ is a column stack of $\mat K_{\matW}$, thus
we can retract the matrices $\mat K_{\matW}$. Similarly $\d_{\matW}$
is proportional to $\m_{\matW}$ which is a column stack of $\mat{\Omega}_{\matW}$,
but only with the the subdiagonal entries of $\mat{\Omega}_{\matW}$.
Since $\mat{\Omega}_{\matW}$ is a skew-symmetric matrix, we can retract
$\mat{\Omega}_{\matW}$. With the matrices $\mat K_{\matW}$ and $\mat{\Omega}_{\matW}$
at hand, we can fully retract $\xi_{\matW}$.

The eigenvalues and corresponding eigenvectors of any linear operator
over a finite dimensional vector space do not depend on the choice
of coordinate chart and basis, thus the eigenvalues and eigenvectors
of $\hess{\flda[\cdot]}:T_{\matW}\Stiefellda\to T_{\matW}\Stiefellda$
which are computed using the Courant Fischer Theorem for compact self-adjoint
linear operators (Eq. (\ref{eq:LDAopt1a-1}) and Eq. (\ref{eq:LDAopt2a-1})),
can be also computed by the Courant Fischer Theorem for symmetric
matrices \cite[Theorem 4.2.6]{horn2012matrix} after applying $\varphi(\cdot)$.
In particular, Eq. (\ref{eq:ldarayleighforMeqB-1}) determines the
signs of the eigenvalues of the Riemannian Hessian at any $\matW\in\Stiefellda$
(in the special case $\matM_{\matW}\coloneqq\ldaB$, the eigenvalues of Eq. (\ref{eq:LDAHessianoptimal-1})
are the eigenvalues of the Riemannian Hessian at $\matW\in\Stiefellda$),
and the bounds of Eq. (\ref{eq:ldarayleighforMeqB-1}) together with
the bounds of the term $\Trace{\xi_{\matW}^{\T}(\ldaB)\xi_{\matW}}/\Trace{\xi_{\matW}^{\T}\matM_{\matW}\xi_{\matW}}$
bound the condition number of the Riemannian Hessian at $\matW^{\star}\in\Stiefellda$. 

To that end, we perform the following computation. The righthand side
of Eq. (\ref{eq:RayleighLDA-1}) is a Rayleigh quotient, so according
to the Courant-Fischer Theorem for symmetric matrices the eigenvalues
(in particular, the maximal and the minimal) of the $pd-p(p+1)/2\times pd-p(p+1)/2$
symmetric matrix
\begin{equation}
\blockdiag{\frac{1}{2}\mat{\Psi},\matA\matN\otimes\matI_{d-p}-\matN\otimes\tilde{\matA}},\label{eq:LDAHessianoptimal-1}
\end{equation}
are determined by optimal (in particular the maximal and the minimal)
values of Eq. (\ref{eq:RayleighLDA-1}). The set of eigenvalues of
Eq. (\ref{eq:LDAHessianoptimal-1}) is equal to the union of the set
of eigenvalues of $\frac{1}{2}\mat{\Psi}$ and $\matA\matN\otimes\matI_{d-p}-\matN\otimes\tilde{\matA}$.
Since Eq. (\ref{eq:LDAHessianoptimal-1}) is a diagonal matrix, its
eigenvalues are simply the diagonal elements:

\[
\frac{1}{2}\left(\mu_{j}-\mu_{i}\right)\left(\alpha_{j}-\alpha_{i}\right),\ 1\leq j<i\leq p,
\]
and
\[
\mu_{j}\left(\alpha_{j}-\alpha_{i}\right),\ j=1,...,p,\ i=p+1,...,d.
\]

Finally, we bound the condition number of the Riemannian Hessian at an optimal point,
$\matW^{\star}\in\Stiefellda$. In such case, $\alpha_{1}=\rho_{1}>...>\alpha_{p}=\rho_{p}$.
Without loss of generality, we can always choose $\matW_{(\ldaB)\perp}$
such that the corresponding eigenvalues to its columns are in a descending
order, thus, $\alpha_{p+1}=\rho_{p+1}\geq...\geq\alpha_{d}=\rho_{d}\geq0$.
Then, we have that Eq. (\ref{eq:ldarayleighforMeqB-1}) is bounded
by the minimal and maximal eigenvalues of (\ref{eq:LDAHessianoptimal-1}).
Thus,
\[
0<\max_{\mat 0\neq\xi_{\matW^{\star}}\in T_{\matW^{\star}}\Stiefellda}\frac{\Trace{\xi_{\matW^{\star}}^{\T}\left[-\mat{S_{B}}\xi_{\matW^{\star}}+(\ldaB)\xi_{\matW^{\star}}\matA\right]\matN}}{\Trace{\xi_{\matW^{\star}}^{\T}(\ldaB)\xi_{\matW^{\star}}}}\leq\mu_{1}\left(\rho_{1}-\rho_{d}\right),
\]
and 
\begin{eqnarray*}
\min_{\mat 0\neq\xi_{\matW^{\star}}\in T_{\matW^{\star}}\Stiefellda}\frac{\Trace{\xi_{\matW^{\star}}^{\T}\left[-\mat{S_{B}}\xi_{\matW^{\star}}+(\ldaB)\xi_{\matW^{\star}}\matA\right]\matN}}{\Trace{\xi_{\matW^{\star}}^{\T}(\ldaB)\xi_{\matW^{\star}}}} & \geq\\
\min\left\{ \mu_{p}\left(\rho_{p}-\rho_{p+1}\right),\min_{1\leq j<p}\frac{1}{2}\left(\mu_{j}-\mu_{j+1}\right)\left(\rho_{j}-\rho_{j+1}\right)\right\}  & > & 0.
\end{eqnarray*}
We use Eq.~(\ref{eq:LDAopt1-1}) and Eq.~(\ref{eq:LDAopt2-1}) to
bound the condition number:
\begin{eqnarray*}
\lambda_{\max}(\hess{\justflda(\matW^{\star})}) & = & \max_{\mat 0\neq\xi_{\matW^{\star}}\in T_{\matW^{\star}}\Stiefellda}\frac{g_{\matW^{\star}}(\xi_{\matW^{\star}},\hess{\justflda(\matW^{\star})[\xi_{\matW^{\star}}]})}{g_{\matW^{\star}}(\xi_{\matW^{\star}},\xi_{\matW^{\star}})}\\
 & \leq & \mu_{1}\left(\rho_{1}-\rho_{d}\right)\tilde{\lambda}_{\text{max}}\ ,
\end{eqnarray*}
and 
\begin{eqnarray*}
\lambda_{\min}(\hess{\justflda(\matW^{\star})}) & = & \min_{\mat 0\neq\xi_{\matW^{\star}}\in T_{\matW^{\star}}\Stiefellda}\frac{g_{\matW^{\star}}(\xi_{\matW^{\star}},\hess{\justflda(\matW^{\star})[\xi_{\matW^{\star}}]})}{g_{\matW^{\star}}(\xi_{\matW^{\star}},\xi_{\matW^{\star}})}\\
 & \geq & \min\left\{ \mu_{p}\left(\rho_{p}-\rho_{p+1}\right),\min_{1\leq j<p}\frac{1}{2}\left(\mu_{j}-\mu_{j+1}\right)\left(\rho_{j}-\rho_{j+1}\right)\right\} \tilde{\lambda}_{\min}\ .
\end{eqnarray*}
Finally,
\[
\kappa(\hess{\justflda(\matW^{\star})})=\frac{\lambda_{\max}(\hess{\justflda(\matW^{\star})})}{\lambda_{\min}(\hess{\justflda(\matW^{\star})})}\leq\kappa_{{\bf FDA}}^{\star}\cdot\kappa\left(\ldaB,\matM_{\matW}\right),
\]
where
\[
\kappa_{{\bf FDA}}^{\star}=\frac{\mu_{1}\left(\rho_{1}-\rho_{d}\right)}{\min\left\{ \mu_{p}\left(\rho_{p}-\rho_{p+1}\right),\min_{1\leq j<p}\frac{1}{2}\left(\mu_{j}-\mu_{j+1}\right)\left(\rho_{j}-\rho_{j+1}\right)\right\} }.
\]
In the special case $\matM_{\matW^{\star}}=\ldaB$, the bound on the condition number
of the Riemannian Hessian at $\matW^{\star}\in\Stiefellda$ is reduced
to 
\[
\kappa(\hess{\justflda(\matW^{\star})})=\kappa_{{\bf FDA}}^{\star}.
\]
\end{proof}

\subsection{Proof of Theorem \ref{thm:morestabilityLDA}}\label{subsec:proofmorestabilityLDA}

\begin{proof}
To show that $\matW^{\star}$ is a strict local minimum of $\flda$ we show that
the eigenvalues of the Riemannian Hessian at $\matW^{\star}$ are
strictly positive (see \cite[Proposition 6.5.]{boumal2022intromanifolds}). Moreover, under the assumption
that for all $i=1,...,d$ the values $\rho_{i}$ are distinct, we
prove that $\matW^{\star}$ is the only local minimum of $\flda$ (and it is also strict)
and all other critical points are either saddle points or strict local maximizers.

Eq. \eqref{eq:LDAproofRayleigh} helps to determine the signs of the eigenvalues of the Riemannian Hessian at
any critical point $\matW\in\Stiefellda$ and in particular at $\matW^{\star}$:
the matrices $\ldaB$ and $\matM_{\matW}$ are both SPD matrices, therefore
for $\xi_{\matW}\neq\mat 0$ the term $$\Trace{\xi_{\matW}^{\T}(\ldaB)\xi_{\matW}}/\Trace{\xi_{\matW}^{\T}\matM_{\matW}\xi_{\matW}}>0,$$
thus only Eq. (\ref{eq:ldarayleighforMeqB-1}),
where $\mat 0\neq\xi_{\matW}\in T_{\matW}\Stiefellda$ determines
the signs. In addition, at a critical point $\matW\in\Stiefellda$
Eq. (\ref{eq:ldarayleighforMeqB-1}) equals to the quotient $R(\xi_{\matW})$
for the choice $\matM_{\matW}\coloneqq\ldaB$, since \cite[Proposition 5.5.6 and Eq. (5.25)]{AMS09}
show that at a critical point the term $g_{\matW}(\xi_{\matW},\hess{\flda[\xi_{\matW}]})$
which is the numerator of $R(\xi_{\matW})$ do not depend on the choice
of Riemannian metric. Thus, the optimal values of Eq. (\ref{eq:ldarayleighforMeqB-1})
satisfying Eq. (\ref{eq:LDAopt1a-1}) or Eq. (\ref{eq:LDAopt2a-1})
are the eigenvalues of the Riemannian Hessian at $\matW\in\Stiefellda$
with the choice $\matM_{\matW}\coloneqq\ldaB$. Obviously, classification of the critical
points does not depend on the Riemannian metric. Therefore, we can
classify the critical points using the signs of the eigenvalues of
the Riemannian Hessian at $\matW\in\Stiefellda$ with the choice $\matM_{\matW}\coloneqq\ldaB$.

Recall from the proof of Theorem \ref{thm:LDAHessiantheorem} that in the special case $\matM_{\matW}\coloneqq\ldaB$, the eigenvalues of Eq. (\ref{eq:LDAHessianoptimal-1}  are also the eigenvalues of the Riemannian Hessian at a critical point $\matW\in\Stiefellda$. The eigenvalues are:

\[
\frac{1}{2}\left(\mu_{j}-\mu_{i}\right)\left(\alpha_{j}-\alpha_{i}\right),\ 1\leq j<i\leq p,
\]
and
\[
\mu_{j}\left(\alpha_{j}-\alpha_{i}\right),\ j=1,...,p,\ i=p+1,...,d.
\]

Now, we can conclude the signs of the eigenvalues of the Riemannian
Hessian at any critical point $\matW\in\Stiefellda$ for the choice
$\matM_{\matW}\coloneqq\ldaB$, and classify these critical points. Recall that $\mu_{i}>0$
for $i=1,...,p$, and $\mu_{j}-\mu_{i}>0$ for $j<i$. Also $\alpha_{i}\geq0$.
Thus, the signs of the eigenvalues of Eq. (\ref{eq:LDAHessianoptimal-1})
are only determined by the differences between $\alpha_{j}$, $j=1,...,p$
which are some generalized eigenvalues of the matrix pencil $(\mat{S_{B}},\ldaB)$,
corresponding to the generalized eigenvectors on ordered $1\leq j<i\leq p$,
but on the other hand there exists at least one pair of indexes $j=1,...,p$
and $i=p+1,...,d$ such that $\alpha_{j}-\alpha_{i}<0$, otherwise
it contradicts the assumption on the columns of $\matW$, and the
$d-j$ trailing $\alpha_{i}$'s, i.e., $\alpha_{j}-\alpha_{i}$ where
$1\leq j<i\leq p$ or $j=1,...,p$ and $i=p+1,...,d$. 

In particular, for $\matW^{\star}$ such
that the columns are the $p$-dominant generalized eigenvectors of
the matrix pencil $(\mat{S_{B}},\ldaB)$, then by the assumption that
$\alpha_{1}=\rho_{1}>...>\alpha_{p}=\rho_{p}$, we have $\alpha_{j}-\alpha_{i}>0$
where $1\leq j<i\leq p$ or $j=1,...,p$ and $i=p+1,...,d$. Therefore,
in this case all the eigenvalues of Eq. (\ref{eq:LDAHessianoptimal-1})
are strictly positive, thus the matrix in Eq. (\ref{eq:LDAHessianoptimal-1})
is SPD, and consequently the eigenvalues of the Riemannian Hessian
at $\matW^{\star}$ are all strictly positive. Thus,
$\matW^{\star}$ is a strict local minimum of $\flda$ on $\Stiefellda$. 

If we prove that $\matW^{\star}$ it is the only local minimum (up
to the signs of the columns), then it is the only asymptotically stable
critical point according to Theorem \ref{thm:stabilityLDA}. In order
to prove it we further assume that for all $i=1,...,d$ the values
$\rho_{i}$ are distinct, then we can conclude the following. Suppose
$\matW$ is any other critical point differs from $\matW^{\star}$
at the optimal value, i.e., such that the columns of $\matW$ are
ordered generalized eigenvectors corresponding to some $p$ singular
values of the matrix pencil $(\mat{S_{B}},\ldaB)$, $\alpha_{1},...,\alpha_{p}$,
which are not the leading $p$ values. We consider the different cases:
\begin{enumerate}
\item Suppose $\alpha_{1},...,\alpha_{p}$ are not ordered in any particular
order (possible only for $p\geq3$), then there exists $j$ such that
$\alpha_{j}$ is larger than some $\alpha_{k}$ and smaller than $\alpha_{m}$
where $j<k,m\leq p$, then $\alpha_{j}-\alpha_{k}>0$ and $\alpha_{j}-\alpha_{m}<0$.
In this case there are both strictly positive and strictly negative
eigenvalues of the Riemannian Hessian at $\matW$ for the choice $\matM_{\matW}\coloneqq\ldaB$,
thus $\matW$ is a saddle point.
\item Suppose $\alpha_{1},...,\alpha_{p}$ are ordered in a descending order.
Since $\matZ$ is not an optimal solution of Problem~(\ref{eq:lda_riemannian}),
then there exists at least one $1\leq j\leq p$ for which $\alpha_{j}\neq\rho_{j}$.
Thus, on the one hand $\alpha_{j}-\alpha_{i}>0$ where $1\leq j<i\leq p$,
but on the other hand there exists at least one pair of indexes $j=1,...,p$
and $i=p+1,...,d$ such that $\alpha_{j}-\alpha_{i}<0$, otherwise
it contradicts the assumption that there exists at least one $1\leq j\leq p$
for which $\alpha_{j}\neq\rho_{j}$. In this case there are both strictly
positive and strictly negative eigenvalues of the Riemannian Hessian
at $\matW$ for the choice $\matM_{\matW}\coloneqq\ldaB$, therefore, $\matW$ is
a saddle point.
\item Suppose $\alpha_{1},...,\alpha_{p}$ are ordered in an ascending order.
Then, $\alpha_{j}-\alpha_{i}<0$ where $1\leq j<i\leq p$. Now, we
consider two sub-cases:
\begin{enumerate}
\item There exists at least one $1\leq j\leq p$ for which $\alpha_{j}\neq\rho_{d-j+1}$.
Then, there exists at least one pair of indexes $j=1,...,p$ and $i=p+1,...,d$
such that $\alpha_{j}-\alpha_{i}>0$, otherwise it contradicts the
assumption that there exists at least one $1\leq j\leq p$ for which
$\alpha_{j}\neq\rho_{d-j+1}$. In this case there are both strictly
positive and strictly negative eigenvalues of the Riemannian Hessian
at $\matW$ for the choice $\matM_{\matW}\coloneqq\ldaB$, therefore, $\matW$ is
a saddle point.
\item Consider the case $\alpha_{1}=\rho_{d}<\alpha_{2}=\rho_{d-1}<...<\alpha_{p}=\rho_{d-p+1}$.
Thus, $\alpha_{j}-\alpha_{i}<0$ where $1\leq j<i\leq p$ or $j=1,...,p$
and $i=p+1,...,d$. In this case all the eigenvalues of the Riemannian
Hessian at $\matW$ are strictly negative for the choice $\matM_{\matW}\coloneqq\ldaB$,
thus, $\matW$ is a local maximizer. Since it is the only strict local
maximizer up to the signs of the columns of $\matW$, and it is also
the global maximizer.
\end{enumerate}
\end{enumerate}
In all cases, $\matW$ is not a local minimizer of $\flda$ on $\Stiefellda$. Thus, $\matW^{\star}$ is the only local minimum
(up to the signs of the columns.
According to Theorem \ref{thm:stabilityLDA} all these critical points
are unstable.
\end{proof}

\subsection{Proof of Corollary~\ref{cor:LDA}}\label{subsec:proofofcorLDA}
\begin{proof}
The condition number bound follows immediately from Lemma \ref{lem:sketching}
and Theorem \ref{thm:LDAHessiantheorem}. For the costs, the same
arguments as in the proof of Corollary \ref{cor:CCA} hold: none of
the operations require forming $\matSb$ or $\matSw$, but instead
require taking product of these matrices with vectors. These products
can be computed in cost proportional to the iterated products of the
matrices $\hat{\matX}$ and/or $\hat{\matY}$ with vectors. Computing
a product of the matrix $\matXhat=\matX-\matY$ with a vector is equivalent
to computing the product of the same vector with the matrices $\matX$
and $\matY$ and subtracting the result. Computing the product of
$\matX$ with a vector is proportional number of non-zeros in $\mat X$,
and the cost of the product of $\matY$ with a vector is $O(ld)$
since $\matY$ has exactly $l$ distinct rows. Computing a product
of the matrix $\matYhat$ with a vector costs $O(ld)$ since $\matYhat$
is a $l\times d$ matrix. The preprocessing steps follow from Table
\ref{tab:costs}. Assuming a bounded number of line-search steps in
each iteration of Riemannian CG, the costs follows from Table \ref{tab:costs},
as each iteration requires a bounded number of computations of each
of the following: objective function evaluation $O\left(p\left(\nnz{\matX}+ld\right)\right)$,
retraction $O\left(\nnz{\hat{\matX}}p+dp^{2}\right)$, vector transport
and Riemannian gradient $O\left(p\left(\nnz{\matX}+ld\right)+\nnz{\hat{\matX}}p+dp^{2}+d^{2}p\right)$. 
\end{proof}

\bibliography{ellipsoid}
\bibliographystyle{plain}

\end{document}